\documentclass{amsart}

\usepackage{amsmath}
\usepackage{amssymb}
\usepackage{amsfonts}
\usepackage{MnSymbol}
\usepackage{graphics} %% add this and next lines if pictures should be in esp format
\usepackage{epsfig} %For pictures: screened artwork should be set up with an 85 or 100 line screen
\usepackage{psfrag}
\usepackage{pdfsync}
\usepackage[usenames]{color}
\usepackage{cancel}

\definecolor{arancio}{rgb}{0.90,0.50,0.20}

\definecolor{blu}{rgb}{0.,0.,1.}
\newcommand{\blu}[1]{{\color{blu}{#1}}}

\definecolor{pavone}{rgb}{0.00,0.00,0.63}

\definecolor{malva}{rgb}{0.10,0.50,0.50}

\definecolor{rosso}{rgb}{1.,0.,0.}

\definecolor{geranio}{rgb}{0.90,0.00,0.20}

\definecolor{cerulean}{rgb}
{0.0, 0.48, 0.65}

\newtheorem{theorem}{Theorem}[section]

\newtheorem{lem}[theorem]{Lemma}
\newtheorem{prop}[theorem]{Proposition}

\theoremstyle{definition}
\newtheorem{definition}{Definition}[section]
\newtheorem{remark}{Remark}[section]

\newcommand{\ep}{\varepsilon}
\newcommand{\N}{\mathbb{N}}
\newcommand{\R}{\mathbb{R}}
\newcommand{\lla}{\left\langle}
\newcommand{\rra}{\right\rangle}
\date{\today}

\newcommand{\bcl}{\begin{center}}
\newcommand{\ecl}{\end{center}}
\newcommand{\brl}{\begin{right}}
\newcommand{\erl}{\end{right}}
\newcommand{\ben}{\begin{enumerate}}
\newcommand{\een}{\end{enumerate}}
\newcommand{\barr}{\begin{array}}
\newcommand{\earr}{\end{array}}
\newcommand{\btab}{\begin{tabular}}
\newcommand{\etab}{\end{tabular}}
\newcommand{\bdoc}{\begin{document}}
\newcommand{\edoc}{\end{document}}
\newcommand{\beqy}{\begin{eqnarray}}
\newcommand{\eeqy}{\end{eqnarray}}

\newcommand{\beq}{\begin{equation}}
\newcommand{\eeq}{\end{equation}}
\newcommand{\beqi}{\begin{eqnarray*}}
\newcommand{\eeqi}{\end{eqnarray*}}
\newcommand{\bitem}{\begin{itemize}}
\newcommand{\brem}{\begin{remark}}
\newcommand{\erem}{\end{remark}}
\newcommand{\eitem}{\end{itemize}}
\newcommand{\nln}{\newline}
\newcommand{\newt}{\newtheorem}

\renewcommand{\a }{\alpha }
\renewcommand{\b }{\beta }
\newcommand{\g }{\gamma}
\newcommand{\G }{\Gamma }
\renewcommand{\d }{\delta }

\newcommand{\D }{\Delta }
\newcommand{\e }{\epsilon }
\newcommand{\z }{\zeta }
\renewcommand{\l }{\lambda }
\renewcommand{\L }{\Lambda }
\newcommand{\m }{\mu }
\newcommand{\n }{\tau }
\renewcommand{\r }{\rho }
\newcommand{\s }{\sigma }
\newcommand{\Sig }{\Sigma }
\renewcommand{\t }{\tau }
\newcommand{\var }{\varphi }
\renewcommand{\o }{\omega }
\renewcommand{\O }{\Omega }

\newcommand{\supp}{\text{\rm supp}\,}
\newcommand{\sgn}{\text{\rm sgn}\,}
\title[Measure-valued solutions]
{Radon measure-valued solutions 
\\  of first order hyperbolic conservation laws }

\author[Bertsch]{Michiel Bertsch}
\address{Dipartimento di Matematica, Universit\`a di Roma "Tor Vergata", 
Via della Ricerca Scientifica, 00133 Roma, Italy \\ and
Istituto per le Applicazioni del Calcolo "M. Picone", CNR, Roma, Italy} 
\email{bertsch.michiel@gmail.com}
\author[Smarrazzo]{Flavia Smarrazzo}
\address{Universit\`a Campus Bio-Medico di Roma\\ Via Alvaro del Portillo 21, 00128 Roma, Italy}
\thanks{}
\email{flavia.smarrazzo@gmail.com}
\author[Terracina]{Andrea Terracina}
\address{Dipartimento di Matematica "G. Castelnuovo", Universit\`a ''Sapienza'' di Roma\\ P.le A. Moro 5, I-00185 Roma, Italy}
\email{terracina@mat.uniroma1.it}
\author[Tesei]{Alberto Tesei}
\address{Dipartimento di Matematica "G. Castelnuovo", Universit\`a ''Sapienza'' di Roma\\ P.le A. Moro 5, I-00185 Roma, Italy, and 
Istituto per le Applicazioni del Calcolo "M. Picone", CNR, Roma, Italy}
\email{albertotesei@gmail.com}

\subjclass{Primary: 35D99, 35K55, 35R25; Secondary: 28A33, 28A50.
 }  

\keywords{ First order hyperbolic conservation laws,
Radon measure-valued solutions, entropy inequalities, uniqueness.}

\date{\today}

\begin{document}

\bibliographystyle{h-elsevier2}
\begin{abstract} 

We study nonnegative solutions of the Cauchy problem 
$$
\left\{\begin{array}{ll}
u_t+ \left[\varphi(u)\right]_x=0 
& \quad\mbox{in}\  \R\times (0,T) \\
u=u_0\ge 0 &\quad\mbox{in}\  \R\times \{0\}\,,
\end{array}\right. \leqno{(P)}
$$ 
where $u_0$ is a Radon measure and $\varphi:[0,\infty)\mapsto \R$ is a  globally Lipschitz continuous function. We construct suitably defined entropy solutions in the space of Radon measures. Under some additional conditions on $\varphi$, we prove their uniqueness if the singular part of $u_0$ is a finite superposition of Dirac masses. In terms of the behaviour of $\varphi$ at infinity we give criteria to distinguish two cases: either all solutions are function-valued for positive times (an instantaneous  regularizing effect), or the singular parts of certain solutions persist until some positive {\em waiting time} (in the linear case $\varphi(u)=u$ this happens for all times). In the latter case we describe the evolution of the singular parts.
\end{abstract}

\maketitle

%%%%%%%%%%%%%%%%%%%%%%%%%%%%%%%%%%%%%%%%%%%%%%%%%%%%%%%%%%%%%%%%%%%%%%%%%%%%%%%%%%%%%%%%%%%%%%%%%%%%%%%%%%%%%%%%%%

\section{Introduction}\label{intro}

%\fl{6) Ho cambiato la definizione di $\varphi_{\ep}$ in modo da avere $\varphi_{\ep}(0)=0$ [di questo l'ultima volta non abbiamo parlato]} \rosso{Non capisco bene il motivo. Comunque mi sembrerebbe pi\`u logico scrivere la definizione come ho fatto io.}

In this paper we consider the Cauchy problem 
$$
\left\{\begin{array}{ll}
u_t+\left[\varphi(u)\right]_x=0 
& \quad\mbox{in}\  \R\times (0,T)=:S \smallskip\\
u=u_0  &\quad\mbox{in}\  \R\times \{0\}\,,
\end{array}\right. \leqno{(P)}
$$ 
where 
$T>0$,   $u_0$ is a nonnegative finite Radon measure on $\R$, and 
$\varphi:[0,\infty)\mapsto\R$, $\varphi(0)=0$, is a Lipschitz continuous function (see assumption $(H_1)$). Therefore, $\varphi$  grows at most linearly. 

\smallskip

\noindent $(a)$ Problem $(P)$ with a {\em superlinear} $\varphi$ of the type $\varphi(u)=u^p$, $p>1$ was studied in \cite{LP}, proving  existence and uniqueness of nonnegative entropy solutions (see also \cite{CP}). %; to our knowledge,  \cite{LP} is the only paper studying $(P)$ with measure-valued initial data and nonlinear $\varphi$). 
By definition, in that paper the solution for positive times takes values  in $L^1(\R)$, although the initial data $u_0$ is a finite Radon measure. Interesting, albeit sparse results concerning $(P)$ with %Radon measure-valued initial data and  
$\varphi$ {\em at most linear at infinity} can be found in the pioneering paper \cite{DeS}, in which the same definition of 
%\blu{(weak)}
 Radon measure-valued solutions used below (see  equality \eqref{ewf}) was proposed.

When $\varphi(u)=Cu$ $(C\in\R)$ problem $(P)$ is the Cauchy problem for the {\em linear transport equation},
$$
\left\{\begin{array}{ll}
u_t+Cu_x=0 
& \quad\mbox{in}\  S \smallskip\\
u=u_0  &\quad\mbox{in}\  \R\times \{0\}\,,
\end{array}\right. 
$$ 
whose solution is trivially the translated of $u_0$ along the lines $x=Ct+x_0$ $(x_0\in\R)$.   
In particular, the singular part $u_s(\cdot,t)$ of the solution is nonzero for $t>0$ if and only if  the same holds for $t=0$.

It is natural to ask what happens if $\varphi$ is {\em sublinear}. To address this case we must consider {\em solutions of problem $(P)$ which for $t>0$ possibly are finite Radon measures on $\R$ as the initial data $u_0$}. Therefore, throughout the paper we consider solutions of problem $(P)$ as maps from $[0,T]$ to the cone of nonnegative finite Radon measures on $\R$, which satisfy $(P)$ in the following sense: for a suitable class of test functions $\zeta$ there holds
 \begin{eqnarray*}%\label{provvi}
&&\iint_{S} \big[u_r\zeta_t+ \varphi(u_r) \zeta_x\big]\,dxdt+\int_0^{T}\lla u_s(\cdot,t), 
\zeta_{\nu}(\cdot,t)\rra_{\R}dt= - \lla u_{0},\zeta(\cdot,0)\rra_{\R}
\end{eqnarray*}
 (see Definition \ref{deso}). Here the measure $u(t)$ is defined for a.e.~$t\in(0,T)$, 
$u_r \in L^1(S)$ is the density of its absolutely continuous part, 
$\lla \cdot,\cdot\rra_{\R}$ denotes the duality map, and
$$
\zeta_{\nu}:=\zeta_t+C_{\varphi}\,\zeta_x, \qquad  C_{\varphi}:=\lim_{u\to \infty}\frac{\varphi(u)}{u}.
$$
Measure-valued entropy solutions are defined similarly (see Definition \ref{deso}). 

We use an approximation procedure to construct measure-valued entropy solutions of problem $(P)$
(see Theorem \ref{th.exi1}). In addition, we prove that the singular part $u_s$ of an entropy solution of problem $(P)$ does not increase along the lines $x=x_0+C_{\varphi}t$ (see Proposition \ref{ac}).  In particular, if $C_{\varphi}=0$ the map $t\mapsto u_s(\cdot,t)$ is nonincreasing.

Concerning the case when $\varphi$ is sublinear, the following example is particularly instructive:
\begin{equation}\label{ester}
\left\{\begin{array}{ll}
u_t+ \left[\varphi(u)\right]_x=0 
& \quad\mbox{in $S
$} 
\smallskip\\
u=\delta_0  &\quad\mbox{in}\  \R\times \{0\}\,
\end{array}\right. 
\end{equation}
with $S:=\R\times (0,T)$, $T>1$ and 
\begin{equation}\label{esterbis}
\varphi(u)=\sgn\,p \left[(1+u)^p-1\right] \qquad \text{($p<1$, $p\ne 0$)\,.}
\end{equation}
The function 
in \eqref{esterbis} is increasing and concave, with
$C_\varphi=0$, and belongs to a class for which the constructed entropy solution of  problem \eqref{ester}-\eqref{esterbis} is unique (see Theorem \ref{exiuni}). Hence the following holds:
\begin{prop}\label{wai1}
$(i)$ Let  $p<0$. Let $\xi(t)$ be defined by 
\begin{equation*}%\label{ester1}
\xi' = -\dfrac{(|p|t\xi^{-1})^{\frac{p}{1-p}}-1}{(|p|t\xi^{-1})^{\frac{1}{1-p}}-1} \quad\mbox{in $(1,T)$}, 
\qquad \xi(1)=0 \,.
\end{equation*}
Let
$A\blu{:}=\{(x,t)\in S \, |\, 0<x\le |p|t, 0\le t \le 1 \}\cup \{(x,t)\in S \, |\, \xi(t)\le x\le |p|t, 1< t \le T  \}$,
and
\begin{equation}\label{sol1}
u_s(t):=\max\{1-t,0\}\delta_0\,,  
\quad u_r(x,t):=\left[(|p|tx^{-1})^{\frac{1}{1-p}}-1\right]\,\chi_A(x,t)
\quad ((x,t)\in S).
\end{equation} 
Then $u=u_r+u_s$ is the unique constructed 
entropy solution of problem \eqref{ester}-\eqref{esterbis}.

\smallskip
\noindent $(ii)$ Let  $0<p<1$. Let $\xi(t)$ be defined by 
\begin{equation*}%\label{ester2}
\xi' =\dfrac{(|p|t\xi^{-1})^{\frac{p}{1-p}}-1}{(|p|t\xi^{-1})^{\frac{1}{1-p}}-1}\quad\mbox{in $(0,T)$}, \qquad 
\xi(0)=0 \,. 
\end{equation*}
If $B\blu{:}=\{(x,t)\in S \, |\, \xi(t)\le x\le |p|t, 0< t \le T  \}$, then
\begin{equation}\label{sol2}
u(x,t)=u_r(x,t):=\left[(|p|tx^{-1})^{\frac{1}{1-p}}-1\right]\,\chi_B(x,t)
\qquad ((x,t)\in S)\,
\end{equation} 
is the unique constructed 
entropy solution of problem \eqref{ester}-\eqref{esterbis}.
\end{prop}

Let us define the {\em waiting time} $t_0\in [0,T]$ for solutions $u$ of $(P)$: 
\begin{equation}\label{dewait}
t_0:=\inf\{\tau\in (0,T]\,|\,\ u_s(\cdot,t)=0, \ u_r(\cdot,t)\in L^\infty(\R) \text{ for a.e.~}t\in (\tau,T)\} 
\end{equation}
 (by abuse of language, we call $t_0$ ``waiting time'' even if $t_0=T$). Then by Proposition \ref{wai1}:
 $$
 \text{
 {\em Positive waiting times occur in problem \eqref{ester}-\eqref{esterbis} if and only if  $p<0$.}
 } \leqno(*)
 $$
 More precisely, if $p<0$  the singular part $u_s(\cdot,t)$ persists until the waiting time $t_0=1$ at which it disappears, whereas for $0<p<1$  the singular part vanishes for all $t>0$, thus $t_0=0$ - an instantaneous regularizing effect. Instantaneous regularization also occurs if $p>1$ (see \cite{LP} and Remark \ref{uLP}), whereas, as already remarked, in the linear case $p=1$ there holds $t_0=T$ if $u_{0s}\ne 0$. 
 
Since $\varphi(u)=\sgn\,p \left[(1+u)^p-1\right]$ $(p<1, p\ne0)$ is bounded if and only if $p<0$, and $C_\varphi=0$, statement $(*)$ could be rephrased as follows: 

\begin{prop}
Positive waiting times occur in problem \eqref{ester}  
if and only if the map $u\mapsto\varphi(u)-C_\varphi u$, with $\varphi$ as in \eqref{esterbis},  is bounded in $[0,\infty)$.
\end{prop}
 
The above result is generalized to problem $(P)$ by Theorem \ref{th.regbis}, for functions $\varphi$ which satisfy for $u$ large a condition implying either concavity or convexity 
(see  assumption $(H_2')$  
and Remark \ref{phi_k}). The proof of  Theorem \ref{th.regbis} makes use of estimates of the density $u_r$ of the solution of $(P)$, which are strongly reminiscent of the {\em Aronson-B\'enilan inequality} for the porous medium equation (see Proposition \ref{bc}). The main results on the waiting time and the regularity of solutions of $(P)$ are collected in Subsection \ref{subs33}. The existence and an upper bound, in terms of $\varphi$ and $u_0$, of a waiting time was already pointed out in [10, Proposition 2.1] (see also Theorem 3.8-(ii)).
%Let us observe that the possible existence of a waiting time was already pointed out in \cite{DeS}, and an upper estimate of it depending on $\varphi$ and $u_0$ was provided (\cite[Proposition 2.1]{DeS}; see also Theorem \ref{th.regbis}-$(ii)$).}

\smallskip

Another interesting feature of the solution of \eqref{ester}-\eqref{esterbis} with $p<0$ is that  for $t\in(0,1)$ - $i.e.$, as long as $u_s(\cdot,t)>0$ - there holds
\begin{equation*}%\label{condition uniq}
\lim_{x\to 0^+} u_r(x,t)=\infty\,. 
\end{equation*}
Namely, the regular part $u_r(\cdot,t)$ diverges when approaching from the right the point $x_0=0$ where $u_s(\cdot,t)$ is concentrated. As we shall see below (see \eqref{crusg}-\eqref{crus bisg}), this property can be generalized to entropy solutions of a larger class of problems, characterized by the concavity/convexity property on $\varphi$ mentioned before. In this class a generalized form of this property  will also be used as a uniqueness criterion, provided that  $\varphi(u)-C_\varphi u$ is bounded in $[0,\infty)$ and $u_{0s}$ is a finite superposition of Dirac masses (see Proposition \ref{dsuter} and Theorem \ref{exiuni}). 
In \cite{DeS} it was already observed that Kruzkov's entropy inequalities do not guarantee the uniqueness of solutions (see also Remark \ref{altresolu} below), and the formulation of an additional uniqueness criterion was left as an open problem. This problem is addressed  in a forthcoming paper  where
more general {\em compatibility conditions} are given, which ensure uniqueness also for non-convex or non-concave functions $\varphi$ (see \cite{BSTT}). 
%are given in a forthcoming paper \cite{BSTT}.

\smallskip

Apart from the intrinsic mathematical interest of problem $(P)$, it is worth pointing out its connection with a class of relevant models. Ion etching is a common technique for the fabrication of semiconductor devices,  also relevant in other fields of metallurgy, in which the material to be etched is bombarded with an ion beam (see \cite{F, R1,R2}). Mathematical modelling of the process leads to the Hamilton-Jacobi equation in  one space dimension
$$
\left\{\begin{array}{ll}
U_t+ \varphi(U_x)=0 
& \quad\mbox{in}\  \R\times (0,T) \smallskip\\
U=U_0  &\quad\mbox{in}\  \R\times \{0\}\,,
\end{array}\right. \leqno{(HJ)}
$$ 
where $U=U(x,t)$ denotes the thickness of the material and $\varphi$ is bounded, non-convex and vanishing at infinity. Formal differentiation with respect to $x$ suggests to describe the problem in terms of the unknown $u:=U_x$, which formally solves $(P)$ with $u_0=U_0'$. In this way discontinuous solutions of $(HJ)$ correspond to Radon measure-valued solutions of $(P)$ having a Dirac mass $\delta_{x_0}$ concentrated at any point $x_0$ where $U(\cdot,t)$ is discontinuous $(t\in(0,T))$.
A rigorous justification of the above argument, relating {\em discontinuous viscosity solutions} of $(HJ)$ to {\em Radon measure-valued entropy solutions} of $(P)$, is to our knowledge an open problem (in this connection see \cite{CS, Ev}).

Let us mention that a number of ideas used in the present paper go back to papers dealing with Radon measure-valued solutions 
of quasilinear parabolic problems, also of forward-backward type (in particular, see \cite{BST1, BST2, BST3, OPS, PST, ST1}).
 
The paper is organized as follows. In Section~\ref{preli} we recall several known results used in the sequel and introduce some notation. 
In Section~\ref{mainre} we present the main results of the paper. 
In Section~\ref{sec.ap} we introduce the approximation procedure needed for the construction of solutions. 
Sections \ref{sec.ninf}-\ref{sec.uni} are devoted to the proofs of existence, qualitative properties and uniqueness of solutions. 

%%%%%%%%%%%%%%%%%%%%%%%%%%%%%%%%%%%%%%%%%%%%%%%%%%%%%%%%%%%%%%%%%%%%%%%%%%%%%%%%%%%%%%%%%%%%%%%%%%%%%%%%%%%%

\section{Preliminaries}\label{preli}
\setcounter{equation}{0}

\subsection{Function spaces and Radon measures}
%\label{subs21}

We denote by $\mathcal{M}(\R)$ the Banach space of finite Radon measures on $\R$, with norm $\|\mu \|_{\mathcal M(\R)}:=|\mu|(\R)$. By $\mathcal{M}^+(\R)$ we denote the cone of nonnegative  finite Radon measures; if $\mu_1,\mu_2 \in \mathcal{M}(\R)$, we write $\mu_1\le\mu_2$ if $\mu_2-\mu_1\in\mathcal{M}^+(\R)$.  We denote the convex set of probability measures on  $\R$ by  $\mathcal{Pl1}(\R)\subset\mathcal{M}^+(\R)$: 
$\|\t \|_{\mathcal M(\R)}=\t(\R)=1$ for  $\t \in\mathcal P (\R)$.

We denote by $C_c(\R)$ the space  of continuous real functions with compact support in $\R$. 
The space of the functions of bounded variation in $\R$ is denoted by $BV(\R):= \left\{u \in L^1(\R) \, | \, u' \in \mathcal{M}(\R) \right\}$, where $u'$ is the distributional derivative of $u$. It is endowed with the norm
$\|u\|_{BV(\R)}:= \|u\|_{L^1(\R)}+ \|u'\|_{\mathcal{M}(\R)}$. We say that $u\in BV_{loc}(\R)$ if $u\in BV(\Omega)$ for every open bounded subset $\Omega\subset\R$.

The Lebesgue measure, either on $\R$ or $S:=\R\times(0,T)$, is denoted by  $|\cdot|$. Integration with respect to the Lebesgue measure on $\R$ or on $S$ will be denoted by the usual symbols $dx$, respectively $dxdt$. A Borel set  $E$ is {\em null} if $|E|=0$. The expression "almost everywhere", or shortly "a.e.", means "up to null sets". For every measurable function $f$ defined on $\R$ and $x_0\in\R$, we say that {\rm ess\,}$\lim_{x\to x_0}  f(x) =l\in\R$ if there is a  null set $E^*\subseteq\R$  such that $f(x_n)\to l$ for any sequence $\left\{x_n\right\}\subseteq\R\setminus (E^*\cup\{x_0\})$, $x_n\to x_0$. We set $f^{\pm}:=\max\{\pm f,0\}$ for every measurable function $f$ on $\R$.

We denote the duality map between $\mathcal M(\R)$  and  $C_c(\R)$ by $\lla \mu, \r\rra_{\R}:=\int_{\R}\r\,d\mu$. By abuse of notation, we extend $\lla \mu, \r\rra_{\R}$ to any $\mu$-integrable function $\r$.  A sequence $\{\mu_n\}$ converges strongly to $\mu$ in $\mathcal{M}(\R)$ if  $\|\mu_n-\mu\|_{\mathcal{M}(\R)}\to 0$ as  $n\to\infty$.
A sequence $\{\mu_n\}$ of (possibly not finite) Radon measures on $\R$ converges weakly* to a (possibly not finite) Radon measure $\mu$, $\mu_n\stackrel{*}\rightharpoonup \mu$, if $\left\langle \mu_n,\rho\right\rangle_{\R}\to \left\langle \mu,\rho\right\rangle_{\R}$ for all $\rho \in C_c(\R)$. Similar definitions are used for (possibly not finite) Radon measures on $\Omega\times (0,T)$ with $\Omega\subseteq\R$.

Every $\mu\in\mathcal{M}(\R)$ has a 
unique decomposition $\mu=\mu_{ac}+\mu_s$, with $\mu_{ac} \in\mathcal{M}(\R)$
absolutely continuous and $\mu_s\in\mathcal{M}(\R)$ singular with respect to the Lebesgue  measure. 
We  denote by $\mu_r\in L^1(\R)$ the density of $\mu_{ac}$. Every function $f\in L^1(\R)$ can be identified to a finite 
absolutely continuous Radon measure on $\R$; we shall denote this measure by the same symbol $f$ used for the function. 

The restriction $\mu \, \lefthalfcup E$ of $\mu\in \mathcal{M}(\R)$ to a Borel set $E\subseteq \R$ is defined by 
$(\mu \, \lefthalfcup E) (A):=\mu(E\cap A)$ for any Borel set $A\subseteq\R$. 
Similar notations are used for the spaces of finite  Radon measures $\mathcal{M}(\Omega)$ with $\Omega \subseteq\R$, 
$\mathcal{M}(S)$ and $\mathcal{M}(S\times\R)$, where $S:=\R\times (0,T)$.

 We shall use measures $u\in \mathcal{M}(S)$ which, roughly speaking, admit a parametrization with respect to the time variable:
\begin{definition}\label{dli} 
We denote by $L^{\infty}(0,T;\mathcal{M}^+(\R))$ the set of finite nonnegative Radon measures $u\in \mathcal{M}^+(S)$ such that 
for a.e.~$t\in (0,T)$ there is a measure $u(\cdot,t)\in \mathcal{M}^+(\R)$ with the following properties:

\noindent $(i)$ if $\zeta\in C([0,T];C_c(\R))$ the map $t\mapsto \left\langle u(\cdot,t),\zeta(\cdot,t)\right\rangle_{\R}$ 
belongs to $L^1(0,T)$ and
\begin{equation}\label{eq.disintegrazioneU}
\left\langle u,\zeta\right\rangle_{S}=\int_0^T\left\langle u(\cdot,t),\zeta(\cdot,t)\right\rangle_{\R}\,dt\,;
\end{equation}

\noindent $(ii)$ the map $t\mapsto \|u(\cdot,t)\|_{\mathcal{M}(\R)}$ belongs to $L^{\infty}(0,T)$.
\smallskip

\noindent Accordingly, we set
$$
\|u\|_{L^{\infty}(0,T;\mathcal{M}(\R))}:=\mbox{ess}\sup_{t\in (0,T)}\|u(\cdot,t)\|_{\mathcal{M}(\R)}
\quad \text{for  $u\in L^{\infty}(0,T;\mathcal{M}^+(\R))$.}
$$
\end{definition} 

\begin{remark} The definition implies that for all  $\rho\in C_c(\R)$ the map  $t\mapsto \left\langle u(\cdot,t),\rho\right\rangle_{\R}$ is measurable, thus the map   $u:(0,T) \to \mathcal{M}(\R)$ is weakly* measurable ($e.g.$, see \cite[Section 6.7]{Pe}).  For simplicity we prefer the notation $L^{\infty}(0,T;\mathcal{M}(\R))$ to the more correct one $L^{\infty}_{w*}(0,T;\mathcal{M}(\R))$ which is used in \cite{Pe}.
\end{remark}

If $u\in L^{\infty}(0,T;\mathcal{M}^+(\R))$, also $u_{ac}, u_s\in L^{\infty}(0,T;\mathcal{M}^+(\R))$ and, by  \eqref{eq.disintegrazioneU},  
\begin{equation}\label{disicomp}
\lla u_{ac}\,, \zeta\rra_{S}\,=\iint_{S} u_r \,\zeta\,dxdt \,,
\quad 
\lla u_s, \zeta \rra_{S} \,=\int_0^T\!\! \lla u_s(\cdot,t),\zeta(\cdot,t) \rra_{\R} dt
\end{equation}
if $\zeta\in C([0,T];C_c(\R))$. One easily checks that for a.e.~$t\in (0,T)$
\begin{equation}\label{us(t)=u(t)s}
u_{ac}(\cdot,t)=[u(\cdot,t)]_{ac}\,, \quad u_s(\cdot,t)=[u(\cdot,t)]_s\,, \quad u_r(\cdot,t)=[u(\cdot,t)]_r\,,
\end{equation}
where $[u(\cdot,t)]_r$ denotes the density of the measure $[u(\cdot,t)]_{ac}$: for   $\rho\in C_c(\R)$ 
$$
\lla [u(\cdot,t)]_{ac}, \rho\rra_{\R}  =\int_{\R} [u(\cdot,t)]_r \,\rho \,dx 
=\int_{\R} u_r(\cdot,t)\, \rho \,dx \quad\text{for a.e.~$t\in(0,T)$.} 
$$  
In view of  \eqref{disicomp}-\eqref{us(t)=u(t)s}, we shall always  identify the quantities which appear 
on either side of equalities \eqref{us(t)=u(t)s}. 

For any $\mu\in\mathcal{M}(\R)$ and $a\in\R$, the {\em translated measure} $\mathcal{T}_a(\mu)$ is defined by 
\begin{equation*}
\lla  \mathcal{T}_a(\mu), \rho\rra_{\R}\,:=\,\lla \mu, \rho_{-a}\rra_{\R} 
\end{equation*}
for any $\rho\in C_c(\R)$, where $\rho_{-a}(x):=\rho(x+a)$ $(x\in\R)$. Clearly, $\mathcal{T}_a(\mu)\in\mathcal{M}(\R)$ and
$$
[\mathcal{T}_a(\mu)]_{ac}=\mathcal{T}_a(\mu_{ac})\,,\quad [\mathcal{T}_a(\mu)]_s=\mathcal{T}_a(\mu_s)\,.
$$

%%%%%%%%%%%%%%%%%%%%%%%%%%%%%%%%%%%%%%%%%%%%%%%%%%%%%%%%%

\subsection{Young measures}
\label{subs23}

\noindent  We recall  the following result \cite{Ba}. 
\begin{theorem}\label{th.Ball}
Let $\Omega\subseteq \R^N$ be Lebesgue measurable, let $K\subseteq \R$ be closed, and let $u_n:\Omega\mapsto\R$  
be a sequence of Lebesgue measurable functions such that 
$$\lim_{n\to\infty} \left|\left\{x\in\Omega\,|\,u_n(x)\notin U\right\}\right|=0\,$$
for  any open neighbourhood $U$ of $K$ in $\R$. Then there exist a subsequence 
$\{u_j\}\equiv\{u_{n_j}\}\subseteq\{u_n\}$ and a family $\{\t_x\}$ 
of nonnegative measures on $\R$, depending measurably on $x\in\Omega$, such that:

\noindent $(i)$ $\|\t_x\|_{\mathcal{M}(\R)}:=\int_{\R}d\t_x\leq 1$ for a.e.~$x\in\Omega$;

\noindent $(ii)$ ${\rm supp}\,\t_x \subseteq K$ for a.e.~$x\in\Omega$;

\noindent $(iii)$ for every continuous function $f:\R\mapsto\R$ satisfying $\lim_{|\xi|\to \infty}f(\xi)=0$ there holds
\begin{equation*}%\label{eq.y.1infty}
f(u_j )\stackrel{*}\rightharpoonup f^*\quad\mbox{in }L^{\infty}(\Omega)\,,
\end{equation*}
where for a.e.~$x\in \Omega$
\begin{equation}\label{eq.y.2}
f^*(x):=\left\langle\t_{x},f \right\rangle_{\R}
=\int_{\R}f(\xi)\,d \t_{x}(\xi) \,.
\end{equation}

Suppose further that $\{u_j\}$ satisfies the boundedness condition
\begin{equation}\label{bou.Ball}
\lim_{k\to\infty} \sup_j \left|\left\{x\in\Omega\cap B_R:\,|u_j(x)|\geq k\right\}\right|=0\,
\end{equation}
for every $R>0$, where $B_R:=\left\{x\in \R^N\,|\, |x|<R\right\}$. 

\noindent $(iv)$ $\t_x$ is a probability measure for a.e.~$x\in\Omega$;

\noindent $(v)$ given any measurable subset $A\subseteq\Omega$ there holds
\begin{equation}\label{eq.y.1}
f(u_j)\rightharpoonup f^*\quad\mbox{in }L^{1}(A)\,
\end{equation}
for all continuous functions $f:\R\mapsto\R$ such that $\{f(u_j)\}$ is sequentially weakly compact in $L^1(A)$.  
\end{theorem}

Below we shall always refer to the family $\{\t_x\}$ of probability measures given by the previous theorem as the {\it disintegration of the Young measure} $\tau$ (or briefly Young measure) associated to the sequence $\{u_j\}$. We denote the set of Young measures on $\Omega\times\R$ by $\mathcal{Y}(\Omega;\R)$; in particular, $\mathcal{Y}(S;\R)$ denotes the set of Young measures on $S\times\R$ with $S:=\R\times(0,T)$.

\begin{remark}\label{rem.th.Ball}
$(i)$ The argument used in the proof of Theorem \ref{th.Ball} shows that, under hypothesis \eqref{bou.Ball}, 
the convergence in \eqref{eq.y.1} holds true for 
Carath$\acute{\mbox{e}}$odory functions $f:A\times \R\mapsto \R$ if $\{f(\cdot,u_j)\}$ is sequentially weakly relatively 
compact in $L^1(A)$.
\smallskip

\noindent $(ii)$ Condition \eqref{bou.Ball} is very weak. It is equivalent to the statement that for any $R>0$ there is a continuous nondecreasing function $g_R:[0,\infty)\mapsto \R$, such that 
\begin{equation*}%\label{equi.bou.Ball}
\lim_{\xi  \to \infty}g_R(\xi)=\infty\,,\qquad \sup_j \int_{\Omega\cap B_R} g_R(|u_j(x)|)\,dx <\infty\,.
\end{equation*}
Therefore, Theorem \ref{th.Ball} applies to bounded sequences $\{u_j\}$ in $L^1(\Omega)$ (in which case $g_R(\xi)=\xi$).
\end{remark}

If $\Omega\subset \R^N$ is bounded and $\{u_j\}$ is a bounded but not uniformly integrable sequence  in $L^1(\Omega)$, 
it is possible to extract a uniformly integrable subsequence "by removing sets of small measure". 
This is the content of the following "Biting Lemma" ($e.g.$, see \cite{GMS, V1} and references therein).

\begin{theorem}\label{co22} 
 Let $\{u_n\}$ be a bounded sequence in $L^1(\Omega)$, where $\Omega\subset \R^N$ is a bounded open set. 
 Moreover, let $\{u_j\}\subseteq \{u_n\}$ and $\{\t_x\}$ be the subsequence and the Young measure given in Theorem \ref{th.Ball}. 
 Then there exist a subsequence $\left\{u_k\right\}\equiv \{u_{j_k}\}\subseteq\left\{u_j\right\}$ 
and a decreasing sequence of measurable sets  $E_k \subseteq \Omega$ of Lebesgue measure $|E_k|\to 0$, such that the 
sequence $\left\{u_k\chi_{\Omega\setminus E_k}\right\}$ is uniformly integrable and 
\begin{equation*}%\label{merc.221}
u_k\chi_{\Omega\setminus E_k}\rightharpoonup Z
:=\int_{\R}\xi\,d \tau
(\xi)\  \mbox{in } \,L^1(\Omega) \, ,
\end{equation*}
where $Z\in L^1(\Omega)$ is called the {\em barycenter} of  the disintegration $\{\t_x\}$.
\end{theorem}

%%%%%%%%%%%%%%%%%%%%%%%%%%%%%%%%%%%%%%%%%%%%%%%%%%%%%%%%%
%%%%%%%%%%%%%%%%%%%%%%%%%%%%%%%%%%%%%%%%%%%%%%%%%%%%%%%%%

\section{Main results}\label{mainre}
\setcounter{equation}{0}

Throughout the paper we assume that  $u_0\in\mathcal{M}^+(\R)$. Concerning $\varphi$, we always suppose that
\begin{equation*} 
\varphi\in C([0,\infty))\,,\;\; \varphi(0)=0\,,\;\; \varphi'\in L^{\infty}(0,\infty)\,,\;\; \text{there exists 
\;$\displaystyle{
 \lim_{u\to\infty} 
 \frac{\varphi(u)}{u} 
 =:C_{\varphi}
 \,.} 
$
}
\leqno(H_1)
\end{equation*}
Hence there exists $M>0$ such that 
\begin{equation}\label{disder} 
|\varphi'(u)| \le M\,, \quad  |\varphi(u)| \le M u\quad\text{ for a.e. $u> 0$}\,.
\end{equation}

%%%%%%%%%%%%%%%%%%%%%%%%%%%%%%%%%%%%%%%%%%%%%%%%%%%%%%%%%

\subsection{Definition of solution}
\label{subs30}

In the following definitions we denote by $\zeta_{\nu}:=\zeta_t+C_{\varphi}\,\zeta_x\,$ the derivative of any $\zeta \in C^1(S)$ along the vector $\underline{\tau}\equiv(C_{\varphi},1)$.
\begin{definition}\label{desoY} 
By a {\em solution} of problem $(P)$ {\em in the sense of Young measures} we mean a pair $(u,\tau)$ such that:

\noindent $(i)$ $u\in L^\infty(0,T;\mathcal{M}^+(\R))$, $\tau\in \mathcal{Y}(S;\R)$; 

\noindent $(ii)$ ${\rm supp}\,\tau_{(x,t)}\subseteq [0,\infty)$ for a.e.~$(x,t)\in S$, and
\begin{equation}\label{nabladis}
u_r(x,t) 
=\int_{[0,\infty)} \xi\,d\tau_{(x,t)}(\xi)\,,
\end{equation}
where  
$\tau_{(x,t)}\in\mathcal{Pl1}(\R)$ is the disintegration of $\tau$;

\noindent $(iii)$ for all $\zeta\in C^1([0,T];C^1_c(\R))$ with $\zeta(\cdot,T)=0$ in $\R$ there holds  
\begin{equation}\label{eqy}
\iint_{S} \big[u_r\zeta_t+ \varphi^* \zeta_x\big]\,dxdt+\int_0^{T}\lla u_s(\cdot,t), 
\zeta_{\nu}(\cdot,t)\rra_{\R}dt= - \lla u_{0},\zeta(\cdot,0)\rra_{\R}\,,
\end{equation}
where 
\begin{equation}\label{eq.q*}
\varphi^*(x,t):= 
\int_{[0,\infty)} \varphi(\xi)\,d\tau_{(x,t)}(\xi)\quad\text{for a.e.~$(x,t)\in S$.}
\end{equation}
By an {\em entropy solution} of problem $(P)$ {\em in the sense of Young measures} we mean a 
solution in the sense of Young measures such that
\begin{eqnarray}\label{misey}
&& \iint_S \left[E^*\zeta_t+F^*\zeta_x\right]\,dxdt +C_E\int_0^T\left\langle u_s(\cdot,t),\zeta_t(\cdot,t)\right\rangle_{\R}\,dt \,+\\
&& \qquad +C_F\int_0^T\left\langle u_s(\cdot,t),\zeta_x(\cdot,t)\right\rangle_{\R}\,dt \geq  - \int_{\R} E(u_{0r})\zeta(x,0)\,dx -C_E \left\langle u_{0s}, \zeta(\cdot,0)\right\rangle_{\R}\, . \nonumber  
\end{eqnarray}
 for all $\zeta$ as above, $\zeta\geq0$, and  for every pair $(E,F)$, $E,F:[0,\infty)\mapsto\R$, such that
\begin{equation}\label{ass.EF}
\left\{\begin{array}{ll}
\mbox{$E$ convex,\, $E',F'\in L^{\infty}(0,\infty)$, \, $F'=E'\varphi'$ in $(0,\infty)$,}
\smallskip\\
\mbox{there exist  
$\lim_{u\to \infty}\frac{E(u)}{u}=:C_E$, $\lim_{u\to \infty} \frac{F(u)}{u}=:C_F$\,.}
\end{array}\right.
\end{equation}
In \eqref{misey} for a.e.~$(x,t)\in S$ we set
\begin{equation*}%\label{EF*}
E^*(x,t):= 
\int_{[0,\infty)} E(\xi)\,d\tau_{(x,t)}(\xi)\,, \quad
F^*(x,t):= 
\int_{[0,\infty)} F(\xi)\,d\tau_{(x,t)}(\xi)\,.
\end{equation*}
{\em Entropy subsolutions} (respectively {\em supersolutions}) of problem $(P)$ {\em in the sense of Young measures} are defined by requiring that inequality \eqref{misey} be satisfied for all $\zeta$ and  $(E,F)$ as above, with $E$ nondecreasing (nonincreasing, respectively). 
\end{definition}
Observe that choosing $E(u)=\pm u$ in the entropy inequality \eqref{misey} plainly gives the weak formulation \eqref{eqy}.

\begin{remark}\label{rentro}
\noindent $(i)$ By \eqref{disder}, \eqref{nabladis} and \eqref{eq.q*}, 
\begin{equation}\label{fi*}
|\varphi^*(x,t)|\leq M\int_{[0,\infty)} \xi\,d\tau_{(x,t)}(\xi)= M u_r(x,t)\quad \text{for a.e.~} (x,t)\in S\,. 
\end{equation}
Since $u_r\in L^{\infty}(0,T;L^1(\R))$, by \eqref{fi*} we have that $\varphi^* \in L^{\infty}(0,T;L^1(\R))$, 

\smallskip

\noindent $(ii)$ By \eqref{ass.EF} the functions $E$, $F$ have at most linear growth. Arguing as in $(i)$, it follows that $E^*$ and $F^*$ belong to $L^{\infty}(0,T;L^1_{\rm loc}(\R))$, and to $ L^{\infty}(0,T;L^1(\R))$ if $E(0)=F(0)=0$.  
\end{remark}
\begin{definition}\label{deso}  
A measure $u\in L^\infty(0,T;\mathcal{M}^+(\R))$ is called a 
{\it solution} of problem $(P)$ if 
for all $\zeta\in C^1([0,T];C^1_c(\R))$, $\zeta(\cdot,T)=0$ in $\R$ there holds
\begin{equation}\label{ewf}
\iint_{S} \big[u_r\zeta_t+ \varphi(u_r) \zeta_x\big]\,dxdt+\int_0^{T}\lla u_s(\cdot,t), 
\zeta_{\nu}(\cdot,t)\rra_{\R}dt= - \lla u_{0},\zeta(\cdot,0)\rra_{\R}\,,
\end{equation}
A  solution of problem $(P)$ is called an {\em entropy solution}, if for all $\zeta\geq0$ as above and for all $(E,F)$ as in \eqref{ass.EF} it 
satisfies the {\em entropy inequality} 
\begin{eqnarray}\label{misei}
&& \iint_S \left[E(u_r)\zeta_t+F(u_r)\zeta_x\right]\,dxdt +C_E\int_0^T\left\langle u_s(\cdot,t),\zeta_t(\cdot,t)\right\rangle_{\R}\,dt +\\
&& \quad + C_F\int_0^T\!\langle u_s(\cdot,t),\zeta_x(\cdot,t)\rangle_{\R}\,dt 
\geq - \!\int_{\R} \! E(u_{0r})\zeta(x,0)\,dx -C_E \left\langle u_{0s}, \zeta(\cdot,0)\right\rangle_{\R}.\nonumber  
\end{eqnarray}
{\em Entropy subsolutions} (respectively {\em supersolutions}) of problem $(P)$ are defined by requiring  \eqref{misei} to be satisfied for all $\zeta$ and  $(E,F)$ as before, with $E$ nondecreasing (nonincreasing, respectively). 
\end{definition}
 
A solution of problem $(P)$ is also a solution  in the sense of Young measures. 
Moreover, it follows from \eqref{disder} that $\varphi(u_r)\in L^{\infty}(0,T;L^1(\R))$. 
Similar remarks hold for  entropy solutions, subsolutions and supersolutions.

\begin{remark}\label{rem.su}
\noindent $(i)$ If $C_{\varphi}=0$, equality \eqref{ewf} reads
\begin{equation*}
\iint_{S} \big[u\zeta_t+ \varphi(u_r) \zeta_x\big]\,dxdt 
= - \lla u_{0},\zeta(\cdot,0)\rra_{\R}\,,
\end{equation*}
whence 
$u_t =-[\varphi(u_r)]_x$ in $\mathcal{D}'(S)$.

\noindent $(ii)$ For the  Kru\v zkov entropies $E(u)=|u-k|$, $F(u)=\sgn(u-k)\left [\varphi(u)-\varphi(k)\right ]$ $(k\in[0,\infty))$ there holds $C_E=1$, $C_F=C_{\varphi}$.  
Then inequality \eqref{misei} reads, for all $k\in[0,\infty)$,
\begin{eqnarray}\label{miseikru}
&&\iint_S \left\{|u_r-k|\,\zeta_t+\sgn(u_r-k)\left [\varphi(u_r)-\varphi(k)\right ]\zeta_x\right\}\,dxdt 
+\\
&&\qquad + \int_0^T\left\langle u_s(\cdot,t),\zeta_{\nu}(\cdot,t)\right\rangle_{\R}\,dt \geq - \int_{\R} |u_{0r}-k|\,\zeta(x,0)\,dx - \left\langle u_{0s}, \zeta(\cdot,0)\right\rangle_{\R}.\nonumber  
\end{eqnarray}
\end{remark}

The following proposition states that for any solution of $(P)$  in the sense of Young measures the map $t\mapsto u(t)$, possibly redefined in a null set, is continuous up to $t=0$ with respect to  the weak$^*$ topology of $\mathcal M^+(\R)$. In particular, it explains in which sense the initial condition is satisfied. 
\begin{prop}\label{condini}
Let  $(H_1)$ be satisfied, let 
$(u,\tau)$ be a solution of problem $(P)$  in the sense of Young measures, and let $\rho\in C_c(\R)$. Then
\begin{equation}\label{cini}
\mbox{{\upshape ess}}\lim_{t\to 0^+}\lla u(\cdot,t),\rho\rra_{\R}=\lla u_0,\rho\rra_{\R}\,,
\end{equation}
\begin{equation}\label{cini bis}
\mbox{{\upshape ess}}\lim_{t\to t_0}\lla u(\cdot,t),\rho\rra_{\R}=\lla u(\cdot,t_0),\rho\rra_{\R}\quad\text{for a.e.~$t_0\in (0,T)$}\,.
\end{equation}
The map $t\mapsto u(t)$ has a representative defined for all $t\in [0,T]$, such that 
\begin{equation}\label{cini tris}
\lim_{t\to t_0}\lla u(\cdot,t),\rho\rra_{\R}=\lla u(\cdot,t_0),\rho\rra_{\R}\quad\text{for all $t_0\in [0,T]$}\,.
\end{equation}
\end{prop}

%%%%%%%%%%%%%%%%%%%%%%%%%%%%%%%%%%%%%%%%%%%%%%%%%%%%%%%%%

\subsection{Existence and monotonicity}
\label{subs31}

The existence of solutions is proven by an approximation procedure. 
If $u_0\in\mathcal{M}^+(\R)$, there exist $u_{0n}\in   L^1(\R)\cap L^{\infty}(\R)$  such that 
\begin{equation}\label{stima.u0n}
u_{0n}\geq 0 \quad\text{in }\R\,,\qquad \|u_{0n}\|_{L^1(\R)}\leq \|u_0\|_{\mathcal{M}(\R)}\,,
\end{equation}
\begin{equation}\label{conv.qo.u0n}
u_{0n}\stackrel{*}\rightharpoonup u_0\,,
\qquad u_{0n}\to u_{0r}\quad\mbox{a.e.~in }\R\,\,,\qquad \|u_{0n}-u_{0r}\|_{L^1_{{\rm loc}}(\R\setminus{\rm supp}\,u_{0s})}\to 0
\end{equation}
($e.g.$, see \cite[Lemma 4.1]{PST}). Consider the approximating problem
$$
\left\{\begin{array}{ll}
u_{nt}+ \left[\varphi(u_n)\right]_x=0
& \quad\mbox{in}\  S \smallskip\\
u_n=u_{0n} &\quad\mbox{in}\  \R\times \{0\}\qquad (n\in\N)\,.
\end{array}\right. \leqno{(P_n)}
$$ 
Let us recall the definition of entropy solution of problem ($P_n$) ($e.g.$, see \cite{Da}). 
\begin{definition}\label{dsl1} A function $u_n\in L^{\infty}(0,T;L^1(\R)) \cap L^{\infty}(S)$ is called an {\em entropy solution} of problem $(P_n)$  if for every $\zeta\in C^1([0,T];C^1_c(\R))$, $\zeta(\cdot,T)=0$ in $\R$, 
$\zeta\geq0$ and for any couple $(E,F)$ with $E$ convex, $F'=E'\varphi'$ there holds
\begin{equation}\label{eq.KlE}
 \iint_S \left[E(u_n)\,\zeta_t+F(u_n)\,\zeta_x\right]\,dxdt \ge  - \int_{\R} E(u_{0n})\,\zeta(x,0)\,dx \,.
\end{equation}
\end{definition}
\noindent
Entropy solutions are weak solutions:  if $\zeta\in C^1([0,T];C^1_c(\R))$, $\zeta(\cdot,T)=0$ in $\R$
\begin{equation}\label{eq.Klc}
\iint_S\left[ u_n \zeta_t+\varphi(u_n) \zeta_x\right]\,dxdt +  \int_{\R}u_{0n}\,\zeta(x,0)\,dx=0\,.
\end{equation}
Studying the limiting points of the sequence $\{u_n\}$ we shall prove the following result. 
\begin{theorem}\label{th.exi1}
$(i)$ Let $(H_1)$ be satisfied. Then problem $(P)$ has a solution $u$, which is 
obtained as a limiting point of the sequence $\{u_n\}$ of entropy solutions to problems $(P_n)$. 
In addition, $u$ is an entropy solution of problem $(P)$  in the sense of Young measures. 

\noindent $(ii)$ Let $(H_1)$ and the following assumption:
\begin{equation}\label{exH_2}
\left\{\begin{array}{ll}
\mbox{$\varphi\in C^1([0,\infty))$, and for every $\bar{u}> 0$ there exist $a,b\geq 0$, $a+b>0$}\\
\mbox{such that $\varphi'$ is strictly monotone in $[\bar{u}-a,\bar{u}+b]$}\,,
\end{array}\right. 
\end{equation}
be satisfied. Then $u$ is an entropy solution of problem $(P$). 
\end{theorem} 
Hypothesis \eqref{exH_2} fails if for example $\varphi$ is affine in an interval $(a,b)\subset (0,\infty)$. 
In that case Proposition \ref{idebis2}-$(iii)$, 
which characterizes the limiting Young measure, gives some additional information.

\medskip

The following proposition shows that the singular part of an entropy subsolution of $(P)$ does not increase  along the lines $x=C_{\varphi}t+x_0$.
\begin{prop}\label{ac}
Let $(H_1)$ be satisfied. 

\noindent $(i)$ Let $u$ be an entropy subsolution of problem $(P)$  in the sense of Young measures. Then 
\begin{equation}\label{absco2}
u_s(\cdot,t_2) 
\le \mathcal{T}_{C_{\varphi}(t_2-t_1)}\,(u_s(\cdot,t_1))\;\;\text{in }\mathcal{M^+}(\R)\quad \text{for a.e.~}0\le t_1\le t_2\le T\,.
\end{equation}
In particular, 
\begin{equation}\label{absco22}
u_s(\cdot,t) 
\le \mathcal{T}_{C_{\varphi}t}\,(u_{0s})\;\;\text{in $\mathcal{M^+}(\R)$} \quad \text{for a.e.~$t\in(0,T)$\,,} 
\end{equation}
whence $\|u_s(\cdot,t)\|_{\mathcal{M}(\R)} 
\le \|u_{0s}\|_{\mathcal{M}(\R)}$  for a.e.~$t\in(0,T)$.

\noindent $(ii)$ Let $u$ be a solution of problem $(P)$. Then there is conservation of mass:
\begin{equation*}%\label{comas}
\|u(\cdot,t)\|_{\mathcal{M}(\R)}= \|u_0\|_{\mathcal{M}(\R)}\quad \text{ for a.e.~$t\in(0,T)$}\,.
\end{equation*}
\end{prop}

The linear case $\varphi(u)=u$ shows that equality may hold in \eqref{absco2}. Morover, if $C_{\varphi}=0$, it follows from \eqref{absco2} that the map $t\mapsto u_s(\cdot,t)$ is nonincreasing.

%%%%%%%%%%%%%%%%%%%%%%%%%%%%%%%%%%%%%%%%%%%%%%%%%%%%%%%%%

\subsection{Waiting time and regularity }
\label{subs33}

 It is convenient to distinguish two cases: $C_\varphi=0$ (sublinear growth at infinity) and $C_\varphi\ne0$ (linear growth at infinity), with $C_\varphi$ defined by $(H_1)$.
 
 \subsubsection{Sublinear growth }
Beside $(H_1)$, we will use the following assumption:
\begin{equation*}
\left\{\begin{array}{ll} \text{$\varphi\in C^{\infty}([0,\infty))$, $C_\varphi=0$\,;} 
\smallskip\\
 \text{there exist $H\ge-1$, $K\in\R$ such that 
 }
\smallskip\\
\text{$\varphi''(u)\,[H\varphi(u)+K] \leq -[\varphi'(u)]^2<0$ for all $u\in[0,\infty)$\,.  }
\end{array}\right.\leqno(H_2)
\end{equation*}
By $(H_2)$ the map $u\mapsto \varphi''(u)\,[H\varphi(u)+K] $ is strictly negative and continuous  in $[0,\infty)$, hence two cases are possible: either $(a)$ $H\varphi+K>0 \,$, $\varphi''<0$, or $(b)$ $H\varphi+K<0 \,$, $\varphi''>0$ in $[0,\infty)$. In case $(a)$ there holds 
$\varphi'>0$ in $[0,\infty)$, since $\varphi''<0$ and there exists $\lim_{u\to\infty}\varphi'(u)=C_\varphi=0$. Similarly, in case $(b)$ there holds 
plainly $\varphi'<0$ in $[0,\infty)$. In particular, in both cases $(H_2)$ implies \eqref{exH_2}. Moreover, if also $(H_1)$ holds, thus $\varphi(0)=0$, there holds $H\varphi+K>0 \,$ in $[0,\infty)$ if and only if $K>0$.

\begin{remark}\label{HK} The following examples show that all values of $H\ge -1$ may occur in $(H_2)$:
\begin{equation*}
\begin{cases}
\,\varphi(u)\!=\!\sgn p\,[(1\!+\!u)^p\!-\!1]\;\; (p\!<\!1, p\!\ne \!0) &\Rightarrow\quad H\!=\!\tfrac{p}{1-p}\!\in \!(-1,0)\!\cap\! (0,\infty), \ K\!=\!|H|\,,
\smallskip\\
\, \varphi(u)=1-e^{-\alpha u}\ (\alpha>0) &\Rightarrow\quad H=-1, \ K=1\,,
\smallskip\\
\, \varphi(u)\!=\!\log(1\!+\!u), \text{\; or \;} \varphi(u)\!=\!1\!-\!\frac 1{\log(e+u)} 
&\Rightarrow\quad H=0,\ K=1\,.
\end{cases}
\end{equation*}
\end{remark} 
The following property of constructed entropy solutions plays an important role 
as a uniqueness criterion (see its generalized form given by Proposition \ref{dsuter}, and Theorem \ref{exiuni} below).

\begin{prop}\label{dsu}
Let $(H_1)$-$(H_2)$ be satisfied, and let $ \varphi$ be bounded in $[0,\infty)$. Then every entropy solution $u$ of problem $(P)$ given by Theorem \ref{th.exi1} satisfies for a.e.~$t\in (0,T)$ and all $x_0\in  {\rm supp}\,u_s(\cdot,t)$: 
\begin{equation}\label{crus}
{\rm ess} \!\!\lim_{x\to x_0^+}
u_r(x,t)=\infty\, \quad\text{if }\varphi'>0 \text{ in }[0,\infty)\,,
\end{equation}
\begin{equation*}%\label{crusbis}
{\rm ess} \!\!\lim_{x\to x_0^-}u_r(x,t)=\infty \quad\text{if }\varphi'<0\text{ in }[0,\infty)\,.
\end{equation*}
\end{prop}

\begin{theorem}\label{th.reg} 
\noindent $(i)$ Let $(H_1)$ be satisfied,
let $u_{0s}(\{x_0\})>0$ for some $x_0\in \R$ and let $u$ be a solution of problem $(P)$.
If $ \varphi$ is bounded in $(0,\infty)$, then the waiting time $t_0$ defined by \eqref{dewait} satisfies
\begin{equation}\label{stib}
t_0\,\ge \min\left\{ \,T,\frac{ u_{0s}(\{x_0\}) }{\|\varphi\|_{L^{\infty}(0,\infty)}}\,\right\}>0\,.
\end{equation}

\smallskip

\noindent $(ii)$ 
Let $(H_1)$-$(H_2)$ be satisfied, and let $u$ be the entropy solution of problem $(P)$ given by Theorem \ref{th.exi1}. 
\begin{itemize}
\item[$(a)$] If $ \varphi$ is bounded in $(0,\infty)$,  and  
moreover 
$H>-1$, $|K|<\lim_{u\to\infty} |\varphi (u)|=:\gamma$, then 
\begin{equation}\label{stiu}
t_0\le  \min\left\{T, \frac{(H+1)\,\|u_0\|_{\mathcal{M}(\R)}}{\gamma-|K|}\right\}\,.
\end{equation}
\item[$(b)$]  If $ \varphi$ is unbounded in $(0,\infty)$, then $t_0=0$. 
\end{itemize}
\end{theorem} 

\begin{remark}\label{stsh}
Concerning the estimates in \eqref{stib} and \eqref{stiu}, it is worth considering the case in which $u_0=\delta_0$ and $\varphi(u)=1-(1\!+\!u)^p$, $p<0$. By explicit calculations, in Proposition \ref{ester} we show that in this case the waiting time defined in \eqref{dewait} is $t_0=1$. Hence in this case estimates \eqref{stib}-\eqref{stiu} are sharp, since 
$$
\frac{ \delta_{0}(\{0\}) }{\|\varphi\|_{L^{\infty}(0,\infty)}}=1\quad\mbox{and}\quad\frac{(H+1)\,\|\delta_0\|_{\mathcal{M}(\R)}}{\gamma-|K|}=\frac{(p/(1-p)+1)\,\|\delta_0\|_{\mathcal{M}(\R)}}{1+p/(1-p)}=1\,.
$$
\end{remark}

\begin{remark}\label{phi_k} In part $(ii)$ of Theorem~\ref{th.reg} it is enough to require condition $(H_2)$ for large values of $u$. More precisely (see Remark \ref{phi_k bis}), Theorem~\ref{th.reg}-$(ii)$ remains valid if instead of $(H_2)$ for some $k>0$ there holds:
\begin{equation*}
\text{the function \quad 
$\,\begin{cases} 
\;\varphi_k:[0,\infty)\to \R\\
\;\varphi_k(u):=\varphi(u+k)-\varphi(k)
\end{cases}$
 satisfies }(H_2)\,.
\leqno(H_{2,k})
\end{equation*}

In this connection, observe that the conditions $H>-1$ and $|K|<\lim_{u\to\infty} |\varphi (u)|$  exclude the function $\varphi(u)=1-e^{-u}$. 
The same  conditions 
also exclude the function $\varphi(u)=1-\tfrac 1{\log(e+u)}$, where $K=1=\gamma$. However, in this case we can use hypothesis $(H_{2,k})$ for $k>0$, which is satisfied with $H=0$ and $K=\log^{-2}(e+k)<\gamma_k=\log^{-1}(e+k)$. 
\end{remark}

Let us finally mention the following regularization result. 
\begin{prop}\label{wai2}
Let  $(H_1)$-$(H_2)$ be satisfied, and let $ \varphi$ be bounded in $[0,\infty)$. Then for $a.e.$ $t\in(0,T)$  {\rm supp}\,$ u_s(t)$ is a null set.  
\end{prop} 
\begin{remark}\label{remmi} 
It suffices to prove Proposition \ref{dsu}, Theorem \ref{th.reg} and Proposition \ref{wai2} by assuming $\varphi''<0$ 
in $(H_2)$ (hence, $K>0$ by $(H_2)$ and the assumption $\varphi(0)=0$). Otherwise it is easily seen that, if $u\in L^{\infty}(0,T;\mathcal{M}^+(\R))$ is a solution of  problem $(P)$, the map $\tilde u$ defined by setting 
$$
\left\langle \tilde u, \zeta \right\rangle_S:=\int_0^T\left\langle u(\cdot,t),\zeta(-\,\cdot,t)\right\rangle_{\R}dt 
$$
for every $\zeta \in C([0,T];C_c(\R))$ is a solution 
of the problem
\begin{equation}\label{traP}
\left\{\begin{array}{ll}
\tilde u_t+\left[\tilde\varphi(\tilde u)\right]_x=0 
& \quad\mbox{in}\  S \smallskip\\
\tilde u=\tilde u_0  &\quad\mbox{in}\  \R\times \{0\}\,.
\end{array}\right. 
\end{equation}
Here 
$\left\langle \tilde u_0,\rho\right\rangle_{\R}:=\left\langle u_0,\rho(-\,\cdot)\right\rangle_{\R}$ for all $\rho \in C_c(\R)$, and the function $\tilde\varphi:=-\varphi$ satisfies  $(H_2)$ with $\tilde K:=-K$. The same holds for entropy solutions. 
\end{remark}

%%%%%%%%%%%%%%%%%%%%%%%%%%%%%%%%%%%%%%%%%%%%%%%%%%%%%%%%%

 \subsubsection{Linear growth } 
Let  $\varphi$ satisfy the following assumption:
\begin{equation*}
\left\{\begin{array}{ll} 
 \text{$\varphi\in C^{\infty}([0,\infty))$; there exist $H\ge-1$, $K\in\R$ such that }
\smallskip\\
\text{$\varphi''(u)\,\{H[\varphi(u)-C_\varphi u]+K\} \leq -[\varphi'(u)-C_\varphi ]^2<0$ for all $u\in[0,\infty)$\, }
\end{array}\right.\leqno(H_2')
\end{equation*}
(observe that $(H_2')$ reduces to $(H_2)$ if $C_{\varphi}=0$). If $(H_2')$ holds, the function $\tilde\varphi:=\varphi(u)-C_\varphi u$ satisfies  $(H_2)$ since $C_{\tilde\varphi}=0$.

\begin{remark}\label{trasl}
It is easily seen that, if $u$ is a solution (respectively, an entropy solution) of problem $(P)$, then $v\in L^{\infty}(0,T;\mathcal{M}^+(\Omega))$ defined by 
\begin{equation*}%\label{eq.vtr}
v(\cdot,t)=\mathcal{T}_{-h}(u(\cdot,t))\quad\mbox{in}\ \,\mathcal{M}(\R)
\end{equation*}
for any $h\in\R$ 
is a solution (respectively, an entropy solution) of $(P)$ with $u_0$ replaced by $v_0:=\mathcal{T}_{-h}(u_0)$. 
Similarly, $\tilde u(\cdot,t):=\mathcal{T}_{-C_{\varphi}t}(u(\cdot,t)) $ is a solution (respectively, an entropy solution)  of problem \eqref{traP} with  $\tilde u_0=u_0$ and $\tilde\varphi(u)=\varphi(u)-C_\varphi u$. 
\end{remark}

By Remark \ref{trasl}, the above results for the case $C_\varphi=0$ can be generalized as follows.
\begin{prop}\label{dsuter}
Let $(H_1)$-$(H_2')$ be satisfied, and let $u\mapsto \varphi(u)-C_\varphi u$ be bounded in $(0,\infty)$. Then every entropy solution $u$ of problem $(P)$ given by Theorem \ref{th.exi1} satisfies for a.e.~$t\in (0,T)$ and all $x_0\in  {\rm supp}\,u_s(\cdot,t)$
\begin{equation}\label{crusg}
{\rm ess} \!\!\lim_{x\to x_0^+}
u_r(x+C_\varphi t,t)=\infty\quad\text{if }\varphi'>C_\varphi \text{ in }[0,\infty)\,,
\end{equation}
\begin{equation}\label{crus bisg}
{\rm ess} \!\!\lim_{x\to x_0^-}u_r(x+C_\varphi t,t)=\infty \quad\text{if }\varphi'<C_\varphi \text{ in }[0,\infty).
\end{equation}
\end{prop}
\begin{theorem}\label{th.regbis} 
\noindent $(i)$ Let $(H_1)$ be satisfied,
let $u_{0s}(\{x_0\})>0$ for some $x_0\in \R$ and let $u$ be a solution of problem $(P)$.
If $u\mapsto \varphi(u)-C_\varphi u$ is bounded in $(0,\infty)$, then  
$$
t_0\,\ge\,\min\left\{T, \frac{ u_{0s}(\{x_0\}) }{\|\varphi-C_\varphi u\|_{L^{\infty}(0,\infty)}}\right\}>0\,.
$$
\noindent $(ii)$ 
Let $(H_1)$ and $(H_2')$ be satisfied, and let $u$ be the entropy solution of problem $(P)$ given by Theorem \ref{th.exi1}. 
\begin{itemize}
\item[$(a)$] Let $u\mapsto \varphi(u)-C_\varphi u$ be bounded in $(0,\infty)$. 
If $H>-1$ and $|K|<\lim_{u\to\infty} |\varphi (u)-C_\varphi u|=:\tilde\gamma$, then 
$$
t_0\le  \min\left\{T, \frac{(H+1)\,\|u_0\|_{\mathcal{M}(\R)}}{\tilde\gamma-|K|}\right\}\,.
$$
\item[$(b)$]  Let $u\mapsto \varphi(u)-C_\varphi u$ be unbounded in $(0,\infty)$. Then $t_0=0$. 
\end{itemize}
\end{theorem} 
Again, Theorem~\ref{th.regbis}-$(ii)$ remains valid if  for some $k>0$ the function $\varphi_k$ defined in Remark \ref{phi_k}  satisfies $(H_2')$.
\begin{prop}\label{wai2bis}
Let  $(H_1)$-$(H_2')$ be satisfied, and let $u\mapsto \varphi(u)-C_\varphi u$ be bounded in $(0,\infty)$. Then for $a.e.$ $t\in(0,T)$  {\rm supp}\,$ u_s(t)$ is a null set.  
\end{prop}

%%%%%%%%%%%%%%%%%%%%%%%%%%%%%%%%%%%%%%%%%%%%%%%%%%%%%%%%%

\subsection{Uniqueness}  
\label{subs34}

In connection with equality \eqref{cini} observe that, if $u_{0s}\ne 0$ and the waiting time $t_0$ is equal to 0, the map $t\mapsto u(\cdot,t)$ is not continuous at $t=0$ in the strong topology of $\mathcal M(\R)$
(otherwise we would have $\lim_{t\to0^+}\|u_s(\cdot,t)\|_{\mathcal{M}(\R)}=0=\|u_{0s}\|_{\mathcal{M}(\R)}$, a contradiction).
Instead,
continuity along the lines $x=x_0+C_{\varphi}t$ may occurs if the waiting time $t_0$ is positive:  

\begin{prop}\label{preco} 
Let $(H_1)$ be satisfied. Let $u\mapsto \varphi(u)-C_\varphi u$ be bounded in $(0,\infty)$, and let $u_0$ satisfy 
\begin{equation}\label{assu0}
u_{0s}= \sum_{l=1}^N c_l\delta_{x_l} \quad \text{with $c_l\in[0,\infty)$, 
$l=1,\dots,N$ for some $N\in\N$.}
\end{equation}
\noindent  $(i)$ If condition \eqref{exH_2} holds, every entropy solution $u$ of problem $(P)$ given by Theorem \ref{th.exi1}-$(ii)$ 
satisfies 
\begin{equation}\label{epc.tot}
{\rm ess}\lim_{t\to 0^+}\|\mathcal T_{-C_\varphi t}(u(\cdot,t)) - u_{0}\|_{\mathcal{M}(\R)}=0 \,.
\end{equation} 
$(ii)$ All entropy solutions $u$ of problem $(P)$ satisfy $\mathcal{T}_{-C_{\varphi}t}(u(\cdot,t))\in C((0,T];\mathcal{M}(\R))$. 
\end{prop}

Let us mention that the above statement $(ii)$ holds for any $u_0\in\mathcal{M}^+(\R)$, if $\varphi$ satisfies $(H_1)$-$(H_2')$ (see Proposition \ref{bc}). 

The following uniqueness result will be proven in Section \ref{sec.uni}.

\begin{theorem}\label{uni} 
Let  $(H_1)$  be satisfied and let $u\mapsto \varphi(u)-C_\varphi u$ be bounded and monotonic 
in $(0,\infty)$. Let $u_0$ satisfy \eqref{assu0}. Then there exists at most one entropy solution $u$ of problem $(P)$ which satisfies either \eqref{crusg} or \eqref{crus bisg}, and the condition
 \begin{equation}\label{epc}
{\rm ess}\lim_{t\to 0^+}\|u_r(\cdot,t) - u_{0r}\|_{L^1(\R)}=0 \,.
\end{equation}
 \end{theorem}
By Propositions \ref{dsuter}, \ref{preco} 
and Theorem \ref{uni} we have the following existence and uniqueness result (observe that $(H_2')$ implies \eqref{exH_2}). 

\begin{theorem}\label{exiuni} 
Let $(H_1)$-$(H_2')$ be satisfied, and let $u\mapsto \varphi(u)-C_\varphi u$ be bounded in $(0,\infty)$. Let $u_0$ satisfy \eqref{assu0}. Then there exists a unique entropy solution of problem $(P)$ which satisfies \eqref{crusg}-\eqref{crus bisg}. 
 \end{theorem}

\begin{remark}\label{altresolu} Conditions \eqref{crusg}-\eqref{crus bisg} in Theorem \ref{exiuni} cannot be omitted. In fact, there exist entropy solutions of problem $(P)$ which do not satisfy either \eqref{crusg} or \eqref{crus bisg}, depending on $\varphi$. Therefore, by Proposition \ref{dsuter} they are different from those given by Theorem \ref{th.exi1}, thus uniqueness fails.

For example, let $u_{0s}\neq 0$ and $u_{0r}\in L^1(\R)\cap L^{\infty}(\R)$.
Let $u\in L^{\infty}(0,T;\mathcal{M}^+(\R))$ be defined 
\begin{equation*}
u(\cdot,t):=u_r(\cdot,t)+\mathcal{T}_{C_{\varphi}t}(u_{0s})\quad \text{for a.e.~$t\in (0,T)$,}
\end{equation*}
where $u_r\in C([0,T];L^1(\R))\cap L^{\infty}(S)$ is the unique entropy solution of problem $(P)$ with $u_0$ replaced by $u_{0r}$. 
Since $u(\cdot,0)=u_r(\cdot,0)+u_{0s}=u_{0r}+u_{0s}=u_0$, one 
easily checks that \eqref{ewf}-\eqref{misei} are satisfied, thus $u$ is an entropy solution of $(P)$. 
On the other hand $u_r\in  L^{\infty}(S)$, so $u_r(\cdot,t)\in L^{\infty}(\R)$ for a.e.~$t\in (0,T)$ and \eqref{crusg}-\eqref{crus bisg} fails. 
\end{remark}
 
\begin{remark}\label{uLP}  If $u\mapsto\varphi(u)-C_\varphi u$ is
unbounded and satisfies assumptions 
$(H_1)$-$(H_2')$, by \cite[Theorem 1.1]{LP} and Theorem \ref{th.regbis}, for every $u_0\in\mathcal{M}^+(\R)$ 
there exists a unique entropy solution of problem $(P)$ with waiting time $t_0$ equal to 0.
In fact, every entropy solution $u$ given by Theorem \ref{th.regbis} is a solution according to \cite{LP}. This follows if we show that
\begin{equation}\label{eq.cfr.LP}
u=u_r\in L^{\infty}(\R\times (\tau,T))\quad\mbox{for every $\tau>0$}
\end{equation}
and ${\rm ess}\lim_{t\to 0} u(\cdot,t)=u_0$ narrowly in $\mathcal{M}(\R)$, i.e.~${\rm ess}\lim_{t\to 0} 
\langle u(\cdot,t), \rho\rangle = \langle u_0,\rho\rangle$ for all bounded $\rho\in C(\R)$.
The latter follows from \eqref{cini} and Proposition \ref{ac}-$(ii)$ (see \cite[Proposition 2, p. 38]{GMS}). 

To prove \eqref{eq.cfr.LP} we fix $\tau>0$. By \eqref{dewait} we may assume that $u_r(\cdot,\tau)\in L^{\infty}(\R)$ and $u(\cdot,t)=u_r(\cdot,t)$ for all $t\geq \tau$. By standard approximation arguments, we may substitute  in the entropy inequality \eqref{misei} $E(u)=[s-k_{\tau}]_+$, with $k_{\tau}=\|u_r(\cdot,\tau)\|_{L^{\infty}(\R)}$, and $\zeta(x,t)\equiv\chi_{[\tau,t]}(t)$. Hence $\int_{\R}[u_r(\cdot,t)-k_{\tau}]_+\,dx \leq 0$ for a.e.~$t\geq \tau$ and \eqref{eq.cfr.LP} follows. 
\end{remark}

%%%%%%%%%%%%%%%%%%%%%%%%%%%%%%%%%%%%%%%%%%%%%%%%%%%%%%%%%%%%%%%%%%%%%%%%%%%%%%%%%%%%%%%%%%%%%%%%%%%%%%%%%%%%%%%%%%

\section{Approximating problems}\label{sec.ap}
\setcounter{equation}{0}

\noindent In this section we consider problem $(P_n)$.
Let  $u_{0n}\in L^1(\R)\cap L^{\infty}(\R)$ satisfy \eqref{stima.u0n} and let $\{u_{0n}^{\ep}\}\subseteq C^{\infty}_c(\R)$, $u_{0n}^{\ep}\ge0$ be any sequence such that  
\begin{equation}\label{est.u0ep}
\|u_{0n}^{\ep}\|_{L^1(\R)}\leq \|u_{0n}\|_{L^1(\R)}\leq \|u_0\|_{\mathcal{M}(\R)}\,,\quad \|u_{0n}^{\ep}\|_{L^{\infty}(\R)}\leq \|u_{0n}\|_{L^{\infty}(\R)}\,,
\end{equation}
\begin{equation}\label{conv.u0ep}
u_{0n}^{\ep} \to u_{0n}\quad\mbox{in} \ \,L^1(\R)\,,\quad u_{0n}^{\ep} \stackrel{*}\rightharpoonup u_{0n}\quad\mbox{in} \ \,L^{\infty}(\R)\,.
\end{equation}
Let $\eta\in C_c^{\infty}(\R)$ be a standard mollifier, let $\eta_{\ep}(u):=\frac{1}{\ep}\eta\left(\frac{u}{\ep}\right)$ for $ \ep>0$, and set
$$
\varphi_{\ep}(u):=(\eta_{\ep}*\overline{\varphi})(u) -(\eta_{\ep}*\overline{\varphi})(0)
=\!\!\int_{\R}\!\eta_{\ep}(u-v)\overline{\varphi}(v)dv -\!\!\int_{\R}\!\eta_{\ep}(-v)\overline{\varphi}(v)dv \quad (u\in\R)\,
$$ 
(here $\overline{\varphi}(u)=\varphi(u)$ for $u\geq 0$ and $\overline{\varphi}(u)=0$ for $u<0$).
The regularized problem associated with  $(P_n)$,
$$
\left\{\begin{array}{ll}
u_{nt}^{\ep}+ \left[\varphi_{\ep}(u_n^{\ep})\right]_x=\ep u_{nxx}^{\ep} 
& \quad\mbox{in}\  S \medskip \\
u_n^{\ep}=u_{0n}^{\ep} &\quad\mbox{in}\  \R \times \{0\} 
\end{array}\right. \leqno{(Q_n^{\ep})}
$$ 
(where $\ep>0$, $n\in\N$), has a unique strong solution $u_n^{\ep}\in C([0,T];H^2(\R))\cap L^{\infty}(S) $, $u_{nt}^{\ep}\in L^2(S)$  ({\it e.g.}, see \cite{MNRR}). Some properties of the family $\{u_n^{\ep}\}$ are collected in the following lemmata. Up to minor changes the proof is standard ($e.g.$, see \cite{Da}), thus is omitted.

\begin{lem}\label{L1} 
Let $u_n^{\ep}$ be the solution of problem $(Q_n^{\ep})$. Then for every $n\in\N$ and $\ep>0$ 
\begin{equation}\label{stilinfty}
u_n^{\ep}\ge0\;\;\text{in $S$}, \quad \|u_n^{\ep}\|_{L^{\infty}(S)}
\le  \|u_{0n}\|_{L^{\infty}(\R)}\,, 
\end{equation}
\begin{equation*}%\label{cm}
\int_{\R} u_n^{\ep}(x,t)\,dx =  \int_{\R} u_{0n}^{\ep}(x)\,dx
\qquad (t\in(0,T))\,,
\end{equation*}
\begin{equation}\label{stil1}
\sup_{t\in(0,T)}\|u_n^{\ep}(\cdot,t)\|_{L^1(\R)}
\le \|u_{0n}\|_{L^1(\R)}\le \|u_0\|_{\mathcal{M}(\R)}\,, 
\end{equation}
\begin{equation}\label{sl12}
\sup_{t\in(0,T)}\|u^{\ep}_n(\cdot+h,t)-u^{\ep}_n(\cdot,t)\|_{L^1(\R)}\leq \|u_{0n}^{\ep}(\cdot+h)-u_{0n}^{\ep} \|_{L^1(\R)}\quad \text{for any $h\in\R$\,.}
\end{equation}
\end{lem}
\begin{lem}\label{L1bis} 
 Let $\varphi$ satisfy \eqref{disder}. Then there exists $C>0$, which only depends on $\|u_0\|_{\mathcal{M}(\R)}$,
  such that for all $n\in\N$, $\ep\in(0,1)$ and $p\in(0,1)$ 
\begin{equation}\label{stint}
  \ep \iint_S 
  (1+u_n^{\ep})^{p-2}\,(u_{nx}^{\ep})^2\,dxdt\le \frac{C}{p\,(1-p)} \,.
  \end{equation}
\end{lem}

\begin{proof} Let $U\in C^2([0,\infty))$, 
$U'\ge 0$ in $(0,\infty)$, and set 
\begin{equation}\label{defithetaU}
\Theta_{U,\ep}(u):=\int_0^u U'(s)\,\varphi'_{\ep}(s)\,ds + \theta_U  \qquad(\theta_U\in\R)\,.
\end{equation}
Multiplying the first equation in $(Q_n^{\ep})$ by $U'(u_n^{\ep})$ gives 
\begin{equation}\label{RS}
\left[U(u_n^{\ep})\right]_t+\left[\Theta_{U,\ep}(u_n^{\ep})\right]_x=\ep\, \left[U(u_n^{\ep})\right]_{xx} - \ep \,U''(u_n^{\ep})(u_{nx}^{\ep})^2
 \quad\mbox{in}\; S\,.
\end{equation}
Hence for all $\zeta\in C^1([0,T];C^2_c(\R))$ 
\begin{eqnarray}\label{RRS}
&&\ep \iint_S U''(u_n^{\ep})(u_{nx}^{\ep})^2\zeta\,dxdt +\int_{\R} U(u^{\ep}_n(x,T))\zeta(x,T)\,dx =\\
&&= \int_{\R}U(u^{\ep}_{0n})\zeta(x,0)\,dx +\iint_S \left\{U(u_n^{\ep})\,\zeta_t+\Theta_{U,\ep}(u_n^{\ep})\, \zeta_x+\ep\, U(u_n^{\ep})\,\zeta_{xx} \right\}dxdt\,.\nonumber
\end{eqnarray}
By \eqref{disder} and the definition of the function $\varphi_{\ep}$, for all $u\ge 0$
\begin{equation}\label{disthetaU}
|\Theta_{U,\ep}(u)| \le\int_0^u U'(s)|\varphi'_{\ep}(s)|\,ds  + \left|\theta_U\right|
\le  M\, \left[U(u)-U(0)\right]+\left|\theta_U\right|\,.
\end{equation}
Choose $\theta_U=0$,
$U(u)=(1+u)^p-1$, with $p\in(0,1)$, and 
\begin{equation*}
\zeta =\rho_k:= \chi_{\{|x|\le k\}}+ \rho(\cdot-k) \chi_{\{k\le x  < k+1\}}+ \rho(\cdot+k) \chi_{\{-(k+1)< x  \le -k\}}
\qquad (k\in\N)\,,
\end{equation*}
with any $\rho\in C^2_c((-1,1))$ such that $\rho(0)=1$, $0\le\rho\le1$ and the derivatives $\rho'$, $\rho''$ vanish at $\{0\}$.
Then $0\leq U(u)\le u$ for $u\ge 0$ and, by \eqref{stil1}, \eqref{RRS} and \eqref{disthetaU}, 
$$
\begin{aligned}
& \ep p(1-p) \iint_S (1+u_n^{\ep})^{p-2}(u_{nx}^{\ep})^2\rho_k\,dxdt 
\leq  \int_{\R} u_{0n}^{\ep}(x)\,dx +\\
&\quad + \iint_S \left\{M\,u_n^{\ep}\,|\rho'_k| + \ep u_n^{\ep}\,|\rho_k''|\,\right\}\,dx dt \leq 
\left\{1+ (M+1) T\|\rho\|_{C^2([-1,1])} \right\}\,\|u_0\|_{\mathcal{M}(\R)} =:C\,
\end{aligned}
$$
 for all $\ep\in(0,1)$ and $k\in\N$. Passing to the limit $k\to\infty$ we obtain \eqref{stint}. 
\end{proof} 

\begin{lem}\label{L2} 
Let  $\varphi$ satisfy \eqref{disder} 
and let $U\in C^2([0,\infty))$ be such that
\begin{equation}\label{T0}
|U''(u)|\le K\,(1+u)^{p-2} \; \text{ for all $u\in[0,\infty)$,  for some $K\ge0$ and $p\in(0,1)$\,.} 
\end{equation}
Then there exists $C_p>0$ such that for all $n\in\N$ and $\ep>0$
\begin{equation}\label{Ma1}
  \ep \iint_S 
  |U''(u_n^{\ep})|\,(u_{nx}^{\ep})^2\,dxdt\le C_p \,.
  \end{equation}
If moreover $U'\in L^{\infty}(0,\infty)$,  
the family $\{U_{n,\rho}^{\ep}\}$, where
\begin{equation}\label{defUne}
U_{n,\rho}^{\ep}(t):=\int_{\R}U(u_n^{\ep})(x,t)\rho(x)\, dx \qquad (t\in(0,T))\,
\end{equation}
and $\rho \in C_c^2(\R)$, 
is bounded in $BV(0,T)$.
\end{lem}

\begin{proof} Inequality \eqref{Ma1} follows immediately from \eqref{stint} and \eqref{T0}. 
To prove  that $\{U_{n,\rho}^{\ep}\}$ is bounded in $BV(0,T)$, observe that by \eqref{RS}
\begin{equation}\label{bv111}
\left( U_{n,\rho}^{\ep}\right)'(t) 
=\int_{\R}\left[\Theta_{U,\ep}(u_n^{\ep}) \rho'+\ep\, U(u_n^{\ep})\rho'' - \ep U''(u_n^{\ep})\,(u_{nx}^{\ep})^2\rho\right](x,t)\,dx\,.
\end{equation}
Since $U'\in L^{\infty}(0,\infty)$, there exists $N>0$ such that $|U(u)|\,\le\, N\,(1+u)$ for $u\ge 0$.
 Hence $\left|U(u_n^{\ep})\right|\le N(1+u_n^{\ep})$ and by \eqref{defithetaU}, \eqref{disder} and the  definition of $\varphi_{\ep}$, there holds 
 $$
 |\Theta_{U,\ep}(u_n^{\ep})| \le
 \|\varphi'_{\ep}U'\|_{L^{\infty}((0,\infty))} 
 |u^{\ep}_n|+\left|\theta_U\right|=:\tilde{M} u^{\ep}_n + \left|\theta_U\right|\,.
 $$
Then it follows from \eqref{bv111} that
\begin{eqnarray*}
 \left|\left( U_{n,\rho}^{\ep}\right)'\right|(t) &\le & \|\rho\|_{C^2(\R)} \int_{{\rm supp}\,\rho}\left\{(\tilde{M}+\ep N)u^{\ep}_n(x,t)+\ep N+|\theta_U|) \,dx\right\}+ \\
 &&+\ep\|\rho\|_{L^{\infty}(\R)}  \int_{\R} \left[|U''(u^{\ep}_n)|\,(u^{\ep}_n)_x^2\right](x,t)\,dx \,,
\end{eqnarray*}
and, by \eqref{stil1} and \eqref{Ma1}, there exists a constant $C_{p,\rho}>0$ such that
\begin{equation}\label{bv16}
\int_0^T \left|\left( U_{n,\rho}^{\ep}\right)'\right|(t)\,dt \le \|\rho\|_{C^2(\R)} \left\{ (\tilde{M}+N)\,T\,\|u_0\|_{\mathcal{M}(\R)} \,+
C_{p,\rho} \right\}\,.
\end{equation}
On the other hand, by \eqref{stil1} and since $\left|U(u_n^{\ep})\right|\le N(1+u_n^{\ep})$, there holds
\begin{equation}\label{bv16bis}
\int_0^T \left| U_{n,\rho}^{\ep}\right|(t)\,dt 
\le NT\|\rho\|_{L^{\infty}(\R)}\,\left(\|u_0\|_{\mathcal{M}(\R)}+ |{\rm supp}\,\rho\,|\right) \,,
\end{equation}
whence the result follows.
\end{proof} 

From the above lemmata we get the following convergence results. 

\begin{lem}\label{prop.conv} 
$(i)$ Let $\varphi\in C([0,\infty))$.  Then
there exist a subsequence $\{u_n^{\ep_m}\}\subseteq\{u_n^{\ep}\}$ and $u_n\in L^{\infty}(S)\cap L^{\infty}(0,T;L^1(\R))$ such that 
 as $\ep_m\to 0$
\begin{equation}\label{cLeb}
u_n^{\ep_m} \stackrel{*} \rightharpoonup u_n\quad \mbox{in}\ \,L^{\infty}(S),
\quad u_n^{\ep_m}\to u_n  \text{ and }
\varphi_{\ep_m}(u_n^{\ep_m})\to \varphi(u_n) \ \text{a.e.~in $ S$}\,,
\end{equation}
\begin{equation}\label{concomp}
u_n^{\ep_m}\to u_n \quad \text{in }L^1((-L,L)\times (0,T))\quad\text{for all }L>0. 
\end{equation}
Moreover, 
$u_n\ge0$ a.e.~in $S$, $\|u_n\|_{L^{\infty}(S)}
\le  \|u_{0n}\|_{L^{\infty}(\R)}$, and  
\begin{equation}\label{stil1lim}
\sup_{t\in(0,T)}\|u_n(\cdot,t)\|_{L^1(\R)}
\le \|u_{0n}\|_{L^1(\R)}\le \|u_0\|_{\mathcal{M}(\R)}\,. 
\end{equation}
$(ii)$ Let  $\varphi$ satisfy \eqref{disder}, let $\rho \in C_c^2(\R)$, and let $U\in C^2([0,\infty))$, with $U'\in L^{\infty}(0,\infty)$, satisfy \eqref{T0}. Let $U_{n,\rho}^{\ep_m}$ be defined by \eqref{defUne} and set 
\begin{equation}\label{defUn}
U_{n,\rho}(t):=\int_{\R}U(u_n)(x,t)\rho(x)\, dx \qquad (t\in(0,T))\,.
\end{equation}
Then 
\begin{equation}\label{convUne}
U_{n,\rho}^{\ep_m}\to U_{n,\rho} \quad \text{in $L^1(0,T)$ and a.e.~in $(0,T)$}\,.
\end{equation}
\end{lem}
\begin{proof} By \eqref{stilinfty},  $u_n^{\ep_m} \stackrel{*} \rightharpoonup u_n$  in $L^{\infty}(S)$, where $u_n\in L^{\infty}(S)$,  
$\|u_n\|_{L^{\infty}(S)} \le  \|u_{0n}\|_{L^{\infty}(\R)}$ and $u_n\ge0$ a.e.~in $S$. 
The a.e.-convergence of $u_n^{\ep_m}$ and part $(ii)$ follow from \eqref{concomp}, and since 
$\varphi_{\ep}$ converges uniformly to the continuous function $\varphi$ on compact subsets of $\R$, we also 
obtain the a.e.-convergence of $\varphi_{\ep_m}(u_n^{\ep_m})$.

It remains to prove \eqref{concomp} and \eqref{stil1lim}. We claim that 
that for a.e.~$t\in (0,T)$ 
\begin{equation}\label{convLinfty.ep.sez} 
u^{\ep_m}_n(\cdot,t) \stackrel{*}\rightharpoonup u_n(\cdot,t) \quad\mbox{in}\ \,L^{\infty}(\R)\quad \mbox{as}\ \,\ep_m\to 0\,.
\end{equation} 
Set
$I^{\ep_m}_{n,\rho}(t):=\int_{\R} u^{\ep_m}_n(x,t)\,\rho(x)\,dx$  for $t\in (0,T))$ and let $\rho \in C^2_c(\R)$.
By Lemma \ref{L2}, with $U(u)=u$, the sequence $\left\{I^{\ep_m}_{n,\rho}\right\}$ is bounded in $BV(0,T)$ and
has a subsequence (not relabelled) $\left\{I^{\ep_m}_{n,\rho}\right\}$  such that 
\begin{equation}\label{conv.Iep}
I^{\ep_m}_{n,\rho} \to I_{n,\rho}\quad\mbox{in}\ \,L^1(0,T)\quad \mbox{as}\ \,\ep_m\to 0
\end{equation}
for some $I_{n,\rho}\in BV(0,T)$. Since $u_n^{\ep_m} \stackrel{*} \rightharpoonup u_n$ in $L^{\infty}(S)$, 
\begin{eqnarray*}
&& \lim_{m\to\infty} \int_0^T I^{\ep_m}_{n,\rho}(t)\,dt 
=\iint_S u_n(x,t)\,\rho(x)\,dxdt =\int_0^T \left( \int_{\R} u_n(x,t)\,\rho(x)\,dx\right)\,dt\,,
\end{eqnarray*}
whence 
$I_{n,\rho}= \int_{\R} u_n(x,t)\,\rho(x)\,dx$ for a.e.~$t\in (0,T))$
and the convergence in \eqref{conv.Iep} is satisfied along the whole sequence $\left\{I^{\ep_m}_{n,\rho}\right\}$. 
Hence for all $\rho\in C^2_c(\R)$ there exists a null set $N\subset (0,T)$  such that  
\begin{equation*}
\lim_{\ep_m\to0}\int_{\R} u_n^{\ep_m}(x,t)\,\rho(x)\,dx =\int_{\R} u_n(x,t)\,\rho(x)\,dx \quad \mbox{for all}\ \,t\in (0,T)\setminus N\,.
\end{equation*}
Since $C^2_c(\R)$ is dense in $L^1(\R)$ and $L^1(\R)$ is separable, the choice of the set $N$ can be made independent of $\rho$.
Hence we have proven \eqref{convLinfty.ep.sez}.

By \eqref{conv.u0ep}, \eqref{stil1}, \eqref{sl12}, 
and the Fr$\acute{\mbox{e}}$chet-Kolmogorov Theorem,  
$\{u_n^{\ep_m}(\cdot,t)\}$ is relatively compact in $L^1((-L,L))$  for all $t\in (0,T)$ and $L>0$. 
Hence, by  \eqref{convLinfty.ep.sez}, 
\begin{equation}\label{conv.strong.sez}
u_n^{\ep_m}(\cdot,t) \to u_n(\cdot,t)\quad\mbox{in}\ L^1((-L,L))\ \mbox{as}\ \,\ep_m\to 0\ \text{for $L>0$ and a.e.~$t\in (0,T)$}\,,
\end{equation}
and \eqref{stil1lim} follows from \eqref{stil1}.
Finally \eqref{concomp} follows from \eqref{stil1}, \eqref{conv.strong.sez} and  the Dominated Convergence Theorem.  
\end{proof} 

\begin{prop}\label{sl3} 
Let  $\varphi\in C([0,\infty))$. For all $n\in\N$ problem $(P_n)$ has an  entropy solution $u_n$, which is unique if $\varphi$  is locally Lipschitz continuous. For a.e.~$t\in (0,T)$ there holds
\begin{equation}\label{sl12lim}
\|u_n(\cdot+h,t)-u_n(\cdot,t)\|_{L^1(\R)}\leq \|u_{0n}(\cdot+h)-u_{0n} \|_{L^1(\R)}\quad \text{for any $h\in\R$\,,}
\end{equation}
\begin{equation}\label{cmlim}
\int_{\R} u_n(x,t)\,dx =  \int_{\R} u_{0n}(x)\,dx\,.
\end{equation}
Moreover, given $\rho\in C^2_c(\R)$ and $U\in C^2([0,\infty))$, with $U'\in L^{\infty}(0,\infty)$, satisfying \eqref{T0}, 
the sequence $\{U_{n,\rho}\}$  defined by \eqref{defUn} is bounded in $BV(0,T)$.
\end{prop}

\begin{proof} Let $\zeta$ and $E$ be as in Definition \ref{dsl1}, and $F_{\ep}'=E'\varphi_{\ep}'$. Then
\begin{equation}\label{visenin}
  \iint_{S}\left\{E(u_n^{\ep})\,\left(\zeta_t+\ep \, \zeta_{xx}\right)+F_{\ep}(u_n^{\ep})\,\zeta_x\right\}\,dxdt \, +   \int_{\R}E(u_{0n}^{\ep})\,\zeta(x,0)\,dx\geq0, 
 \end{equation}
where $u_n^{\ep_m}$ is defined by Lemma \ref{prop.conv}.  
By \eqref{stilinfty}, it is not restrictive to assume that $E(u)= |u-k|$, $F_{\ep}(u)=\sgn(u-k)\left [\varphi_{\ep}(u)-\varphi_{\ep}(k)\right ]$\; $(k\in[0,\infty))$.
By \eqref{stilinfty}, 				
\begin{equation*}
 \|\varphi_{\ep_m}(u_n^{\ep_m})\|_{L^{\infty}(S)}
\le  \sup_{|v|\le\|u_{0n}\|_{L^{\infty}(\R)}}\left|\varphi_{\ep_m}(v)\right|
\le  \sup_{|v|\le\|u_{0n}\|_{L^{\infty}(\R)}+1}\left|\varphi(v)\right|\,.
\end{equation*}
Since  $\varphi_{\ep_m}(u_n^{\ep_m})\to \varphi(u_n)$ a.e.~in $ S$ (see \eqref{cLeb}), it follows from 
\eqref{concomp} and the Dominated Convergence Theorem that  
\begin{equation*}
 \iint_{S}F_{\ep_m}(u_n^{\ep_m})\,\zeta_x\,dxdt \, \to    \iint_{S} F(u_n)\,\zeta_x\,dxdt\
 \quad \text{as $\ep_m\to0$}\,.
\end{equation*}
The remaining terms in 
\eqref{visenin} (with $\ep=\ep_m$) are dealt with similarly.
Letting $\ep_m\to0$
we obtain \eqref{eq.KlE}, so $u_n$ is an entropy solution of problem $(P_n)$. 
Its uniqueness follows from Kru\v zkov's Theorem (\cite{Se}).

Inequality \eqref{sl12lim} follows from \eqref{sl12} and \eqref{conv.strong.sez}. Concerning \eqref{cmlim}, 
it follows from \eqref{eq.Klc} that for all $\rho \in C^1_c(\R)$ and a.e.~$t\in (0,T)$
\begin{equation}\label{fd.sez}
\int_{\R} u_n(x,t)\,\rho(x)\,dx -\int_{\R} u_{0n}(x)\,\rho(x)\,dx = \int_0^t\!\!\!\int_{\R} \varphi(u_n)(x,s)\,\rho'(x)\,dxds\,.
\end{equation} 
Let $\{\rho_k\}\subseteq C^1_c(\R)$ be such that 
$\rho_k(x)=1$ for $x\in [-k,k]$, $\rho_k(x)=0$ if $|x|\geq k+1$, and $\|\rho'_k\|_{L^{\infty}(\R)}\leq 2$.  
Setting $\rho=\rho_k$ in \eqref{fd.sez} and letting
$k\to \infty$ we get
\begin{eqnarray*}
\left|  
\int_0^t\!\!\!\int_{\R} \varphi(u_n)(x,s)\,\rho'_k(x)\,dxds\,\right| 
\leq 
2M \int_0^t\!\!\!\int_{\{x\in\R\, |\, k\leq |x|\leq k+1\}} |u_n(x,s)|\,dxds\;\to\; 0\,, 
\end{eqnarray*} 
since $u_n\in L^1(S)$. On the other hand, by the Monotone Convergence Theorem, 
$$
\int_{\R} u_n(x,t)\,\rho_k(x)\,dx\;\to\; \int_{\R} u_n(x,t) \,dx\,,
\quad 
\int_{\R} u_{0n}(x)\,\rho_k(x)\,dx \;\to\; \int_{\R} u_{0n}(x)\,dx \,,
$$
and \eqref{cmlim} follows from 
\eqref{fd.sez}.

Finally, let us show that $\{U_{n,\rho}\}$ is bounded in $BV(0,T)$. By \eqref{bv16bis} and \eqref{convUne} 
\begin{eqnarray*}
\int_0^T\left|U_{n,\rho}(t)\right| \,dt = \lim_{\ep_m\to 0} \int_0^T\left|U_{n,\rho}^{\ep_m}(t)\right| \,dt \,\le 
N\|\rho\|_{L^{\infty}(\R)}\,\left(T\,\|u_0\|_{\mathcal{M}(\R)}+ |\,{\rm supp}\,\rho\,|\right),
\end{eqnarray*}
and, 
by \eqref{bv16} and the lower semicontinuity of the total variation 
in $L^1(0,T)$ (\cite[Theorem 1, Subsection 5.2.1]{EG}), we get
$$
\|  U_{n,\rho}'\|_{\mathcal{M}(0,T)}\leq   \|\rho\|_{C^2(\R)} \{ (\tilde{M}+N)\,T\,\|u_0\|_{\mathcal{M}(\R)} \,+
C_{p,\rho} \}\,
$$
with $C_{p,\rho}>0$ as in \eqref{bv16}. This completes the proof. 
\end{proof}

%%%%%%%%%%%%%%%%%%%%%%%%%%%%%%%%%%%%%%%%%%%%%%%%%%%%%%%%%%%%%%%%%%%%%%%%%%%%%%%%%%%%%%%%%%%%%%%%%%%%%%%%%%%%%%%%%%

\section{Existence and monotonicity: Proofs}
\label{sec.ninf}
\setcounter{equation}{0}

We proceed with the proof of Theorem \ref{th.exi1}. 

\begin{prop}\label{th.biting1}
Let $(H_1)$ hold and let $u_n$ be the entropy solution of problem $(P_n)$. Then there exist 
a sequence $\{u_{n_j}\}$ 
and $u\in L^{\infty}(0,T;\mathcal{M}^+(\R))$ such that
\begin{equation}\label{conv.un.t}
u_{n_j}\stackrel{*}\rightharpoonup u \ \ \mbox{in }\mathcal{M}(S).
\end{equation}
For all $L>0$ 
there exists a decreasing sequence $\left\{E_j\right\}\subset (-L,L)\times (0,T)$ of Lebesgue measurable sets with $|E_j|\to 0$ as $j\to\infty$, such that
 \begin{equation}\label{merc.2}
u_{n_j} \chi_{((-L,L)\times (0,T))\setminus E_j}\rightharpoonup u_b
:= \int_{[0,\infty)}\xi\,d \tau
(\xi)\quad\mbox{in}\ \, L^1((-L,L)\times (0,T)) \,,
\end{equation}
where $\tau\in \mathcal{Y}(S;\R)$ is the Young measure associated with $\{u_{n_j}\}$, and
\begin{equation}\label{eq.dom.2}
u_{n_j}\chi_{E_j}\stackrel{*}\rightharpoonup \mu:=u-u_b\quad \mbox{in}\ \, \mathcal{M}((-L,L)\times (0,T)).
\end{equation}
\end{prop} 

\begin{proof} By \eqref{stil1lim}, there exist $u\in \mathcal{M}^+(S)$ and a sequence  $\{u_{n_j}\}$ such that
$u_{n_j}\stackrel{*}\rightharpoonup u$ in $\mathcal{M}(S)$. Arguing as in \cite[Proposition 4.2]{ST1} we obtain that $u\in L^{\infty}(0,T;\mathcal{M}^+(\R))$. 
 
Since by \eqref{stil1lim} the sequence $\{u_{n_j}\}$ is  bounded in $L^1(S)$, by  Theorem \ref{th.Ball} there exist a subsequence of $\{u_{n_j}\}$ (not relabeled) and a Young measure  $\tau\in \mathcal{Y}(S;\R)$ such that:
 
\noindent  $(i)$ for every measurable set $A\subseteq S$, 
\eqref{eq.y.2}-\eqref{eq.y.1} are valid  
for any $f\in C(\mathbb{R})$ such that the sequence $\left\{f (u_{n_j}) \right\} $ is sequentially weakly relatively compact in $L^1(A)$;

\noindent $(ii)$ 
supp$\,\tau_{(x,t)}\subseteq [0,\infty)$ for a.e.~$(x,t)\in S$ (here $\tau_{(x,t)}$ is the disintegration of $\tau$).

\noindent Then the result follows by Theorem \ref{co22} and a standard diagonal procedure. 
\end{proof}

\begin{remark}\label{rem.ub}
The function $u_b$ in \eqref{merc.2} is defined for a.e.~in $(x,t)\in S$, since $\tau$ is globally defined in $S$.
In addition, by \eqref{stil1lim} and the arbitrariness of $L$ in Proposition \ref{th.biting1}, a routine proof shows that $u_b\in L^{\infty}(0,T;L^1(\R))$ and $u_b\geq 0$ a.e.~in $S$. Therefore the Radon measure $\mu\ge 0$ (see \eqref{eq.dom.2}) is defined on $S$, $\mu\in L^{\infty}(0,T;\mathcal{M}^+(\R))$, and 
\begin{equation}\label{eq.stru}
\mu=u-u_b\ \,\Rightarrow\ \,u =  u_b+\mu\quad \hbox{in } \mathcal{M}(S)\,.
\end{equation}
\end{remark}

\begin{prop}\label{th.biting2}
Let $(H_1)$ hold, let $\mu$ be as in 
\eqref{eq.stru} and let $U\in C([0,\infty))$. If 
\begin{equation}\label{limU}
\lim_{u\to\infty}\frac{U(u)}{u}=:C_U\in[0,\infty)\, ,
\end{equation}
for all $L>0$ 
\begin{equation}\label{merc.22}
U(u_{n_j})\, \stackrel{*}\rightharpoonup\, U^* +\, C_U \mu\,
\quad\mbox{in}\ \, \mathcal{M}((-L,L)\times (0,T)) \, ,
\end{equation}
where $U^*\in L^{\infty}(0,T;L^1_{{\rm loc}}(\R))$ is defined by 
\begin{equation*}%\label{eq.U*}
U^*(x,t):= \int_{[0,\infty)}U(\xi)\,d \tau_{(x,t)}(\xi) \quad \text{for a.e.~}(x,t)\in S\,.
\end{equation*}
\end{prop}

\begin{remark}\label{rem.U*}
If $U\in C([0,\infty))$ satisfies \eqref{limU}, there exists $N>0$ such that 
\begin{equation}\label{est.limU}
|U(u)|\leq N(1+u)\quad\mbox{for}\ \,u\geq 0\,.
\end{equation}
Moreover $U^*\in L^{\infty}(0,T;L^1(\R))$ if  $|U(u)|\leq Nu$, since $u_b\in L^{\infty}(0,T;L^1(\R))$ and 
$$
|U^*(x,t)|\!\leq \! \!\int_{[0,\infty)}\! |U(\xi)|d\tau_{(x,t)}\!(\xi)\! \leq\! N\!\!\!\int_{[0,\infty)} \!\xi d\tau_{(x,t)}\!(\xi)\!=\!N\,u_b(x,t) \ 
\text{for a.e.~}(x,t)\!\in\! S.
$$
\end{remark}

\smallskip

\noindent {\em Proof of Proposition \ref{th.biting2}.} 
For all $\ep>0$ there exist $m_\ep>0$ such that 
\begin{equation}\label{assini}
-\ep u< U(u)-C_Uu<\ep u \quad\text{if }u>m_\ep\,.
\end{equation}
For any $m\in\N$, $m > m_\ep$ let $l_{1m},\, l_{2m}\in C([0,\infty))$ be such that $0\le l_{1m}\le 1$, $0\le l_{2m}\le 1$,  $l_{1m}+l_{2m}=1$ in $[0,\infty)$, ${\rm supp}\, l_{1m} \subseteq [0,m+1]$, ${\rm supp}\, l_{2m} \subseteq [m,\infty)$. Then, by \eqref{assini},
\begin{equation}\label{cons}
\left|U(u_{n_j})-\left[U(u_{n_j})\,l_{1m}(u_{n_j}) +C_U u_{n_j}l_{2m}(u_{n_j})\right]\right|
<\ep u_{n_j}l_{2m}(u_{n_j})\quad\text{for $j\in\N$}.
\end{equation}
Since $\sup_{S}\left[|U(u_{n_j})|\,l_{1m}(u_{n_j})\right]\le \sup_{u\in[0,m+1]}|U(u)|<\infty$, 
 $\{U(u_{n_j})\,l_{1m}(u_{n_j})\}$
is uniformly integrable in $(-L,L)\times (0,T)$. Hence, by Theorem \ref{th.Ball}, for all $L>0$
\begin{equation}\label{eq.y.11}
\,U(u_{n_j})\,l_{1m}(u_{n_j})\rightharpoonup \,U_{1m}^*:= \int_{[0,\infty)}U(\xi)\,l_{1m}(\xi)\,d \tau(\xi)
\end{equation}
in $L^1((-L,L)\times (0,T))$. 
Here $U_{1m}^*$ belongs to $L^{\infty}(0,T;L^1_{{\rm loc}}(\R))$ since, by \eqref{est.limU}, 
\begin{equation}\label{stU*}
| U_{1m}^*|\leq  \int_{[0,\infty)}|U(\xi)|\,l_{1m}(\xi)\,d \tau(\xi) \leq N\int_{[0,\infty)} (1+ \xi)\,d\tau(\xi) \leq N(1+u_b).
\end{equation}
Similarly, by \eqref{conv.un.t}, \eqref{merc.2}, \eqref{eq.stru} and \eqref{eq.y.11} with $U(u)=u$, 
\begin{equation}\label{eq.y.12}
\begin{aligned}
&u_{n_j}\,l_{2m}(u_{n_j}) 
= u_{n_j}-u_{n_j}\,l_{1m}(u_{n_j})
\stackrel{*}\rightharpoonup   u- \int_{[0,\infty)}\xi\,l_{1m}(\xi)\,d \tau(\xi)=\\
&\qquad = u_b-\int_{[0,\infty)}\xi\,l_{1m}(\xi)\,d \tau(\xi) +\mu 
=\int_{[0,\infty)}\xi[1- l_{1m}(\xi)]\,d \tau(\xi)+\mu\,= \\
&\qquad =\int_{[0,\infty)}\xi\, l_{2m}(\xi)\,d \tau(\xi)+\mu=:l_{2m}^*+\mu\quad\mbox{in } \mathcal{M}((-L,L)\times (0,T))\,. 
\end{aligned}
\end{equation}
From \eqref{cons}-\eqref{eq.y.12} for any  $\zeta\in C_c((-L,L)\times (0,T))$, $\zeta\ge0$, and $m$ as above we get
\begin{eqnarray}\label{ineq.rip}
&&\iint_{(-L,L)\times (0,T)} \left[U_{1m}^* + (C_U-\ep)\,l_{2m}^*\right]\zeta\,dxdt + (C_U-\ep)\lla \mu,\zeta\rra_{(-L,L)\times (0,T)}\,\le \\
&\le& \liminf\limits_{n_j\to \infty} \iint_{(-L,L)\times (0,T)} U(u_{n_j})\,\zeta\,dxdt \, \le   \limsup\limits_{n_j\to \infty}
\iint_{(-L,L)\times (0,T)} U(u_{n_j})\,\zeta\,dxdt \, \le\nonumber\\
&\le&\iint_{(-L,L)\times (0,T)} \left[U_{1m}^* + (C_U+\ep)\,l_{2m}^*\right]\zeta\,dxdt + (C_U+\ep)\lla \mu,\zeta\rra_{(-L,L)\times (0,T)} \,.\nonumber
\end{eqnarray}
Since $U_{1m}^*\in L^{\infty}(0,T;L^1_{{\rm loc}}(\R))$,  
\begin{eqnarray*}
&& 0\leq l_{2m}^*
\leq \int_{[m,\infty)} \xi \,d\tau
(\xi)\leq u_b \in  L^{\infty}(0,T;L^1(\R))\,, 
\end{eqnarray*}
and 
$$\lim_{\ep_m\to0}l^*_{2m}(x,t) = 0\,,\ \ \lim_{\ep_m\to0}U^*_{1m}(x,t) = U^*(x,t)\quad \mbox{for}\ \,a.e.\ \,(x,t)\in S\,,$$
letting $m\to\infty$ in \eqref{ineq.rip} we get plainly
\begin{eqnarray*}
&&\iint_{(-L,L)\times (0,T)} U^* \zeta\,dxdt + (C_U-\ep)\lla \mu,\zeta\rra_{(-L,L)\times (0,T)} \,\le \\
&\le& \liminf\limits_{n_j\to \infty} \iint_{(-L,L)\times (0,T)} U(u_{n_j})\,\zeta\,dxdt \, \le   \limsup\limits_{n_j\to \infty}
\iint_{(-L,L)\times (0,T)} U(u_{n_j})\,\zeta\,dxdt \, \le\nonumber\\
&\le &\iint_{(-L,L)\times (0,T)} U^*\zeta\,dxdt + (C_U+\ep)\lla \mu,\zeta\rra_{(-L,L)\times (0,T)} \,, \nonumber
\end{eqnarray*}
whence
\begin{eqnarray*}
0&\le& \limsup\limits_{n_j\to \infty}\iint_{(-L,L)\times (0,T)} U(u_{n_j})\,\zeta\,dxdt - \liminf\limits_{n_j\to \infty} \iint_{(-L,L)\times (0,T)} U(u_{n_j})\,\zeta\,dxdt \, \le \\
&\le&
2\ep\lla \mu,\zeta\rra_{(-L,L)\times (0,T)}\,.\nonumber
\end{eqnarray*}
From  
the above inequalities the conclusion follows. 
\hfill$\square$

\begin{prop}\label{corb2}
Let $(H_1)$ hold. Let $\mu$, $U$ and $U^*$  be as in Proposition \ref{th.biting2}. Then
\begin{equation}\label{conv.Uun}
\int_0^T \left| \int_{\R} U(u_{n_j})(x,t)\,\rho(x)\,dx - \int_{\R} U^*(x,t)\,\rho(x)\,dx -C_U\left\langle \mu(\cdot,t),\rho\right\rangle_{\R} \right|\,dt  \to 0
\end{equation}
as $j\to\infty$ for  $\rho\in C_c(\R)$. 
Moreover, for all $L>0$ there exist  a null set $N\subset(0,T)$ and a subsequence of $\{u_{n_j}\}$ (not relabelled), 
such that for all $t\in(0,T)\setminus N$
\begin{equation}\label{merc.222}
U(u_{n_j})(\cdot,t)\, \stackrel{*}\rightharpoonup\, U^*(\cdot,t) +\, C_U \mu(\cdot,t)\quad 
\text{in }\mathcal{M}((-L,L)).
\end{equation}
\end{prop} 

\begin{remark}  
Choosing $U(u)=u$ in \eqref{merc.222}, we obtain that
\begin{equation}\label{conv.un.tt}
u_{n_j}(\cdot,t)\stackrel{*}\rightharpoonup u(\cdot,t) 
\quad\mbox{in}\ \,\mathcal{M}((-L,L))\quad \mbox{for}\ a.e.\ t\in (0,T) \text{ and }L>0.
\end{equation}

If $U\in C([0,\infty))$ 
satisfies \eqref{limU}, $U^*\in L^{\infty}(0,T;L^1_{{\rm loc}}(\R))$ and 
$\{U(u_{n})\}$ is bounded in $L^{\infty}(0,T;L^1_{{\rm loc}}(\R))$ (see \eqref{stil1lim} and \eqref{est.limU}).  
Since every  $\zeta\in C(\R^2)\cap L^{\infty}(\R^2)$ can be uniformly approximated in bounded sets by finite sums  
$\sum_{i=1}^p f^{i,p}(x)g^{i,p}(t)$ with $f^{i,p}$, $g^{i,p}$ bounded and continuous functions  
($1\leq i\leq p\,$; $e.g.$, see \cite[Th\'eor\`eme D.1.1]{D}), it follows from \eqref{conv.Uun} that, as $j\to\infty$, 
for all  $\zeta\in C([0,T];C_c(\R))$
\begin{equation}\label{convUun2}
\int_0^T \left| \int_{\R} [U(u_{n_j})\zeta](x,t)\,dx - \int_{\R} [U^*\zeta](x,t)\,dx -C_U\left\langle \mu(\cdot,t),\zeta(\cdot,t)\right\rangle_{\R} \right|\,dt  \to 0.
\end{equation}
\end{remark}

\smallskip

\noindent {\em Proof of Proposition \ref{corb2}.}  $(i)$ Let us first prove \eqref{conv.Uun} if
$U\in C^2([0,\infty))$, with $U'\in L^{\infty}(0,\infty)$, and 
satisfies \eqref{T0} and \eqref{limU}. 
Let $\rho\in C_c(\R)$, $h \in C_c(0,T)$ and fix any $L>0$ such that ${\rm supp}\,\rho \subset (-L,L)$. Then by \eqref{merc.22}
\begin{equation}\label{usag}
\int_0^TU_{n_j,\rho}(t)h(t)\,dt \,\to 
\int_0^TU^*_{\rho}(t)h(t)\,dt + C_U\int_0^T h(t)\lla \mu(\cdot,t),\rho\rra_{\R}\,dt\,,
\end{equation}
where $U_{n_j,\rho}$ is defined by \eqref{defUn} and
%\begin{equation*}%\label{defU*r}
$U^*_{\rho}(t):=\int_{\R}U^*(x,t)\rho(x)\,dx$. 
%\end{equation*}
Since, by Proposition \ref{sl3},   $\{U_{n_j,\rho}\}$ is bounded in  $BV(0,T)$ if $\rho\in C^2_c(\R)$,  there exists a subsequence
which converges in $L^1(0,T)$. Combined with \eqref{usag} this yields that    
$U_{n_j,\rho} \,\to\, U^*_{\rho} + C_U \lla \mu(\cdot,\cdot),\rho\rra_{\R}$ in $\mathcal{D}(0,T)$ and in $L^1(0,T)$
for all $\rho\in C^2_c(\R)$. 
Since the sequence $\{U(u_{n_j})\}$ is bounded in $L^{\infty}(0,T;L^1((-L,L)))$ and $U^*\in L^{\infty}(0,T;L^1((-L,L)))$,
the condition $\rho\in C^2_c(\R)$ may be relaxed to $\rho\in C_c(\R)$, and we have found \eqref{conv.Uun}.

\smallskip

\noindent $(ii)$ Next we prove \eqref{conv.Uun} for all $U\in C([0,\infty))\cap L^{\infty}((0,\infty))$ (in this case $C_U=0$). We set $U_k: [0,\infty)\mapsto[0,\infty),\, U_k(u):=(U \chi_{[0,k]}*\theta_k)(u)$ for any $u\ge0$, where $\theta_k\ge0$ is a sequence of standard mollifiers $(k\in\N$). Then  $\{U_k\}\subseteq C^2_c([0,\infty))$, $U_k\to U$ uniformly on compact subsets of $[0,\infty)$ and $\|U_k\|_{L^{\infty}(\R)}\leq \|U\|_{L^{\infty}(\R)}$. By  part $(i)$ and 
\eqref{stil1lim}, for all $\rho\in C_c(\R)$ and $k\in\N$, $M>0$ 
$$\begin{aligned}
&\limsup_{j\to\infty}\int_0^T \!\!dt\,\left| \int_{\R} U(u_{n_j})\,\rho(x)\,dx - \int_{\R} U^*(x,t)\,\rho(x)\,dx \,\right| \, \leq \\
&\quad \leq \limsup_{j\to\infty} \iint_{\{0\,\le\, u_{n_j}\leq M\}} |U(u_{n_j})-U_k(u_{n_j})|\,|\rho|\,dxdt +  \\
&\quad +\limsup_{j\to\infty} \iint_{\{u_{n_j}>M\}} |U(u_{n_j})-U_k(u_{n_j})|\,|\rho|\,dxdt +  \iint_{S} |U^*-U_k^*|\,|\rho|\,dxdt\;\le\\
&\quad  \le 
 \|\rho\|_\infty |\supp \rho| \,T\,\|U-U_k\|_{L^{\infty}(0,M)} +  \|\rho\|_\infty \left\{\frac{2\,T}{M }\|u_0\|_{\mathcal{M}(\R)}\,\|U\|_{L^{\infty}(\R)}\right.  +\\
&\quad +\!\!\left.
 \iint_{{\rm supp}\rho \times (0,T)} dxdt \!\int_{\,[0,\infty)} |U_k(\xi)-U(\xi)|\,d\tau_{(x,t)}(\xi) \right\} \\
&\quad \le   2\, \|\rho\|_\infty |\supp \rho|\, T\, \|U-U_k\|_{L^{\infty}(0,M)} +\\
&\quad +  2\, \|\rho\|_\infty \|U\|_{L^{\infty}(\R)} \left\{\frac{\,T\,\|u_0\|_{\mathcal{M}(\R)}}{M} +
 \iint_{{\rm supp}\rho \times (0,T)} dxdt \!\int_{\{\xi>M\}}d\tau_{(x,t)}(\xi) \right\},
\end{aligned}
$$
where we have used Chebychev's inequality and the inequality
$$
\begin{aligned}
& \int_{\{0\,\le \xi\leq M\}} |U_k(\xi)-U(\xi)|\,d\tau_{(x,t)}(\xi)+\int_{\{\xi > M\}} |U_k(\xi)-U(\xi)|\,d\tau_{(x,t)}(\xi)\leq \\
&\qquad  \qquad \le \|U_k-U\|_{L^{\infty}(0,M)} 
+ 2\|U\|_{L^{\infty}(\R)} \int_{\{\xi > M\}} d\tau_{(x,t)}(\xi).
\end{aligned}
$$
Letting $k\to\infty$ we obtain, since $U_k\to U$ uniformly on compact sets in $[0,\infty)$, 
\begin{eqnarray}\label{ci2}
&& \limsup_{j\to\infty}\int_0^T \!\!dt\,\left| \int_{\R} U(u_{n_j})\,\rho(x)\,dx - \int_{\R} U^*(x,t)\,\rho(x)\,dx \,\right| \, \leq \\
& \leq&  2\, \|\rho\|_{C(\overline{\R})}\|U\|_{L^{\infty}(\R)} \left\{\frac{\,T\,\|u_0\|_{\mathcal{M}(\R)}}{M} \, + \iint_{{\rm supp}\rho \times (0,T)} dxdt \!\int_{\{\xi>M\}}d\tau_{(x,t)}(\xi) \right\}\, .  \nonumber
\end{eqnarray}
Since $\tau_{(x,t)}$ is a probability measure, there holds
$\int_{\{\xi > M\}} d\tau_{(x,t)}(\xi)\to 0$ as $M\to\infty$ for a.e.~$(x,t)\in S$, thus by the Dominated Convergence Theorem
$$\iint_{{\rm supp}\rho \times (0,T)} dxdt\!\int_{\{\xi > M\}} d\tau_{(x,t)}(\xi)\to 0 \quad \text{as $M\to\infty$\,.}
$$
Then letting $M\to\infty$ in \eqref{ci2} we obtain \eqref{conv.Uun}. 

\smallskip

\noindent $(iii)$ Now let $U\in C([0,\infty))$ be any function satisfying \eqref{limU}. Arguing as in the proof of Proposition \ref{th.biting2}, let $l_{1m},\, l_{2m}\in C^2([0,\infty))$ ($m\in\N$) satisfy $l_{1m}, l_{2m}\ge 0$ and  $l_{1m}+l_{2m}=1$ in $[0,\infty)$, $\supp l_{1m} \subseteq [0,m+1]$, and $\supp l_{2m} \subseteq [m,\infty)$. Then 
\begin{equation}\label{decoU}
U(u_{n_j})=U(u_{n_j})\,l_{1m}(u_{n_j})+ U(u_{n_j})\,l_{2,m}(u_{n_j})\,,
\end{equation}
and, by \eqref{assini}, for all $\ep>0$  and $m>m_\ep$ 
\begin{equation}\label{estU2}
(C_U-\ep) u_{n_j}\,l_{2m}(u_{n_j}) \leq U(u_{n_j}) l_{2m}(u_{n_j}) \leq  (C_U+\ep) u_{n_j}\,l_{2m}(u_{n_j})\,.
\end{equation}
Since $\|Ul_{1m}\|_{L^{\infty}(\R)} \leq  \|U\|_{C([0,m+1])}<\infty$, the function $Ul_{1m}$ belongs to $C([0,\infty))\cap L^{\infty}(\R)$. Then by part $(ii)$ 
\begin{equation}\label{convU1s}
\int_0^T\left|\int_{\R} [U(u_{n_j})\,l_{1m}(u_{n_j})](x,t)\,\rho(x)\,dx -\int_{\R} U^*_{1m}(x,t)\,\rho(x)\,dx \,\right|\,dt \to 0\,
\end{equation}
as $j\to \infty$, where $\rho \in C_c(\R)$ and $U^*_{1m}$ is defined by \eqref{eq.y.11}.  
By \eqref{estU2} and \eqref{stil1lim} 
\begin{eqnarray*}
&& \int_0^T\left| \int_{\R} \big[U(u_{n_j})\,l_{2m}(u_{n_j})- C_U \,u_{n_j}\,l_{2m}(u_{n_j})\big](x,t)\,\rho(x)\,dx \,\right|\,dt \le\\
&&\qquad \qquad \leq \ep \iint_S |u_{n_j}|\,|\rho(x)|\,dx\,
\leq \ep \,T \|u_0\|_{\mathcal{M}(\R)}\|\rho\|_\infty\,.
\end{eqnarray*}
Then we obtain that 
$$
\begin{aligned}
& \int_0^T \left|  \int_{\R} [U(u_{n_j})\,l_{2m}(u_{n_j})- C_Ul_{2m}^*](x,t) \,\rho(x)\,dx  - C_U\left \langle\mu(\cdot,t),\rho\right\rangle_{\R} \right|\,dt \leq \\
&\le  \ep \,T \|u_0\|_{\mathcal{M}(\R)}\|\rho\|_\infty \!+\!C_U \!\!\!\int_0^T \!\left| \int_{\R} \!\big[u_{n_j}l_{2m}(u_{n_j}) \!- \!l_{2m}^*\big](x,t)\rho(x)\,dx - \left\langle \mu(\cdot,t),\rho\right\rangle_{\R} \right|dt\,,
\end{aligned}
$$
with $l_{2m}^*$ defined as in \eqref{eq.y.12}. The map $u \mapsto u\,l_{2m}(u)$ belongs to $C^2([0,\infty))$, has bounded derivative 
and satisfies \eqref{T0} and \eqref{limU}, with $C_U=1$. 
Then by part $(i)$, \eqref{decoU} and \eqref{convU1s} 
\begin{equation}\label{conv.U4s}
\begin{aligned}
&\limsup_{j\to\infty} \int_0^T \left|  \int_{\R} [U(u_{n_j}) -U^*_{1m}-C_U\,l_{2m}^*](x,t)\rho(x)\,dx  - 
C_U\left \langle\mu(\cdot,t),\rho\right\rangle_{\R} \right|\,dt \,\leq \\ 
&\qquad \qquad  \leq  \ep \,T \,\|u_0\|_{\mathcal{M}(\R)}\|\,\rho\|_\infty\quad\text{if }m>m_\ep.
\end{aligned}
\end{equation}

To complete the proof of \eqref{conv.Uun} we show that
\begin{equation}\label{leb.f}
\lim_{m\to\infty} \iint_{S} |U^*-U^*_{1m}-C_U\,l_{2m}^*|(x,t)\,|\rho(x)|\,dxdt =0\,.
\end{equation} 
By \eqref{estU2},  
\begin{eqnarray*}
&& |U^*-U^*_{1m}-C_U\,l_{2m}^*|(x,t) \! \leq \! \int_{[0,\infty)}\! |U(\xi)-U(\xi)l_{1m}(\xi)-C_U\,\xi\,l_{2m}(\xi)|\,d\tau_{(x,t)}(\xi) =\\
&& = \int_{[0,\infty)} 
| U(\xi)l_{2m}(\xi)-C_U\,\xi\,l_{2m}(\xi)|\,d\tau_{(x,t)}(\xi) \leq \ep \int_{[m,\infty)} \xi\,d\tau_{(x,t)}(\xi)
\leq \ep u_b(x,t)
\end{eqnarray*}
for a.e.~$(x,t)\in S$. 
Since $u_b\in L^{\infty}(0,T;L^1(\R))$ and
$\int_{[m,\infty)} \xi\,d\tau_{(x,t)}(\xi)\to 0$ as $m\to \infty$ for a.e.~$(x,t)\in S$,
\eqref{leb.f} follows from  the Dominated Convergence Theorem. 

Letting $m\to\infty$ in \eqref{conv.U4s}, it follows from \eqref{leb.f} that 
\begin{eqnarray*}
&&\limsup_{j\to \infty}\int_0^T\left| \int_{\R} \left[U(u_{n_j}) -U^*\right](x,t)\,\rho(x)\,dx - C_U\left \langle\mu(\cdot,t),\rho\right\rangle_{\R}\right|\,dt \leq \\
&& \leq \limsup_{m\to\infty} \left(\! \limsup_{j\to \infty}\!\!\int_0^T\!\left| \int_{\R} \!\left[U(u_{n_j}\!) \!-\!U^*_{1m}\!-\!C_Ul_{2m}^*\right]\rho \,dx 
\!-\! C_U\langle\mu(\cdot,t),\rho\rangle_{\R}\right|dt\right)\le  \\
&& \leq \ep \,T\, \|u_0\|_{\mathcal{M}(\R)}\,\|\rho\|_\infty\,,
\end{eqnarray*}
and \eqref{conv.Uun} follows from the arbitrariness of $\ep$. 

Finally \eqref{merc.222} follows from \eqref{conv.Uun}, the separability of $C_c(\R)$ and a diagonal argument; we leave the details to the reader.  
\hfill$\square$

\begin{prop}\label{lebdec}  Let $(H_1)$ hold. Then \eqref{eq.stru} is the Lebesgue decomposition of $u$:
\begin{equation}\label{lebeq}
u_b =  u_r\;\; \text{a.e.~in $S$}\,, \quad\mu=u_s\;\; \text{in $\mathcal{M}(S)$\,.}
\end{equation}
\end{prop}
\begin{proof} Let  
$U$ be a convex function with $U(0)=0$ and $U'\in L^{\infty}(0,\infty)$. 
By \eqref{eq.KlE},  
\begin{eqnarray}\label{vv4}
&&\int_{\R} U(u_{n_j})(x,\bar{t})\,\zeta(x,\bar{t})\,dx - \int_{\R} U(u_{0n_j})(x)\,\zeta(x,0)\,dx \, \le \\
&&\qquad \le\iint_{\R\times (0,\bar{t})} \left\{U(u_{n_j})\,\zeta_t +\Theta_U(u_{n_j})\,\zeta_x\right\}dxdt  \, \nonumber
\end{eqnarray}
for all $\zeta \in C^1([0,T];C^1_c(\R))$ and a.e.~$\bar{t}\in (0,t)$, where
\begin{equation}\label{vv3}
\Theta_U(u):=\int_0^u U'(s)\,\varphi'(s)\,ds + \theta_U  \qquad(\theta_U\in\R)\,.
\end{equation}
Let $U_m(u)= (u-m)\chi_{[m,\infty)}(u)$  and $\theta_{U_m}=0$ ($m\in\N$).
Since 
$U_m(u)/u\to C_{U_m}=1$ and  $\Theta_{U_m}(u)/u\to C_\varphi$ as $u\to\infty$
(with $C_{\varphi}$ as in $(H_1)$), it follows from \eqref{convUun2} that
\begin{equation*}
\int_0^{\bar{t}} \left| \int_{\R} [U_m(u_{n_j})\zeta_t](x,t)\,dx - \int_{\R} [U_m^*\zeta_t](x,t)\,dx -\left\langle \mu(\cdot,t),\zeta_t(\cdot,t)\right\rangle_{\R} \right|\,dt  \to 0
\end{equation*}
and
\begin{equation*}
\int_0^{\bar{t}} \left| \int_{\R} [\Theta_{U_m}(u_{n_j})\zeta_x](x,t)\,dx\! -\! \! \!  \int_{\R} [\Theta_{U_m}^*\zeta_x](x,t)\,dx \! -\! C_{\varphi}\left\langle \mu(\cdot,t),\zeta_x(\cdot,t)\right\rangle_{\R} \right|dt  \to 0
\end{equation*}
as $j\to\infty$, where
\begin{equation*}
U_m^*(x,t):= \int_{[0,\infty)}U_m(\xi)\,d \tau_{(x,t)}(\xi)\,, \quad 
\Theta_{U_m}^*(x,t):= \int_{[0,\infty)}\Theta_{U_m}(\xi)\,d \tau_{(x,t)}(\xi)
\end{equation*}
belong to  $ L^{\infty}(0,T;L^1_{\text{loc}}(\R))$. In particular, setting $\zeta_{\nu}:=\zeta_t+C_{\varphi}\,\zeta_x$, we have that 
\begin{eqnarray}\label{vv7}
&&\iint_{\R\times (0,\bar{t})}  \left\{U_m(u_{n_j})\,\zeta_t +\Theta_{U_m}(u_{n_j})\,\zeta_x\right\}dxdt \to\\
&&\to \iint_{\R\times (0,\bar{t})}  \left\{U_m^*\,\zeta_t +\Theta_{U_m}^*\zeta_x\right\} dxdt +  \int_0^{\bar{t}} \lla\mu(\cdot,t),
\zeta_{\nu}(\cdot,t)\rra_{\R}dt.  \nonumber
\end{eqnarray}

By  \eqref{merc.222} and a diagonal argument, there exist a null set $N\subset (0,T)$ and a subsequence, 
denoted again by $\{u_{n_j}\}$, such that  for all $\bar{t}\in (0,T)\setminus N$ and $m\in\mathbb{N}$ 
\begin{equation}\label{vv6}
\lim_{n_j\to\infty} \int_{\R} U_m(u_{n_j})(x,\bar{t})\,\zeta(x,\bar{t})\,dx= \int_{\R} U_m^*(x,\bar{t})\,\zeta(x,\bar{t})\, dx
+ \lla\mu(\cdot,\bar{t}),\zeta(x,\bar{t})\rra_{\R}.
\end{equation}
Since $\left\{U_m(u_{0n_j})-u_{0n_j}\right\}$ is bounded in $L^{\infty}(\R)$ and converges 
a.e.~to $U_m(u_{0r})-u_{0r}$, it follows from  \eqref{conv.qo.u0n} that
\begin{equation}\label{vv8}
\lim_{n_j\to\infty}\int_{\R} U_m(u_{0n_j})(x)\zeta(x,0)\,dx= \int_{\R} U_m(u_{0r})(x)\zeta(x,0)\,dx
+ \lla u_{0s},\zeta(\cdot,0)\rra_{\R} \,.
\end{equation}

Setting $U=U_m$ in \eqref{vv4}  and letting $j\to\infty$, we obtain from \eqref{vv7}-\eqref{vv8} that 
\begin{eqnarray}\label{vv9bis}
&&\int_{\R} U_m^*(x,\bar{t})\,\zeta(x,\bar{t})\, dx
+ \lla\mu(\cdot,\bar{t}),\zeta(\cdot,\bar{t})\rra_{\R}
\le \iint_{\R\times (0,\bar{t})}  \left\{U_m^*\,\zeta_t +\Theta_{U_m}^*\zeta_x\right\} dxdt +\\
&&\qquad +  \int_0^{\bar{t}} \lla\mu(\cdot,t),
\zeta_{\nu}(\cdot,t)\rra_{\R}dt\, + \int_{\R} U_m(u_{0r})(x)\zeta(x,0)\,dx + \lla u_{0s},\zeta(\cdot,0)\rra_{\R} \, \nonumber
\end{eqnarray}
for all $\bar{t}\in (0,T)\setminus N$ and $m\in\mathbb{N}$. Since
for all $u\ge 0$ (see \eqref{disder})
$$
\text{$0\leq U_m(u)\leq u \chi_{[m,\infty)}(u)$\,, \quad  $|\Theta_{U_m}(u)| = |\varphi(u)-\varphi(m)|\chi_{[m,u)}(u) \leq Mu\,\chi_{[m,\infty)}(u)$ }
$$   we have that 
$ |U^*_m| \leq u_b$, $|\Theta^*_{U_m}|\leq M u_b$, 
$U^*_m\to 0$ and $\Theta^*_{U_m}(x,t)\to 0$ (as $m\to\infty$) a.e.~in $S$. 
Thus, by the Dominated Convergence Theorem and \eqref{vv9bis}, for all $\bar{t}\in (0,T)\setminus N$
\begin{equation}\label{vv9biss}
\lla\mu(\cdot,\bar{t}),\zeta(\cdot,\bar{t})\rra_{\R}\,\le \int_0^{\bar{t}} \lla\mu(\cdot,t),
\zeta_{\nu}(\cdot,t)\rra_{\R}dt + \lla u_{0s},\zeta(\cdot,0)\rra_{\R}.
\end{equation}

Let  $\rho\in C_c^1(\R)$ and $\zeta(x,t)=\rho(x-C_{\varphi}t)$, so $\zeta_{\nu}\equiv 0$. By \eqref{vv9biss},
$\lla\mu(\cdot,\bar{t}),\rho(\cdot\,- C_{\varphi}\bar{t})\rra_{\R}\\ \le \lla u_{0s},\rho\rra_{\R}$.
Hence $\mu(\cdot,\bar{t})$ is singular with respect to the Lebesgue measure and, 
since $\mu(\cdot,\bar{t})=[\mu(\cdot,\bar{t})]_s=\mu_s(\cdot,\bar{t})$ for a.e.~$\bar{t}\in (0,T)$ (see \eqref{us(t)=u(t)s}),
\eqref{lebeq} follows from the uniqueness of the Lebesgue decomposition. 
\end{proof}

The following result is based on the concept of compensated compactness ($e.g.$, see \cite{E}).
\begin{prop}\label{ide} Let $(H_1)$ hold.
Then $\varphi (u_r)= \int_{[0,\infty)}\varphi (\xi)\,d \tau(\xi)$ a.e.~in $S$.
\end{prop}
\begin{proof} 
 Let $U,V\in C^2([0,\infty))\cap L^{\infty}((0,\infty))$ satisfy \eqref{T0}, and assume that $\Theta_U,\,\Theta_V$, 
defined by \eqref{vv3}, belong to $L^{\infty}((0,\infty))$. 
By \eqref{Ma1} there holds 
$$
\text{$\ep \,\|U''(u_n^{\ep})(u_{nx}^{\ep})^2\|_{L^1(S)} \le C_p$
\; and \;\;$\ep \,\|V''(u_n^{\ep})(u_{nx}^{\ep})^2\|_{L^1(S)} \le C_p$}
$$ 
for all $\ep\in(0,1)$ and $n\in\N$, and up to a subsequence 
\begin{equation}\label{C1}
\ep \,U''(u_n^{\ep})(u_{nx}^{\ep})^2\stackrel{*}\rightharpoonup \lambda_n,\quad 
\ep \,V''(u_n^{\ep})(u_{nx}^{\ep})^2\stackrel{*}\rightharpoonup \mu_n\qquad  \mbox{in } \mathcal{M}(S)\quad\text{as  $\ep\to 0$,}
\end{equation} for some $\lambda_n,\mu_n\in\mathcal{M}(S)$.
By the lower semicontinuity of the norm,
\begin{equation}\label{B1}
\|\lambda_n\|_{\mathcal{M}(S)} \le C_p,\quad \|\mu_n\|_{\mathcal{M}(S)} \le C_p
\quad \text{for }n\in\N\,. 
\end{equation}

Let $\zeta\in C_c^2(S)$. 
  Then (see \eqref{RS})
\begin{equation}\label{RRR}
  \ep \! \! \iint_S U''(u_n^{\ep})(u_{nx}^{\ep})^2\zeta\,dxdt =\! \! 
\iint_S \{U(u_n^{\ep})\zeta_t+\Theta_{U,\ep}(u_n^{\ep})\zeta_x+\ep U(u_n^{\ep})\zeta_{xx} \}dxdt,
\end{equation}
where
$\Theta_{U,\ep}(u)=\int_0^u U'(s)\,\varphi'_{\ep}(s)\,ds + \theta_U$,   $\theta_U\in\R.$
By \eqref{disder} and \eqref{stilinfty}, for all $n\in\N$   
\begin{equation*}
|\Theta_{U,\ep}(u_n^{\ep})| \le\int_0^{\|u_{0n}\|_\infty} |U'(s)\varphi'_{\ep}(s)|\,ds  + \left|\theta_U\right|
\le M\!\! \int_0^{\|u_{0n}\|_\infty} |U'(s)|\,ds  +\left|\theta_U\right| \le \gamma_{n,U} 
\end{equation*}
for some $\gamma_{n,U}\ge0$, so for fixed $n\in\N$ the family $\left\{\Theta_{U,\ep}(u_n^{\ep})\right\}_\ep$ 
is uniformly bounded in $L^{\infty}(S)$.
Similar results hold for $V$ and $\Theta_{V,\ep}(u)=\int_0^u V'(s)\,\varphi'_{\ep}(s)\,ds+\theta_V$, 
and letting $\ep\to 0$ in \eqref{RRR} 
along some subsequence $\{\ep_m\}$ (see the proof of Proposition \ref{sl3}) 
it follows from by \eqref{C1} that for all $n\in\N$ and $\zeta\in C_c^1(S)$ 
\begin{equation}\label{N1}
\iint_S\! \{U(u_n)\zeta_t+\Theta_U(u_n) \zeta_x\}dxdt =
\left\langle \lambda_n,\zeta\right\rangle_S\!,
\ \iint_S \!\{V(u_n)\zeta_t+\Theta_V(u_n) \zeta_x \}dxdt=
\left\langle \mu_n,\zeta\right\rangle_S
\end{equation}
where $u_n$ is the entropy solution of the approximating problem $(P_n)$ (see \eqref{cLeb}). 

Let $A\subset\subset S$ be a bounded open set and let $Y_n,Z_n:A\mapsto \R^2$ be defined by
$$
Y_n:=(\Theta_U(u_n),U(u_n)),\quad
 Z_n:=(V(u_n),-\Theta_V(u_n)).
$$
By \eqref{N1}, 
\begin{equation}\label{N3}
{\rm div} \, Y_n=-\lambda_n\,,\;\; {\rm curl} \, Z_n=-\mu_n \quad \text{in $\mathcal{D}'(A)$}\,.
\end{equation}

Since $U,\,\Theta_U,\, V,\,\Theta_V$ are bounded in $(0,\infty)$, the sequences $U(u_n)$, $\Theta_U(u_n)$, $V(u_n)$ and $\Theta_V(u_n)$ are bounded in $L^1(A)$ and uniformly integrable, and, by  Theorem \ref{th.Ball}, 
$$
U(u_n)\rightharpoonup U^*:= \int_{[0,\infty)}U(\xi)\,d \tau_{(\cdot,\cdot)}(\xi),\quad 
\Theta_U(u_n)\rightharpoonup \Theta_U^*:= \int_{[0,\infty)}\Theta_U(\xi)\,d \tau_{(\cdot,\cdot)}(\xi),
$$
$$
V(u_n)\rightharpoonup V^*:= \int_{[0,\infty)}V(\xi)\,d \tau_{(\cdot,\cdot)}(\xi),\quad 
\Theta_V(u_n)\rightharpoonup \Theta_V^*:= \int_{[0,\infty)}\Theta_V(\xi)\,d \tau_{(\cdot,\cdot)}(\xi)
$$
in $L^1(A)$, where $\tau_{(\cdot,\cdot)}$ denotes the disintegration of  the Young measure $\tau$ 
associated with $\{u_n\}$.
Since the sequences $U(u_n)$, $\Theta_U(u_n)$, $V(u_n)$ and $\Theta_V(u_n)$ are bounded in $L^{\infty}(A)
\subset L^2(A)$, they also converge weakly in $L^2(A)$, so 
\begin{equation*}%\label{C3}
 Y_n\rightharpoonup Y^*:=(\Theta_U^*,U^*) \,, \;\;  Z_n\rightharpoonup Z^*:=(V^*,-\Theta_V^*) \quad \text{in $\left[L^2(A)\right]^2$}\,.
\end{equation*}
By a similar argument 
\begin{eqnarray}\label{C4}
&& Y_n\cdot Z_n := \Theta_U(u_n)V(u_n)-\Theta_V(u_n)U(u_n)
\rightharpoonup \\
&& \qquad \rightharpoonup \int_{[0,\infty)}\left[\Theta_U(\xi)V(\xi)-\Theta_V(\xi)U(\xi)\right]d \tau_{(\cdot,\cdot)}(\xi)\,
\quad \text{in $L^2(A)$}\,. \nonumber
\end{eqnarray}
By \eqref{B1} and \eqref{N3}, $\{{\rm div} \, Y_n\}$ and $\{{\rm curl} \, Z_n\}$ are precompact in $W^{-1,2}(A)$ (see \cite[Chapter 1, Corollary 1]{E}) and, by the div-curl lemma, 
\begin{equation}\label{C5}
 Y_n\cdot Z_n \to Y^*\cdot Z^*= \Theta_U^*V^*-\Theta_V^*U^* \quad \text{in $\mathcal{D}'(A)$}\,. 
  \end{equation}
By \eqref{C4} and \eqref{C5},  
\begin{equation}\label{ug1}
 \int_{[0,\infty)} [\Theta_U(\xi)-\Theta_U^*]V(\xi)d \tau(\xi)
= \int_{[0,\infty)}[U(\xi)-U^*]\Theta_V(\xi)\, d \tau(\xi)\quad \text{a.e.~in $ A$}. 
  \end{equation}

For every $U$ as above with $U'>0$ in $(0,\infty)$, by a standard approximation argument we may choose 
$V(u)=|U^*-U(u)|$, so $\Theta_V(u)=\sgn (U(u)-U^*)[\Theta_U(u)- \Theta_U(U^{-1}(U^*))]$ and, by \eqref{ug1},
\begin{equation}\label{ug2}
\left[\Theta_U^*- \Theta_U(U^{-1}(U^*)) \right] \int_{[0,\infty)}|U^*-U(\xi)|\, d \tau(\xi)=0\,.
\end{equation}
Let $U_k\in C^2([0,\infty))\cap L^{\infty}((0,\infty))$ satisfy \eqref{T0} and
\begin{equation}\label{Tk}
U_k(0)=0,\quad 0<U_k'\le U_{k+1}'\le 1 \;  \text{ in $[0,\infty)$}, \quad U_k'(u)\to 1 \; \text{for $u\ge 0$ as $k\to\infty$\,.}
\end{equation}
By \eqref{disder},
$|\Theta_{U_k}(u)|\leq \int_0^u U'_k(s)|\varphi'(s)|\,ds + |\theta_{U_k}|\leq M\,U_k(u) + |\theta_{U_k}|$,
 thus $\Theta_{U_k}$ is bounded in
 $(0,\infty))$ for every $k\in\N$. 
We claim that, as $k\to\infty$,
\begin{equation}\label{CT1}
 U_k^*:= \int_{[0,\infty)}U_k(\xi)\,d \tau(\xi)\to 
 u_r \quad\mbox{ a.e.~in $A$}\,,
 \end{equation}
\begin{equation}\label{CT2} 
\Theta_{U_k}^*-\, \Theta_{U_k}({U_k}^{-1}(U_k^*))\,
\to \int_{[0,\infty)}\varphi(\xi)\,d \tau(\xi) - \varphi(u_r) \quad\mbox{ a.e.~in $A$}\,,
 \end{equation}
where 
$\Theta_{U_k}^*:= \int_{[0,\infty)}\Theta_{U_k}(\xi)\,d \tau(\xi)$
(recall that $\varphi\in L^1([0,\infty);d\tau_{(x,t)})$, see Remark \ref{rentro}). 
By \eqref{CT1} and the Dominated Convergence Theorem, 
for a.e.~$(x,t)\in A$
$$
 \int_{[0,\infty)}|U_k^*(x,t)-U_k(\xi)|\, d \tau_{(x,t)}(\xi)\, \to \,\int_{[0,\infty)}| u_r(x,t)-\xi|\, d \tau_{(x,t)}(\xi)\quad\text{as }k\to\infty,
 $$
 since $0\le U_k(\xi)\le\xi$ for all $k\in\N$ and $I(\xi):=\xi$ belongs to $L^1([0,\infty), d \tau_{(x,t)})$ (recall that, by \eqref{lebeq} and the definition of $u_b$ in \eqref{merc.2},   
 $u_r(x,t)=\int_{[0,\infty)}\xi\,d \tau_{(x,t)}(\xi)<\infty$ for a.e.~$(x,t)\in S$).
Letting $k\to\infty$ \eqref{ug2}, with $U=U_k$, we obtain that for a.e.~$(x,t)\in A$
\begin{equation*}
\left[ \int_{[0,\infty)}\varphi(\xi)\,d \tau_{(x,t)}(\xi)- \varphi(u_r)(x,t) \right] \int_{[0,\infty)}| u_r(x,t)-\xi|\, d \tau_{(x,t)}(\xi)=0,
 \end{equation*}
and Proposition \ref{ide} follows from the arbitrariness of $A$.

It remains to prove \eqref{CT1} and \eqref{CT2}. By \eqref{Tk} and the Monotone Convergence Theorem, 
$U_k(\xi)\to \xi$ for any $\xi\in[0,\infty)$, and   
\eqref{CT1} follows (recall that  $I(\xi)=\xi \in L^1([0,\infty),d\tau)$). Concerning \eqref{CT2} we observe  that
\begin{equation}\label{CT0} 
\Theta_{U_k}^*\!\!- \Theta_{U_k}({U_k}^{-1}\!(U_k^*))
\!=\!\!\int_{[0,\infty)}\!\!\left(\int_0^{\xi} \!\!U_k'(s)\varphi'(s)ds\!\right)\!d \tau(\xi)\! - \!\!\!\int_0^{{U_k}^{-1}({U_k}^*)} \!\!U_k'(s)\varphi'(s)ds.
 \end{equation}
Since $U_k'(\xi)\to 1$ and $|U_k'(\xi)\varphi'(\xi)|\le M$  for $\xi\ge0$ (see \eqref{Tk} and \eqref{disder}),
it follows from the Dominated Convergence Theorem that 
\begin{equation}\label{CT3} 
\int_{[0,\infty)}\left(\int_0^{\xi} U_k'(s)\,\varphi'(s)\,ds\right)d \tau_{(x,t)}(\xi) \to   \int_{[0,\infty)}\varphi(\xi)\,d \tau_{(x,t)}(\xi).
 \end{equation}
On the other hand, 
\begin{eqnarray}\label{CT4} 
&&\int_0^{{U_k}^{-1}(U_k^*(x,t))} U_k'(s)\,\varphi'(s)\,ds - \varphi(u_r)(x,t)  \,= \\
&&\qquad =\int_0^{ u_r(x,t)} [U_k'(s)-1]\,\varphi'(s)\,ds + \int_{ u_r(x,t)}^{{U_k}^{-1}(U_k^*(x,t))} U_k'(s)\,\varphi'(s)\,ds\,. \nonumber
 \end{eqnarray}
Arguing as before one shows that the first term in the right-hand side of \eqref{CT4} vanishes as $k\to\infty$. As for the second
term we observe that, by \eqref{Tk} and  \eqref{CT1},  
\begin{eqnarray*}
&&\left|\int_{ u_r(x,t)}^{{U_k}^{-1}(U_k^*(x,t))} U_k'(s)\,\varphi'(s)\,ds\right|  \,\le\, M \left|\,  u_r(x,t)- {U_k}^{-1}(U_k^*(x,t)) \right| \,\le \\
&&\qquad \le  M\left( \left|\,  u_r(x,t)- U_k^{-1}( u_r(x,t)) \right|\,+\,\sup_{s\in I_{\delta}( u_r(x,t))}\frac{1}{U_1'(s)}\left|\,  u_r(x,t)- U_k^*(x,t) \right| \right) 
 \nonumber
 \end{eqnarray*}
for some $\delta>0$ and all $k\in\N$ sufficiently large, where $I_\delta(q)\equiv (q-\delta,q+\delta)$. Hence
\begin{equation}\label{CT5} 
\int_0^{{U_k}^{-1}(U_k^*(x,t))} U_k'(s)\,\varphi'(s)\,ds\,\to \varphi(u_r)(x,t)  \quad\mbox{ for a.e.~$(x,t)\in A$},
 \end{equation}
and we obtain \eqref{CT2} from \eqref{CT0}, \eqref{CT3} and \eqref{CT5}. 
 \end{proof}

To prove the second part of Theorem \ref{th.exi1} we need the following result, which  
characterizes the disintegration of the Young measure $\tau$. 

\begin{prop}\label{idebis2} 
Let $(H_1)$ hold and $\varphi\in C^1([0,\infty))$ satisfy for all $\bar{u}>0$ either \eqref{exH_2} or
\begin{equation}\label{exH_2''}
\begin{cases}
\mbox{$\exists$ $a>0$, $b\in (0,\infty]$ such that $\varphi'$ is constant in $I_{a,b}=[\bar{u}-a,\bar{u}+b]$ and, if $b<\infty$,
$\varphi'$}
\\
\mbox{is strictly monotone in 
$[\bar{u}+b,\bar{u}+\tilde{b}]$ and $[\bar{u}-\tilde{a},\bar{u}-a]$ for some $\tilde{b}>b$ and $\tilde{a}\in (a,\bar{u})$.}
\end{cases}
\end{equation}
Then for a.e.~$(x,t)\in S$ the following holds:

 \noindent $(i)$ if $u_r(x,t)=0$, then
$\tau_{(x,t)}=\delta_0$;

\noindent $(ii)$ if $\varphi'$ is strictly monotone in $I_{a,b}\!=\![u_r(x,t)\!-\!a,u_r(x,t)\!+\!b]$
with $a,b\geq 0$, $a+b>0$, 
\begin{equation}\label{delta'}
\tau_{(x,t)}=\delta_{u_r(x,t)}\,;
\end{equation}
\noindent $(iii)$ if $\varphi'$ is constant in the above interval $I_{a,b}$ for some $a>0$, $b>0$, 
 then
\begin{equation}\label{deltab}
{\rm supp}\,\tau_{(x,t)}\subseteq I_{(x,t)} \quad \text{for a.e.~$(x,t)\in S$}\,,
\end{equation}
where $I_{(x,t)}\supseteq I_{a,b}$ is the maximal interval 
where $\varphi'(\cdot)\equiv\varphi'(u_r(x,t))$. 
\end{prop} 

\begin{proof} Let $(x,t)\in S$ be fixed. If $u_r(x,t)=0$ it follows from \eqref{lebeq} and the definition of $u_b$ in \eqref{merc.2} that
$\int_{[0,\infty)}\xi\,d \tau_{(x,t)}(\xi)=0$, which implies  part $(i)$: $\tau_{(x,t)}=\delta_0$. 

So let $u_r(x,t)>0$. Let $l_1:=u_r(x,t)$,  $l_2>l_1$ and  
\begin{equation*}%\label{defVk'}
V_k(u):= k(u-l_1)\chi_{\left(l_1,l_1+\frac1k\right)}(u) + \chi_{\left[l_1+\frac1k,l_2\right)}(u)+ k\left(l_2+\frac1k-u\right)\chi_{\left[l_2,l_2+\frac1k\right)}(u)
\end{equation*}
for  $u\geq 0$ and sufficiently large $k\in\N$. Then $V_k(u)\to \chi_{\left(l_1,l_2\right]}(u)$ as $k\to\infty$, and
$$
\Theta_{V_k}(u)\!=\!\!\int_0^u\!V_k'(s)\varphi'(s)ds\to \varphi'(l_1)\chi_{\left(l_1,l_2\right]}(u)+ [\varphi'(l_1)\!-\!\varphi'(l_2)]\chi_{\left(l_2,\infty\right)}(u) \quad (u\ge0)\,.
$$
By standard approximation arguments, 
\eqref{ug1} is satisfied with $U=U_k$ and $V=V_k$, where 
$\{U_k\}$ is the sequence in the proof of Proposition \ref{ide} (see \eqref{Tk}):
\begin{equation*}
 \int_{[0,\infty)}\![\Theta_{U_k}(\xi)-\Theta_{U_k}^*(x,t)]V_k(\xi)d \tau_{(x,t)}(\xi)
=\! \! \int_{[0,\infty)}\![U_k(\xi)-U^*_k(x,t)]\Theta_{V_k}(\xi) d \tau_{(x,t)}(\xi). 
\end{equation*} 
Letting  $k\to\infty$ 
and arguing as  in the proof of Proposition \ref{ide}, we obtain that
\begin{eqnarray*}
	&&U_k(\xi)-U^*_k(x,t) \to\, \xi-\int_{[0,\infty)}\xi\,d\tau_{(x,t)}(\xi)=\xi -u_r(x,t)=\xi-l_1\,,\\
	&&\Theta_{U_k}^*\!(x,t)\!-\!\Theta_{U_k}\!(\xi)\!  
	\to\!\! \int_{[0,\infty)}\!\!\varphi(\xi)d\tau_{(x,t)}\!(\xi)\!-\!\varphi(\xi)\!=\!\varphi(u_r)(x,t)\!-\!\varphi(\xi)\!=\!\varphi(l_1)\!-\!\varphi(\xi)
\end{eqnarray*}
 for all $\xi\geq 0$ (see \eqref{lebeq} and Proposition \ref{ide}). This implies that  
\begin{eqnarray*}
	&&\int_{[0,\infty)}\left[\Theta_{U_k}(\xi)-\Theta_{U_k}^*(x,t)\right]V_k(\xi)\, d \tau_{(x,t)}(\xi)\to \int_{(l_1,l_2]} [\varphi(\xi)-\varphi(l_1)]
	\,d\tau_{(x,t)}(\xi)\,,\\
 	&&\int_{[0,\infty)}\left[U_k(\xi)-U^*_k(x,t)\right]\Theta_{V_k}(\xi)\, d \tau_{(x,t)}(\xi)\to  \int_{\left(l_1,l_2\right]}\varphi'(l_1)(\xi-l_1)\, d \tau_{(x,t)}(\xi)+\\
	&&\qquad \qquad \qquad \qquad \qquad \qquad \qquad +\left[\varphi'(l_1)-\varphi'(l_2)\right]\int_{\left(l_2,\infty\right)}\left(\xi-l_1\right)\, d \tau_{(x,t)}(\xi)\,,
\end{eqnarray*}
whence 
\begin{equation}\label{ln'}
\int_{(l_1,l_2]}[\varphi(\xi)-\varphi(l_1)-\varphi'(l_1)(\xi-l_1)] d \tau_{(x,t)}(\xi)
\!=\! [\varphi'(l_1)-\varphi'(l_2)]\int_{(l_2,\infty)}(\xi-l_1) d \tau_{(x,t)}(\xi).
\end{equation}

Similarly, let $l_0\in (0,l_1)$ and set 
\begin{equation*}
\tilde{V}_k(u):= k(u-l_0)\chi_{\left[l_0,l_0+\frac1k\right]}(u) + \chi_{\left(l_0+\frac1k,l_1-\frac1k\right)}(u)+ k\left(l_1-u\right)\chi_{\left[l_1-\frac1k,l_1\right]}(u)\,.
\end{equation*}
Then 
$\tilde{V}_k(u)\to \chi_{\left(l_0,l_1\right)}(u)$, 
and
$$
\Theta_{\tilde{V}_k}(u)=\int_{l_1}^u \tilde{V}'_k(s)\,\varphi'(s)\,ds\to \varphi'(l_1)\chi_{\left(l_0,l_1\right)}(u)+ [\varphi'(l_1)-\varphi'(l_0)]\chi_{\left[0,l_0\right]}(u)\quad (u\ge0)\,.
$$
Letting $k\to \infty$ in \eqref{ug1} with $U=U_k$ as above and $V=\tilde{V}_k$, 
we obtain that 
\begin{equation}\label{lnbis}
\int_{(l_0,l_1)} \! [\varphi(\xi)\!-\!\varphi(l_1)\!-\!\varphi'(l_1)(\xi-l_1\!)\!]d\tau_{(x,t)}(\xi)
\!=\![\varphi'(l_1)\!-\!\varphi'(l_0)]\!\!\int_{[0,l_0]} \!(\xi-l_1)d\tau_{(x,t)}(\xi).
\end{equation}

By \eqref{exH_2} and \eqref{exH_2''}, 
we can distinguish two cases.

\smallskip
\noindent  $(a)$ If $\varphi$ is strictly convex or strictly concave in $[l_1,l_2]$, 
it follows from \eqref{ln'} that 
$$
\int_{(l_1,l_2]}\!|\varphi(\xi)\!-\!\varphi(l_1)\!-\!\varphi'(l_1)(\xi\!-\!l_1)| d \tau_{(x,t)}(\xi)\! 
+\!|\varphi'(l_1)\!-\!\varphi'(l_2)|\!\int_{(l_2,\infty)}\! |l_1\!-\!\xi| d \tau_{(x,t)}(\xi)\!=\!0,
$$  
where
$$
\text{$ \chi_{(l_1,l_{2}]}(\xi)|\varphi(\xi)-\varphi(l_1)-\varphi'(l_1)(\xi-l_1)|>0$\;\; and \;\;$| \varphi'(l_1)-\varphi'(l_{2})|>0$\,.}
$$ 
This implies that 
${\rm supp}\,\tau_{(x,t)}\subseteq\left[0,l_1\right]$. 
 Since  $\tau_{(x,t)}$ is a probability measure and $l_1:=u_r(x,t)$, 
$$
u_r(x,t)=\int_{\left[0,u_r(x,t)\right]}\xi\,d \tau_{(x,t)}(\xi)=\int_{\left[0,u_r(x,t)\right]}\left[\xi-u_r(x,t)\right]\,d \tau_{(x,t)}(\xi) + u_r(x,t)
$$
(see \eqref{merc.2} and \eqref{lebeq}), thus
$$
\int_{\left[0,u_r(x,t)\right]}\left|\xi-u_r(x,t)\right|\,d \tau_{(x,t)}(\xi) \quad \Rightarrow \quad \tau_{(x,t)}(\left[0,u_r(x,t)\right))=0\,.
$$
Hence ${\rm supp}\,\tau_{(x,t)}=\{ u_r(x,t)\}$ and \eqref{delta'} follows since  $\tau_{(x,t)}$ is a probability measure.

Similarly, if $\varphi$ is strictly convex or strictly concave in $(l_0,l_1)$, 
it follows from \eqref{lnbis} that $\tau_{(x,t)}(\left [0,l_1 \right))=0$ (we omit the details). Thus, supp$\,\tau_{(x,t)}\subseteq [l_1,\infty)$ and arguing as above we obtain  \eqref{delta'}. 

\smallskip

\noindent  $(b)$ If $\varphi$ is affine in $[l_1-c,l_1+c]$ for some $c>0$, let $I=[\bar{l}_0,\bar{l}_2]$ be the maximal interval containing $l_1$ where $\varphi'(\xi)=\varphi'(l_1)$.
If $I=[0,\infty)$, \eqref{deltab}  
is satisfied. If $\bar{l}_2<\infty$, by \eqref{exH_2''} and the maximality of $I$, 
$\varphi$ is strictly convex (or concave) in $[\bar{l}_2,\bar{l}_2+b]$ for some  $b>0$  
(and affine in $[l_1,\bar{l}_2]$). By \eqref{ln'}, with $l_2\in (\bar{l}_2,\bar{l}_2+b)$, we obtain that 
$$
\int_{(\bar{l}_2,l_2]}\! |\varphi(\xi)\!-\!\varphi(l_1)\!-\!\varphi'(l_1)(\xi\!-\!l_1)| d \tau_{(x,t)}(\xi) \!+\! 
|\varphi'(l_1)\!-\!\varphi'(l_2)|\int_{(l_2,\infty)}\! |l_1\!-\!\xi| d \tau_{(x,t)}(\xi)\!=\!0\,,
$$ 
where 
$$\text{$ \chi_{\left(\bar{l}_2,l_{2}\right]}(\xi)\,\left|\,\varphi(\xi)-\varphi(l_1)-\varphi'(l_1)(\xi-l_1)\right|>0$ and 
$ | \varphi'(l_1)-\varphi'(l_{2})|>0$\,.}
$$
It follows that $\tau_{(x,t)}\big((\bar{l}_2,\infty)\big)=0$, 
whence  supp$\,\tau_{(x,t)}\subseteq [0,\bar{l}_2]$. 
Similarly, if $\bar{l}_0>0$, by \eqref{exH_2''} and the maximality of $I$,  
$\varphi$ is strictly convex (or concave) in $[\bar{l}_0-a,\bar{l}_0]$ for some $a>0$   (and affine in  $[\bar{l}_0,l_1]$). 
Arguing as before, we obtain from \eqref{lnbis}, with $l_0\in (\bar{l}_0-a,\bar{l}_0)$, that supp$\,\tau_{(x,t)}\subseteq [\bar{l}_0,\infty)$ (we omit the details). Summing up we obtain  \eqref{deltab}: supp$\,\tau_{(x,t)}\subseteq [0,\bar{l}_2]\cap [\bar{l}_0,\infty)=I$.
\end{proof}

\begin{remark}\label{come} 
If \eqref{exH_2} is satisfied for all $\bar{u}>0$, it follows from \eqref{delta'} and standard properties of narrow convergence of Young measures (see \cite{V1}) that $u_{n_j}\to u_r$ in measure, where $\{u_{n_j}\} $ is the subsequence in Proposition \ref{th.biting1}. 
Therefore, up to a subsequence, 
$u_{n_j}\to u_r$ a.e.~in  $S$.  
Hence, if $\varphi$ is  bounded, it follows from  the Dominated Convergence Theorem that
$\varphi(u_{n_j})\to \varphi(u_r)$ in $L^1((-L,L)\times (0,T))$ for all $L>0$.
\end{remark}

Now we can prove Theorem \ref{th.exi1}. 

\smallskip

\noindent {\em Proof of Theorem \ref{th.exi1}.} Let $\zeta\in C^1([0,T];C^1_c(\R))$, $\zeta(\cdot,T)=0$ in $\R$, and 
let $L>0$ be such that supp$\,\zeta\subset (-L,L)\times [0,T]$. 
By \eqref{convUun2}, with $U(u)=u$ and $U(u)=\varphi(u)$, 
\begin{equation*}%\label{aa1}
\iint_S  u_{n_j} \zeta_t\,dxdt \to \iint_S  u_r \zeta_t\,dxdt + \int_0^T \lla u_s(\cdot,t),\zeta_t(\cdot,t)\rra_{\R}dt \,,
\end{equation*}
\begin{equation*}%\label{aa3}
\iint_S\varphi(u_{n_j}) \zeta_x \,dxdt \to \iint_S\varphi^* \zeta_x \,dxdt +
C_{\varphi} \!\!\int_0^T \lla u_s(\cdot,t),\zeta_x(\cdot,t)\rra_{\R}dt\,
\end{equation*}
(see also \eqref{lebeq}).
Letting $j\to\infty$ in \eqref{eq.Klc}, with $n=n_j$, we obtain \eqref{eqy}. 
Inequality  \eqref{misey} is proven similarly, since by arguing as in Proposition \ref{th.biting2} we get
$$E(u_{0n_j}) \stackrel{*}\rightharpoonup E(u_{0r})+C_Eu_{0s}\quad\mbox{in}\ \,\mathcal{M}(\R)$$
(in this regard, see also \eqref{conv.qo.u0n}).
Thus the function $u\in L^{\infty}(0,T;\mathcal{M}^+(\R))$ 
given by Proposition \ref{th.biting1} is an entropy solution of problem $(P)$  in the sense of Young measures. 
By Proposition \ref{ide}, it is also a solution in the sense of Definition \ref{deso}. This proves the first part of  the theorem. The second part is an immediate consequence of Proposition \ref{idebis2}:
in fact,  \eqref{misei} follows from \eqref{misey} and \eqref{delta'}.
\hfill$\square$

\smallskip

Let us end this section by proving Proposition \ref{ac}.
\smallskip

\noindent {\em Proof of Proposition \ref{ac}.} For every $\tilde{\zeta}\in C^1([0,T];C^1_c(\R))$, $\tilde{\zeta}(\cdot,T)=0$, we set $E(u)=U_m(u)=(u-m)\chi_{\{u>m\}}(u)$ and $F(u)=F_m(u)=\int_0^uU'_m(\xi)\,\varphi'(\xi)\,d\xi= (\varphi(u)-\varphi(m))\chi_{\{u>m\}}(u)$ in the entropy inequalities \eqref{misey} ($m\in\mathbb{N}$). Then we get
\begin{eqnarray*}%\label{deri.us1}
&&\iint_{S}\left\{U^*_m\tilde{\zeta}_t+F^*_m\tilde{\zeta}_x\right\}\,dxdt +\int_0^T\left\langle u_s(\cdot,t),\tilde{\zeta}_t(\cdot,t)\right\rangle_{\R}dt + \\
&&+ C_{\varphi} \int_0^T \left\langle u_s(\cdot,t),\tilde{\zeta}_x(\cdot,t)\right\rangle_{\R}dt \geq  - \int_{\R} U_m(u_{0r})\tilde{\zeta}(x,0)\,dx -\left\langle u_{0s},\tilde{\zeta}(\cdot,0)\right\rangle_{\R}\,,\nonumber
\end{eqnarray*}
where, for a.e.~$(x,t)\in S$,  
$$U^*_m(x,t):=\int_{[0,\infty)} U_m(\xi)\,d\tau_{(x,t)}(\xi)\,,\quad F^*_m(x,t):=\int_{[0,\infty)}F_m(\xi)\,d\tau_{(x,t)}(\xi)\,.$$
As in the proof of Proposition \ref{lebdec}, there holds 
$\iint_{S}\left\{U^*_m\tilde{\zeta}_t+F^*_m\tilde{\zeta}_x\right\}\,dxdt\to0$ and 
$\int_{\R} U_m(u_{0r})\tilde{\zeta}(x,0)\,dx\to0$ as $m\to\infty$,
whence
\begin{eqnarray}\label{deri.us2}
&&\int_0^T\left\langle u_s(\cdot,t),\tilde{\zeta}_t(\cdot,t)\right\rangle_{\R}dt + C_{\varphi} \int_0^T \left\langle u_s(\cdot,t),\tilde{\zeta}_x(\cdot,t)\right\rangle_{\R}dt \geq -\left\langle u_{0s},\tilde{\zeta}(\cdot,0)\right\rangle_{\R}\,.
\end{eqnarray}

Let $\zeta\in C([0,T];C_c(\R))$. By definition of $L^{\infty}(0,T;\mathcal{M}(\R))$ (see Definition \ref{dli}), the map
$t \mapsto \left\langle u_s(\cdot,t), \zeta(\cdot,t)\right\rangle_{\R}$ belongs to $L^{\infty}(0,T)$. Hence  
\begin{equation}\label{us.pl}
\lim_{h\to 0} \frac{1}{h} \int_{\bar{t}}^{\bar{t}+h}\left\langle u_s(\cdot,t),\zeta(\cdot,t)\right\rangle_{\R}dt =\left\langle u_s(\cdot,\bar{t}), \zeta(\cdot,\bar{t})\right\rangle_{\R}\quad\mbox{for every}\ \,\bar{t}\in (0,T)\setminus N\,
\end{equation}
for some null set $N\subset (0,T)$ (by separability arguments, $N$ is independent of $\zeta$; 
see the proof of \cite[Lemma 3.1]{PST}). 
Let $t_1,t_2\in (0,T)\setminus N$, $0<t_1<t_2<T$. 
By standard approximation arguments we can choose 
$\tilde{\zeta}(x,t)= g_h(t)\zeta(x,t)$ in \eqref{deri.us2}, where 
\begin{equation}\label{defigh}
g_h(t):=\frac{1}{h}(t-t_1)\chi_{\{t_1\leq t\leq t_1+h\}}(t)+\chi_{\{t_1+h<t<t_2\}}(t)+ \frac{1}{h}(t_2+h-t)\chi_{\{t_2\leq t\leq t_2+h\}}(t)
\end{equation}
and $h\in (0,\min\{t_2-t_1,T-t_2\})$.
Letting $h\to 0$ in \eqref{deri.us2} we obtain that 
\begin{equation}\label{absco1}
\lla u_s(\cdot,t_2),\zeta(\cdot,t_2)\rra_{\R} \le
\int_{t_1}^{t_2}\left\langle  u_s(\cdot,t), 
\zeta_{\nu}(\cdot,t)\right\rangle_{\R}dt + \lla u_s(\cdot,t_1),\zeta(\cdot,t_1)\rra_{\R}.
\end{equation}

Similarly, let $f_h(t):= \chi_{\{0\leq t<t_2\}}(t)+ \frac{1}{h}(t_2+h-t)\chi_{\{t_2\leq t\leq t_2+h\}}(t)$.
Setting $\tilde{\zeta}(x,t)= f_h(t)\zeta(x,t)$
in \eqref{deri.us2} and letting $h\to 0^+$ we
obtain that 
\begin{equation}\label{absco1bis}
\lla u_s(\cdot,t),\zeta(\cdot,t)\rra_{\R} \le
\int_{0}^{t}\left\langle  u_s(\cdot,\tau), 
\zeta_{\nu}(\cdot,\tau)\right\rangle_{\R}d\tau +
\lla u_{0s},\zeta(\cdot,0)\rra_{\R}\,.
\end{equation}
Arguing as in the last part of the proof of Proposition \ref{lebdec}, we obtain \eqref{absco2} and \eqref{absco22} from 
respectively \eqref{absco1} and \eqref{absco1bis} (we omit the details). 

\smallskip

\noindent $(ii)$ It follows from \eqref{ewf} that for a.e. $\tau\in (0,T)$ and $m\in\N$
\begin{equation}\label{eq.cm1}
\left\langle u(\cdot,\tau),\rho_m\right\rangle_{\R}-\left\langle u_0,\rho_m\right\rangle_{\R} =\int_0^{\tau} \left\{\int_{\Omega_m} \varphi(u_r)(x,t)\,\rho_m'\,dx +C_{\varphi}\left\langle u_s(\cdot,t)\lefthalfcup \Omega_m,\rho_m'\right\rangle_{\R}\right\}\,dt\,,
\end{equation}
where $\{\rho_m\}\subset C^1_c(\R)$ is such that 
 $\rho_m=1$ in $[-m,m]$, supp$\,\rho_m\subseteq [-m-1,m+1]$, $0\leq \rho_m\leq 1$ and $|\rho_m'|\leq 2$ in $\R$,
 and $\Omega_m:=[-m-1,-m]\cup [m,m+1]$. 
Since $u_s\in L^{\infty}(0,T;\mathcal{M}^+(\R))$ and $\varphi(u_r)\in L^{\infty}(0,T;L^1(\R))$, a routine proof shows that 
$$
\lim_{m\to\infty}\int_0^{\tau}\!\!\!\! \int_{\R} \varphi(u_r)(x,t)\,\rho_m'(x)\,dxdt=
\lim_{m\to\infty}\int_0^{\tau}\left\langle u_s(\cdot,t)\lefthalfcup \Omega_m,\rho_m'\right\rangle_{\R}\,dt =0\,.
$$
Since $\rho_m(x)\to 1$ for all $x\in\R$, we also get that
$\langle u(\cdot,\tau),\rho_m\rangle_{\R}\to \|u(\cdot,\tau)\|_{\mathcal{M}(\R)}$ and 
$\langle u_0,\rho_m\rangle_{\R}\to\|u_0\|_{\mathcal{M}(\R)}$ as $m\to\infty$. Letting $m\to\infty$ in \eqref{eq.cm1} we obtain claim $(ii)$. 
\hfill$\square$

%%%%%%%%%%%%%%%%%%%%%%%%%%%%%%%%%%%%%%%%%%%%%%%%%%%%%%%%%%%%%%%%%%%%%%%%%%%%%%%%%%%%%%%%%%%%%%%%%%%%%%%%%%%%%%%%%%

\section{Regularity: Proofs
}\label{sec.ps1}
\setcounter{equation}{0}

The first regularity result which we prove is Proposition \ref{condini}. Hence we need
\begin{lem}\label{lefdp}
Let  $(H_1)$ be satisfied.  Let $(u,\tau)$ be a Young measure solution of problem $(P)$. Then there exists a null set  $F^*\subset (0,T)$ such that for every $t_0,t_1\in (0,T)\setminus F^*$, $t_0<t_1$ and any $\rho\in C^1_c(\R)$ there holds 
\begin{equation}\label{fdp}
\lla u(\cdot,t_1),\rho\rra_{\R}-\lla u_0,\rho\rra_{\R}=\int_0^{t_1}\left\{\int_{\R} \varphi^*(x,t)\rho'(x)\,dx + C_{\varphi}\lla u_s(\cdot,t),\rho'\rra_{\R} \right\}\,dt\,,
\end{equation} 
\begin{equation}\label{fdp bis}
\lla u(\cdot,t_1),\rho\rra_{\R}-\lla u(\cdot,t_0),\rho\rra_{\R}=\int_{t_0}^{t_1}\left\{\int_{\R} \varphi^*(x,t)\rho'(x)\,dx + C_{\varphi}\lla u_s(\cdot,t),\rho'\rra_{\R} \right\}\,dt\,.
\end{equation} 
\end{lem}

\begin{proof} Since $u \in L^{\infty}(0,T;\mathcal{M}^+(\R))$, there exists a null set $F_0\subseteq (0,T)$
such that the spatial disintegration $u(\cdot,t)\in\mathcal{M}^+(\R)$ is defined for every $t\in (0,T)\setminus F_0$.  Arguing as in the proof of \cite[Lemma 3.1]{PST} we can show that there exists a null set $F^*\subset (0,T)$, $F_0\subseteq F^*$, such that for every $\rho\in C_c(\R)$ and $t\in (0,T)\setminus F^*$ 
\begin{eqnarray}\label{eq.punti.lebesgue}
&& 
\lim_{q\to\infty}\left\{2q \,\int_{t-\frac{1}{q}}^{t+\frac{1}{q}}\big|\lla u(\cdot,s),\rho\rra_{\R}-\lla u(\cdot,t),\rho\rra_{\R}\big|\,ds \right \}=0\,.
\end{eqnarray}

The proof of \eqref{fdp} is based on \eqref{eqy} and \eqref{eq.punti.lebesgue}. Let $\rho\in C^1_c(\R)$ and $t_1\in(0,T)\setminus F^*$. By standard regularization arguments we can set $\zeta= \rho(x)k_q(t)$ in \eqref{eqy}, with $q \ge \frac{1}{T-t_1} +1$  ($q \in \N$) and 
$k_q(t):=\min\{1,q(t_1+\frac{1}{q}-t)_+\}\to \chi_{(0,t_1]}$ in $(0,T)$ as $q\to\infty$:
\begin{equation*}
q\!\int_{t_1}^{t_1+\frac{1}{q}}\!\!   \lla u(\cdot,t),\rho\rra_{\R}\,dt - \lla u_0,\rho\rra_{\R}=\! \!\!\int_0^T\!\!
\left\{\int_{\R} \varphi^*(x,t)\rho'(x)\,dx \!+ \!C_{\varphi}\lla u_s(\cdot,t),\rho'\rra_{\R} \right\}k_q(t)dt
\end{equation*}
Letting $q\to\infty$  
we obtain \eqref{fdp} from \eqref{fi*} and \eqref{eq.punti.lebesgue}. Subtracting from \eqref{fdp} the same inequality with $t_1$ replaced by $t_0$, we obtain \eqref{fdp bis}.
\end{proof}

\smallskip

\noindent {\em Proof of Proposition \ref{condini}.} Let $F^*\subset (0,T)$ be the null set  
given by Lemma  \ref{lefdp}. Let $\{\tau_n\}\subseteq (0,T)\setminus F^*$, $\tau_n\to 0^+$ as $n \to \infty$. Since, by  \eqref{fi*}, $u \in L^{\infty}(0,T;\mathcal{M}^+(\R))$ and  $\varphi^* \in L^{\infty}(0,T;L^1(\R))$, it follows from \eqref{fdp} that  
$\langle u(\cdot,\tau_n),\rho\rangle_{\R} \to \langle u_0,\rho\rangle_{\R}$ for all $\rho\in C^1_c(\R)$.
Since, by Definition \ref{dli}-$(ii)$,  
$\sup_n \|u(\cdot,\tau_n)\|_{\mathcal{M}(\R)} \le C$,
there exist $\mu_0\in\mathcal{M}^+(\R)$ and a subsequence $\{\tau_{n_k}\}$ such that
$u(\cdot,\tau_{n_k})\stackrel{*}\rightharpoonup \mu_0$ 
in $\mathcal{M}(\R)$ as $k\to\infty$.
By standard density arguments, this implies that $\mu_0=u_0$. Hence 
$u(\cdot,\tau_{n})\stackrel{*}\rightharpoonup u_0$  along the whole sequence $\{\tau_n\}$, and 
\eqref{cini} follows from \eqref{fdp} and the arbitrariness of $\{\tau_n\}$. 

Similarly,  it follows from \eqref{fdp bis} that
$\langle u(\cdot,\tau_n),\rho\rangle_{\R}\to \langle u(\cdot,t_0),\rho\rangle_{\R}$ for all $\rho\in C^1_c(\R)$
as $\tau_n\to t_0$ if  $t_0, \tau_n\in (0,T)\setminus F^*$ and we obtain \eqref{cini bis}.  

To prove \eqref{cini tris} we observe that, given $t_0\in [0,T]$ and two sequences $\tau_n^1$ and $\tau_n^2$ contained in 
$(0,T)\setminus F^*$ and converging to $t_0$, $\langle u(\cdot,\tau_n^1)-u(\cdot,\tau_n^2),\rho\rangle_{\R}\to 0$
for all $\rho\in C_c(\R)$. Hence, if $t_0\not\in F^*$, the continuous extension of $u(\cdot,t)$ from $(0,T)\setminus F^*$ with respect to  the weak$^*$ topology is well-defined. 
\hfill$\square$

\smallskip

Let us now prove the results of Subsection \ref{subs33}. As explained there, replacing $x$ by $x-C_\varphi t$ we may assume without loss of generality that $C_\varphi=0$ - namely, it suffices to prove Proposition \ref{dsu}, Theorem \ref{th.reg} and Proposition \ref{wai2}. Moreover, replacing $x$ by $-x$  and $\varphi$ by $-\varphi$, it suffices to do so by assuming that $(H_2)$ is satisfied with $\varphi''<0$, $\varphi'>0$ in $(0,\infty)$ (see Remark \ref{remmi}).  Therefore, we make use of the following assumption:
\begin{equation*}
\left\{\begin{array}{ll} \text{$\varphi\in C^{\infty}([0,\infty))$, $C_\varphi=0$\,;}
\smallskip\\
 \text{$\varphi''(u)<0$, and there exist $H\ge-1$, $K>0$ such that 
 }
\smallskip\\
\text{$\varphi''(u)\,[H\varphi(u)+K] \leq -[\varphi'(u)]^2<0$ for all $u\in[0,\infty)$\,  }
\end{array}\right.\leqno(\tilde H_2)
\end{equation*}
(recall that in this case $\varphi'>0$ and $H\varphi(u)+K>0$ in $[0,\infty)$).

First we  prove some estimates of the constructed entropy solutions.
As already said, these estimates  are analogous to the Aronson-B\'enilan inequality for the convex case $u^p$, $p>1$ (see \cite{AB}).   

\begin{prop}\label{bc} 
Let $(H_1)$ and $(\tilde H_2)$ be satisfied, and let
$u$ be an entropy  solution of problem $(P)$ given by Theorem \ref{th.exi1}. Then  for a.e.~$0< t_1<t_2\leq T$ 
\begin{equation}\label{morex}
\varphi(u_r)(\cdot,t_2)+\frac KH\, \le \, \left(\frac{t_2}{t_1}\right)^H\!\! \left[ \varphi(u_r)(\cdot,t_1)+\frac KH\right]
\quad\text{a.e.~in $\R$\; if  $H\ne 0$}\,,
\end{equation}
\begin{equation}\label{morex'}
\varphi(u_r)(\cdot,t_2)-K\log(t_2)\le \varphi(u_r)(\cdot,t_1)-K\log t_1 \quad\text{a.e.~in $\R$\; if $H=0$\,.}
\end{equation}
Moreover, if there exists $L>0$ such that
\begin{equation*}
H\varphi(u)+K\leq  L(1+u)\varphi'(u)\quad\mbox{for }u\geq 0\,, \leqno(H_3)
\end{equation*}
then $u_t\in \mathcal{M}(\Omega\times (\tau,T))$, $[\varphi(u_r)]_t\in \mathcal{M}(\Omega\times (\tau,T))$, and  
$u\in C((0,T];\mathcal{M}(\Omega))$ for every bounded open set $\Omega\subset\R$ and $\tau>0$.
\end{prop}
\begin{remark} If $\varphi(u)=\sgn p\,[(1+u)^p-1]$ $(p<1, p\ne0)$,   
\eqref{morex} becomes 
\begin{equation*}
u_r(\cdot,t_2)\, \le \, \left(\frac{t_2}{t_1}\right)^{\frac{1}{1-p}}\!\! [ 1+u_r(\cdot,t_1)]-1\quad\text{a.e.~in $\R$, for  a.e.~$0<t_1\le t_2\le T$} 
\end{equation*}
(see Remark \ref{HK}). 
Similarly, if $\varphi(u)={\rm log}(1+u)$, \eqref{morex'} becomes 
\begin{equation*}
u_r(\cdot,t_2)\, \le \, \left(\frac{t_2}{t_1}\right) [ 1+u_r(\cdot,t_1)]-1\quad\text{a.e.~in $\R$, for  a.e.~$0<t_1\le t_2\le T$}\,.
\end{equation*}
\end{remark}

Let $(\tilde{H}_2)$ hold. To prove Proposition \ref{bc} we use a different regularization of $(P_n)$:
$$
\left\{\begin{array}{ll}
y_{nt}^{\ep}+[\varphi(y_n^{\ep})]_x=\ep [\varphi(y_n^{\ep})]_{xx}
& \quad\mbox{in}\  S 
\smallskip\\
y_n^{\ep}=u_{0n}^{\ep} &\quad\mbox{in}\  \R \times \{0\}\,,
\end{array}\right. \leqno{(V_n^{\ep})}
$$ 
 where $\{u_{0n}^{\ep}\}$ satisfies \eqref{est.u0ep}-\eqref{conv.u0ep}. 
The existence, uniqueness and regularity results recalled in Section \ref{sec.ap} for problem $(Q_n^{\ep})$, 
as well as the a priori  estimates in Lemma \ref{L1} and the convergence results in Lemma \ref{prop.conv}-$(i)$, 
continue to hold for solutions of $(V_n^{\ep})$  (see \cite{LSU}). 
In particular, there exist a sequence $\left\{y_n^{\ep_m}\right\}$ and  $y_n\in L^{\infty}(S)\cap L^{\infty}(0,T;L^1(\R))$  
such that  $y_n^{\ep_m}\stackrel{*}\rightharpoonup y_n$ in $L^{\infty}(S)$ and for all $L>0$  
\begin{equation}\label{convyn}
y_n^{\ep_m}\to y_n\quad \text{ in $L^1((-L,L)\times (0,T))$ \; as $\ep_m\to 0$}.
\end{equation}

From $(V_n^{\ep})$, for every   $E$ convex, $F'=E'\varphi'$, and $\zeta$ as in Definition \ref{dsl1}, we get 
$$
\iint_{S}\!\left\{E(y_n^{\ep_m})\zeta_t\!+\!F(y_n^{\ep_m})\zeta_x\right\}dxdt  \!+\!   \int_{\R}\!E(u_{0n}^{\ep})\zeta(x,0)\,dx\ge 
\ep_m\!\iint_{S} \!F'(y_n^{\ep_m}) y_{nx}^{\ep_m}\zeta_xdxdt\,.
$$
Arguing as in the proof of Proposition \ref{sl3} and letting $\ep_m\to 0$, we obtain that 
 $$
  \iint_S \left[E(y_n)\,\zeta_t+F(y_n)\,\zeta_x\right]\,dxdt \ge  - \int_{\R} E(u_{0n})\,\zeta(x,0)\,dx \,.
$$
So $y_n$ satisfies \eqref{eq.KlE} and, by Kru\v zkov's uniqueness theorem, $y_n=u_n$. Hence we have shown: 

\begin{lem}\label{66b}  Let $(H_1)$ and $(\tilde H_2)$ be satisfied, and let $u_n$ be the unique entropy solution of problem 
$(P_n)$ given by Proposition \ref{sl3}. Then there exists a subsequence $\left\{y_n^{\ep_m}\right\}$ of solutions of $(V_n^{\ep})$ 
such that  $y_n^{\ep_m}\stackrel{*}\rightharpoonup u_n$ in $L^{\infty}(S)$ and satisfies \eqref{convyn}.
\end{lem}

\begin{lem}\label{lemarbe}
Let  $(H_1)$ and $(\tilde H_2)$ be satisfied. Then 
\begin{equation}\label{arbe}
\frac{\partial}{\partial t}\left[ \frac{H\varphi(y_n^{\ep})(\cdot,t)+K}{t^H} \right] 
\begin{cases}\le 0 &\text{in $\R$ if $H>0$\,,}\\
\ge 0&\text{in $\R$ if $H<0$\,,}
\end{cases}
\end{equation}
\begin{equation}\label{arbe'}
\frac{\partial}{\partial t}\big[ \varphi(y_n^{\ep})(\cdot,t)-K\log t \big] \le 0
\quad\text{in $\R$\; if $H=0$\,.}
\end{equation}
for all $t\in(0,T)$, $\ep>0$ and $n\in\N$. 
Moreover, if $(H_3)$ is satisfied, then 
\begin{equation}\label{est.ut.ab}
ty_{nt}^{\ep}\le L(1+y_n^{\ep})
\quad\text{in $S$}\,.
\end{equation}
\end{lem}

\begin{proof} 
For convenience we set $A\equiv \ep\frac{\partial^2}{\partial x^2} -\frac{\partial}{\partial x}$, thus 
$y_{nt}^{\ep}=A [\varphi(y_n^{\ep})]$ in  $S$.
Let 
\begin{equation*}%\label{eq.zn}
z_n^{\ep}:=ty_{nt}^{\ep}-g(y_n^{\ep})\,, \quad\text{where }\ g(y_n^{\ep}):=\frac{H\varphi(y^{\ep}_n)+K}{\varphi'(y_n^{\ep})}
\quad (n\in \N).
\end{equation*}
It follows from $(\tilde H_2)$ and a straightforward calculation that 
\begin{equation*}%\label{pde.zn}
z_{nt}^{\ep} \! 
= \! A [\varphi'(y_n^{\ep})z_n^{\ep}]\! +\! [\underbrace{H\! +\! 1\! -\! g'(y_n^{\ep})}_{\le0}] 
\frac{z_n^{\ep}+g(y_n^{\ep})}{t}\! \le\!  
A [\varphi'(y_n^{\ep})z_n^{\ep}]\! +\! [H\! +\! 1\! -\! g'(y_n^{\ep})]\frac{z_n^{\ep}}{t}
\end{equation*}
in $S$. Since $z_n^{\ep}=-g(u_{0n}^{\ep})\le0$ in  $\R \times \{0\}$, it follows from the 
comparison principle for parabolic equations that $z_n^{\ep}\le0$ in $S$ for all $n\in\N$. Hence  
$ty_{nt}^{\ep}(\cdot,t)\le g(y_n^{\ep})(\cdot,t)$ in $\R$ for all $t\in(0,T)$,
which implies \eqref{arbe}, \eqref{arbe'} and, if $(H_3)$ is satisfied, \eqref{est.ut.ab}. 
\end{proof}

\noindent {\em Proof of Proposition \ref{bc}.}  Let $\left\{y_n^{\ep_m}\right\}$ be as  
in the proof of Lemma \ref{66b}. By \eqref{arbe}-\eqref{arbe'},  
\begin{equation*}
\varphi(y_n^{\ep_m})(\cdot,t_2)+\frac KH\, \le \, \left(\frac{t_2}{t_1}\right)^H\!\! \left[ \varphi(y_n^{\ep_m})(\cdot,t_1)+\frac KH\right]
\quad\text{in $\R$\; if $H\ne 0$\,,}
\end{equation*}
\begin{equation*}
\varphi(y^{\ep}_n)(x,t_2)-K\log(t_2)\le \varphi(y^{\ep}_n)(x,t_1)-K\log t_1 \quad\text{in $\R$\; if $H=0$}
\end{equation*}
for all $0<t_1\le t_2\le T$ and $n\in\N$. 
Hence, by Lemma \ref{66b}, 
\begin{equation}\label{more1}
\varphi(u_n)(\cdot,t_2)+\frac KH\, \le \, \left(\frac{t_2}{t_1}\right)^H\!\! \left[ \varphi(u_n)(\cdot,t_1)+\frac KH\right]
\quad\text{a.e.~in $\R$\; if $H\ne 0$\,,}
\end{equation}
\begin{equation}\label{more1'}
\varphi(u_n)(\cdot,t_2)-K\log(t_2)\le \varphi(u_n)(\cdot,t_1)-K\log t_1 \quad\text{a.e.~in $\R$\; if $H=0$}
\end{equation}
for a.e.~$0<t_1\le t_2\le T$. Since $\varphi'$ is strictly decreasing in $[0,\infty)$ (recall that $\varphi$ is concave by assumption $(\tilde{H}_2)$), 
possibly extracting another subsequence (denoted again by $\left\{n_j\right\}$), 
 $\varphi\left(u_{n_j}\right)\to \varphi\left(u_r\right)$ a.e.~in  $S$ (see Remark \ref{come}).
Letting $j\to\infty$ in \eqref{more1}-\eqref{more1'} (with $n=n_j$) we obtain \eqref{morex}-\eqref{morex'}.

Let $\Omega=(-L,L)$. 
If $(H_3)$ is satisfied, it follows from \eqref{est.ut.ab} and \eqref{stil1} that
\begin{equation}\label{est.ut+}
t\int_{\Omega}\left[y_{nt}^{\ep}\right]^+(x,t)\,dx\le L\,|\Omega|+\|u_0\|_{\mathcal{M}(\R)}\quad\text{for all $t\in (0,T]$.}
\end{equation}
 Since $|y_{nt}^{\ep}|=2[y_{nt}^{\ep}]^+-y_{nt}^{\ep}$ a.e.~in $S$, there exists $C_\Omega>0$ such that 
$$\int_{\tau}^{T}\!\!\!\int_{\Omega}\left|y_{nt}^{\ep}\right|(x,t)\,dxdt
\leq 2(T-\tau)\frac{L|\Omega|+\|u_0\|_{\mathcal{M}(\R)}}{\tau}+\int_{\Omega}\left\{y^{\ep}_n(x,\tau)-y_n^{\ep}(x,T)\right\}\leq \frac{C_\Omega}{\tau}
$$
for all $\tau>0$,  $\ep>0$ and $n\in \N$, and, by \eqref{disder},
 \begin{equation}\label{est.utloc}
\int_{\tau}^T\!\!\!\int_{\Omega} |y_{nt}^{\ep}|\,dxdt\leq \frac{C_\Omega}{\tau}\,,\quad \int_{\tau}^T\!\!\!\int_\Omega \left|[\varphi(y_{n}^{\ep})]_t\right|dxdt \leq \frac{MC_\Omega}{\tau}\,.
\end{equation}

Let $\{\ep_m\}$ and $\{n_j\}$ be as in Lemma \ref{66b} and \eqref{conv.un.t}. Then 
$$
\lim_{n_j\to\infty}\,\lim_{\ep_m\to0}\langle y_{n_j}^{\ep_m},\zeta_t\rangle_{\Omega\times(\tau,T)}=
\left\langle u,\zeta_t\right\rangle_{\Omega\times(\tau,T)}\quad \text{for all }\zeta\in C^1_c(\Omega\times(\tau,T))\,,
$$
whence, by \eqref{est.utloc} and the lower semicontinuity of the total variation,
 \begin{equation*}
\|u_t\|_{\mathcal{M}(\Omega\times(\tau,T))}\leq \frac{C_\Omega}{\tau}\,.
\end{equation*}
Similarly, by \eqref{merc.22}, \eqref{lebeq} and Proposition \ref{ide},  
$$
\lim_{n_j\to\infty}\,\lim_{\ep_m\to0}\left\langle \varphi(y_{n_j}^{\ep_m}),\zeta_t\right\rangle_{\Omega\times(\tau,T)}=
\int_{\tau}^T\!\!\!\int_\Omega\varphi(u_r)\,\zeta_t\,dxdt \quad\text{for all $\zeta\in C^1_c(\Omega\times(\tau,T))$},
$$
and, by \eqref{est.utloc} and the lower semicontinuity of the total variation,
 \begin{equation*}
\|[\varphi(u_r)]_t\|_{\mathcal{M}(\Omega\times(\tau,T))}\leq  \frac{MC_\Omega}{\tau}\,.
\end{equation*}

It remains to prove that $u\in C((0,T];\mathcal{M}(\Omega))$. 
Observe that for all $t_1,\,t_2\in(0,T]$, $0<\tau<t_1<t_2$, and $\rho\in C^2_c(\R)$, $0\leq \rho\leq 1$ in $\R$, $\rho=1$ in $\Omega$, 
$$
\begin{aligned}
&\int_{\Omega}|y_n^{\ep}(x,t_2)-y_n^{\ep}(x,t_1)|\,dx\leq \int_{\R}|y_n^{\ep}(x,t_2)-y_n^{\ep}(x,t_1)|\,\rho(x)\,dx\leq \\
& \qquad \leq \int_{t_1}^{t_2}\!\!\!\int_{\R}|y^{\ep}_{nt}|\,\rho\,dxdt=
\int_{t_1}^{t_2}\!\!\!\!\!\int_{\R} (2\,[y^{\ep}_{nt}]^+-y_{nt}^{\ep})\,\rho\,dxdt=\\
&\qquad =2\int_{t_1}^{t_2}\!\!\!\int_{\R}[y^{\ep}_{nt}]^+\,\rho\,dxdt-\int_{t_1}^{t_2}\!\!\!\int_{\R}\{\varphi(y_n^{\ep})\,\rho'(x)+\ep\varphi(y_n^{\ep})\,\rho''(x)\}\,dxdt \le \\
& \qquad\leq  2\,\frac{L\,|\,{\rm supp}\,\rho\,|+\|u_0\|_{\mathcal{M}(\R)}}{\tau}\,(t_2\!-\!t_1)
\!-\!\int_{t_1}^{t_2}\!\!\!\!\int_{\R}\{\,\varphi(y_n^{\ep})\,\rho'+\ep \varphi(y_n^{\ep})\,\rho''\}\,dxdt\,,
\end{aligned}
$$
where we have used  \eqref{est.ut+}.
We let $\ep=\ep_m\to 0$ and use \eqref{disder} and \eqref{stil1lim}:
\begin{eqnarray*}
&& \int_{\Omega}|u_n(x,t_2)-u_n(x,t_1)|\,dx\leq  \\
&&\qquad \leq 2\,\frac{L\,|\,{\rm supp}\,\rho\,|+\|u_0\|_{\mathcal{M}(\R)}}{\tau}\,(t_2-t_1)-\int_{t_1}^{t_2}\!\!\!\int_{\R}\varphi(u_n)\,\rho'(x)\,dxdt\leq \nonumber\\
&&\qquad \leq \left(2\,\frac{L|\supp\rho|\!+\!\|u_0\|_{\mathcal{M}(\R)}}{\tau}\!+\! 
M\|u_0\|_{\mathcal{M}(\R)}\|\rho'\|_{L^{\infty}(\R)}\!\right)\!(t_2\!-\!t_1) =: \frac{\tilde{C}}{\tau}(t_2\!-\!t_1). 
\nonumber
\end{eqnarray*}
By \eqref{conv.un.tt}  and the lower semicontinuity of the total variation, 
\begin{equation*}%\label{cont3}
\|u(\cdot,t_2)-u(\cdot,t_1)\|_{\mathcal{M}(\Omega)}\leq \frac{\tilde{C}}{\tau}|t_1-t_2|\quad\text{for a.e.~$0<\tau<t_1<t_2\leq T$}.
\end{equation*}
So we may define $u(\cdot,t)$ for all $t\in [\tau,T]$ such that $u\in C([\tau,T];\mathcal{M}(\Omega))$. 
Since $\tau>0$ is arbitrary, the  
proof is complete.
\hfill $\square$

\smallskip

To prove Proposition \ref{dsu} we need the following lemma.

\begin{lem}\label{,.-}
Let $(H_1)$ be satisfied, and let $u$ be the solution of problem $(P)$ given by Theorem \ref{th.exi1}.
Let $\{u_{n_j}\}$ be  as in the proof of Theorem~\ref{th.exi1}. 
Then for a.e.~$t\in(0,T)$ and all  $x_0\in {\rm supp}\,u_s(\cdot,t)$ there exist a sequence $\{x_{0k}\}\subset \R$ 
and a subsequence $\{u_{n_k}\}$ of $\{u_{n_j}\}$ such that $x_{0k}\to x_0$ and  $u_{n_k}( x_{0k},t)\to \infty$ as $k\to \infty$.
\end{lem} 
\begin{proof}
Let $x_0\in {\rm supp}\,u_s(\cdot,t)$. 
We may assume that the convergence in \eqref{conv.un.tt} is satisfied for this  
$t$.
Since $x_0\in {\rm supp}\,u_s(\cdot,t)$, there is no neighbourhood $I_{\delta}(x_0)$ such that the sequence $\left\{u_{n_j}(\cdot,t)\right\}$ lies in a bounded subset of $L^{\infty}(I_{\delta}(x_0))$. Otherwise, up to a subsequence, 
$u_{n_j}(\cdot,t)\stackrel{*}\rightharpoonup f_t$ in $L^{\infty}(I_{\delta}(x_0))$ for some $f_t\in L^{\infty}(I_{\delta}(x_0))$, $f_t\geq 0$. However, this would imply that $u_s(\cdot,t)=0$ in $I_{\delta}(x_0)$,  a contradiction.

Setting $\delta=1/k$, we obtain that 
$\sup_{n_j\in\N}\|u_{n_j}(\cdot,t)\|_{L^{\infty}(I_{1/k}(x_0))}=\infty$ for all $k\in\N$.
Hence for all $k\in \N$ there exists $x_{0k}\in I_{1/k}(x_0)$ such that $u_{n_k}(x_{0k},t) \ge k$.
\end{proof}

%\smallskip

\noindent {\em Proof of Proposition \ref{dsu}.} 
As pointed out above, it suffices to prove equality \eqref{crus} by assuming $(\tilde H_2)$. Let $\{u_{n_j}\}$ be  as in the proof of Lemma \ref{,.-}.
By Lemma \ref{66b}, for every $n_j\in\N$ there exists   $\ep_m\to 0$ such that 
\begin{equation}\label{convunj}
y_{n_j}^{\ep_m}(\cdot,t)\to u_{n_j}(\cdot,t)\quad \text{ in $L^1_{{\rm loc}}(\R)$  as $\ep_m\to 0$ for a.e.~$t\in(0,T)$}\,.
\end{equation}
By the proof of Lemma \ref{lemarbe}, for all $t\in(0,T)$
\begin{equation}\label{ax}
\ep_m [\varphi(y_{n_j}^{\ep_m})(\cdot,t)]_{xx} - [\varphi(y_{n_j}^{\ep_m})(\cdot,t)]_x=(y_{n_j}^{\ep_m})_t\le\frac{g(y_{n_j}^{\ep_m})(\cdot,t)}{t}
\quad\text{in $\R$}\,,
\end{equation}
where $g(u)=\frac{H\varphi(u)+K}{\varphi'(u)}>0$.
For every $\underline{x}<\overline{x}$, let $\rho\in C_c^1((\underline{x},\bar{x}))$, $\rho\ge0$. Multiplying \eqref{ax} by $\rho/g(y_{n_j}^{\ep_m}(\cdot,t))$, 
integrating by parts and setting $\Psi(y):= \int_0^y\frac{\varphi'(u)}{g(u)}\,du$, we find that 
\begin{eqnarray*}
&&\int_{\underline{x}}^{\bar{x}} \Psi(y_{n_j}^{\ep_m})(x,t) [\ep_m\rho''(x)+\rho'(x)]\,dx 
\le \frac 1t\int_{\underline{x}}^{\bar{x}}\rho(x)\,dx\,-\\
&-&\ep_m  \int_{\underline{x}}^{\bar{x}} \frac{\varphi'(y_{n_j}^{\ep_m})\,g'(y_{n_j}^{\ep_m})\,[(y_{n_j}^{\ep_m})_x]^2}{[g(y_{n_j}^{\ep_m})]^2}\,(x,t)\,\rho(x)\,dx \le \frac 1t\int_{\underline{x}}^{\bar{x}}\rho(x)\,dx
\end{eqnarray*}
(observe that by $(\tilde H_2)$ there holds $g'(u)\ge H+1\ge0$ and $\Psi$ is bounded). Hence, by \eqref{convunj},
\begin{equation}\label{axy}
\int_{\underline{x}}^{\bar{x}} \Psi(u_{n_j})(x,t)\rho'(x)\,dx \le \frac 1t\int_{\underline{x}}^{\bar{x}}\rho(x)\,dx.
\end{equation}

Let $x_0\in {\rm supp}\,u_s(\cdot,t)$, and let 
$\left\{x_{0k}\right\}\subset \R$,    $\{u_{n_k}\}$ be as in Lemma \ref{,.-}, for a.e.~$t\in (0,T)$. Let   $\bar{x}>x_0$ be fixed. Since $x_{0k}\to x_0$, there exists $\bar{k}\in\N$ such that   $\bar{x}>x_{0k}$ for all $k>\bar{k}$. Consider any sequence $\{\rho_m\}\subset C_c^1((x_{0k},\bar{x}))$, $0\le\rho_m\le1$, $\rho_m\to \chi_{(x_{0k}, \bar{x})}$ 
in $\R$. Without loss of generality, we may assume that both $x_{0k}$ and $\bar{x}$ are  Lebesgue points of $u_{n_k}(\cdot,t)$ for all $k\in\mathbb{N}$. 
Setting $\rho=\rho_m$ and $\underline{x}=x_{0k}$ in \eqref{axy}, letting $m\to\infty$ we find that 
\begin{equation*}
 \Psi(u_{n_k})( x_{0k},t)   \le \Psi(u_{n_k})( \bar{x},t)+  \frac 1t ( \bar{x}-x_{0k}) \quad\text{for all $n_k$}\,.
\end{equation*}
Since $\Psi$ is continuous, by Lemma \ref{,.-} and Remark \ref{come} (recall that $\varphi$ satisfies \eqref{exH_2} since $\varphi$ is strictly concave by assumption $(\tilde{H}_2)$), letting $n_k\to\infty$ gives 
$$
\Psi(u_r)( \bar{x},t) + \frac 1t ( \bar{x}-x_0)\ge \Psi(\infty) \quad\text{ for  a.e.~$\bar{x}>x_0$\,,}
$$
whence by the invertibility of $\Psi$ 
\begin{equation}\label{Psi}
u_r( \bar{x},t)\ge \Psi^{-1}\left(\Psi(\infty)- \frac 1t ( \bar{x}-x_0)\right)  \quad\text{for a.e.~}\bar{x}>x_0\,.
\end{equation}
Letting $ \bar{x}\to x_0^+$ in the previous inequality we obtain \eqref{crus}. 
\hfill $\square$
\medskip

To prove Theorem \ref{th.reg} we need the following result. 
\begin{prop}\label{cc}
 Let $(H_1)$ be satisfied.  Let $C_{\varphi}=0$, and let $u$ be a solution of problem $(P)$. Then  for $a.e.$ $0\le t_1\le t_2\le T$:

\smallskip

\noindent $(i)$ the map $x\mapsto  \Phi(x,t_1,t_2):=\int_{t_1}^{t_2}\varphi(u_r)(x,t)\,dt$ belongs to $BV(\R)$; 

\smallskip

\noindent $(ii)$ for all $x_0,x_1\in\R$, $x_0\le x_1$, 
\begin{equation}\label{dist1}
u(\cdot,t_2)([x_0,x_1])- u(\cdot,t_1)([x_0,x_1]) =\Phi(x_0^-,t_1,t_2)-\Phi(x_1^+,t_1,t_2)\,,
\end{equation}
\begin{equation}\label{dist1'}
u(\cdot,t_2)([x_0,x_1])- u_0([x_0,x_1]) =\Phi(x_0^-,0,t_2)-\Phi(x_1^+,0,t_2)\,.
\end{equation}
\end{prop}

\begin{remark}
It is easily seen that, if $C_{\varphi}\ne0$, equalities \eqref{dist1}-\eqref{dist1'} are replaced by
\begin{equation}\label{dist1g}
u(\cdot,t_2)([x_0,x_1])- \mathcal{T}_{C_{\varphi}(t_2-t_1)} (u(\cdot,t_1))([x_0,x_1]) = 
\Phi(x_0^-,t_1,t_2)-\Phi(x_1^+,t_1,t_2)\,, 
\end{equation}
\begin{equation*}
u(\cdot,t_2)([x_0,x_1])- \mathcal{T}_{C_{\varphi}t_2} (u_0)([x_0,x_1]) = \Phi(x_0^-,0,t_2)-\Phi(x_1^+,0,t_2)\,,
\end{equation*}
where now
\begin{equation}\label{deffg}
\Phi(x,t_1,t_2):=\int_{t_1}^{t_2} 
[\varphi(u_r)-C_\varphi u_r]\big(x+C_{\varphi}(t-t_1),t\big)\,dt\,.
\end{equation}
\end{remark} 
\noindent {\em Proof of Proposition \ref{cc}.}
\noindent $(i)$  By \eqref{disder}, 
$\left|\int_{t_1}^{t_2} \varphi(u_r)(x,t)\,dt\right|\leq M\int_{t_1}^{t_2} u_r(x,t)\,dt \in L^1(\R)$.
We argue as in the proof of Proposition \ref{ac} (see \eqref{us.pl}): 
there exists a null set $N\subset (0,T)$ such that
\begin{equation}\label{us.tot}
\lim_{h\to 0} \frac{1}{h} \int_{\bar{t}}^{\bar{t}+h}\left\langle u(\cdot,t),\rho\right\rangle_{\R}dt =\left\langle u(\cdot,\bar{t}), \rho\right\rangle_{\R}\quad\text{for all $\rho \in C_c(\R)$ and }
\bar{t}\in (0,T)\setminus N\,.
\end{equation}
Let $t_1,t_2\in (0,T)\setminus N$, $0<t_1<t_2<T$, $\rho\in C^1_c(\R)$, and $\zeta(x,t)=g_h(t)\rho(x)$, with $g_h$ as in \eqref{defigh}. 
Since $C_{\varphi}=0$, we obtain from \eqref{ewf} that 
$$
\frac{1}{h}\int_{t_1}^{t_1+h} \!\langle u(\cdot,t),\rho\rangle_{\R}dt-\frac{1}{h}\int_{t_2}^{t_2+h} \!\langle u(\cdot,t),\rho\rangle_{\R}\,dt + \int_0^T \!\!\!\int_{\R}g_h(t)\,\rho'(x)\,\varphi(u_r)(x,t)\,dxdt =0\,.
$$
Letting $h\to 0$, it follows from \eqref{us.tot} that
\begin{eqnarray}\label{equatf}
\lla u(\cdot,t_2), \rho \rra_{\R}-\lla u(\cdot,t_1), \rho \rra_{\R}
=\int_{\R}\Phi(x,t_1,t_2)\rho'(x)\,dx\,. 
\end{eqnarray}
Hence the distributional derivative $\Phi_x(x,t_1,t_2)$ belongs to $\mathcal{M}(\R)$. 

\smallskip

\noindent $(ii)$ We set, for $m\in\N$ and $x\in \R$,
\begin{equation*}
\rho_m(x):= m\left(x\!-\!x_0\! +\!\tfrac 1m\right) \chi_{[x_0 -\frac 1m, x_0]}\!+\! \chi_{(x_0, x_1)}(x)
\!+\!m\left(-\!x\!+\!x_1\! +\!\tfrac 1m\right) \chi_{[x_1,x_1+\frac 1m, ]}.
\end{equation*}
By standard regularization arguments we can choose $\rho=\rho_m$ in \eqref{equatf}: 
\begin{equation}\label{equatf1}
\lla u(\cdot,t_2), \rho_m \rra_{\R}\!-\!\lla u(\cdot,t_1), \rho_m \rra_{\R}
\!=\!m\!\int_{x_0 -\frac 1m}^{x_0}\!\! \Phi(x,t_1,t_2)\,dx\!-\!m\!\int^{x_1+\frac 1m}_{x_1}\!\! \Phi(x,t_1,t_2)\,dx\,. 
\end{equation}
By the Dominated Convergence Theorem,
$\lla u(\cdot,t_i), \rho_m \rra_{\R} \to 
u(\cdot,t_i)([x_0, x_1])$ 
as $m\to\infty$ $(i=1,2)$, whereas, by part $(i)$, 
$$
m\int_{x_0 -\frac 1m}^{x_0}\Phi(x,t_1,t_2)\,dx \to\Phi(x^-_0,t_1,t_2)\,, 
\quad m\int^{x_1+\frac 1m}_{x_1}\Phi(x,t_1,t_2)\,dx \to \Phi(x^+_1,t_1,t_2)\,.
$$
Hence \eqref{dist1} follows from 
\eqref{equatf1}. The proof of \eqref{dist1'} is similar.
\hfill$\square$
\begin{remark}
Observe that, by \eqref{absco2} and \eqref{dist1g} with $x_0=x_1=x$, all entropy solutions of problem $(P)$ satisfy for a.e.~$0\le t_1\le t_2\le T$
\begin{equation*}
\Phi(x^-,t_1,t_2)\le \Phi(x^+,t_1,t_2)\quad \text{for all }x\in\R
\end{equation*}
with $\Phi$ defined by \eqref{deffg}.
\end{remark}

Now we are ready to prove Theorem \ref{th.reg} and Proposition \ref{wai2}. As pointed out at the beginning of this section, in doing so it is not restrictive to assume that $(\tilde H_2)$ holds.

\smallskip

\noindent {\em Proof of Theorem \ref{th.reg}.} 
$(i)$ By \eqref{dist1'}, for a.e.~$0\le t\le T$
\begin{equation*}
u_s(t)(\{x_0\})=u_{0s}(\{x_0\}) 
+\Phi(x_0^-,0,t)-\Phi(x_0^+,0,t)
\ge  u_{0s}(\{x_0\}) - \|\varphi\|_{L^{\infty}(0,\infty)}t, 
\end{equation*}
whence $u_s(t)(\{x_0\})>0$ if $t\in \left(0, \tfrac{u_{0s}(\{x_0\})}{ \|\varphi\|_{L^{\infty}(0,\infty)}}\right)$. Hence \eqref{stib} follows.
\smallskip

\noindent $(ii)$ Let $u_n$ be the entropy solution of problem $(P_n)$ given by Proposition \ref{sl3}. We argue as in the proof of Proposition \ref{cc}: for all $n\in\N$ the map $x\mapsto \Phi_n(x,t_1,t_2):=\int_{t_1}^{t_2}\varphi(u_n)(x,t)\,dt$ belongs to $BV(\R)$, and for a.e.~$0\le t_1\le t_2\le T$ and  a.e.~$x_0\le x_1\in\R$
\begin{equation*}%\label{kpl}
\int_{x_0}^{x_1}u_n(x,t_2)\,dx- \int_{x_0}^{x_1}u_n(x,t_1)\,dx = \Phi_n(x_0^-,t_1,t_2)-\Phi_n(x_1^+,t_1,t_2)\,.
\end{equation*}
Letting $x_1\to\infty$, 
it follows from  \eqref{cmlim} and \eqref{stima.u0n}  that
\begin{equation}\label{funt}
 \int_{t_1}^{t_2}\varphi(u_n)(x,t)\,dt \le \|u_0\|_{\mathcal{M}(\R)} \quad \text{for $n\in\N$ and a.e.~}x\in\R.
\end{equation}

Let $\{y_n^{\ep_m}\}$ be the subsequence used in the proof of Lemma \ref{66b}. By \eqref{arbe} and \eqref{arbe'}, for every $0<t_1\leq t \leq T$ and $x\in\R$ 
$$
\begin{aligned}
\int_{t_1}^{t}\varphi(y^{\ep_m}_n)(x,s)\,ds&=\frac 1H \int_{t_1}^{t}\frac{H\varphi(y^{\ep_m}_n)(x, s)+K}{ s^H}\, s^Hd s-\frac KH (t-t_1)\ge \\
&\ge\frac{ H\varphi(y^{\ep_m}_n)(x,t)+K }{H\,t^H}\;\frac{t^{H+1}-t_1^{H+1}}{H+1}
-\frac KH (t-t_1)\quad \text{if  $H\ne 0$},
\end{aligned}
$$
$$
\begin{aligned}
&\int_{t_1}^t \varphi(y_n^{\ep_m})(x,s)\,ds= \int_{t_1}^t \left[\varphi(y_n^{\ep_m})(x,s)-K\log s\right]\,ds + K \int_{t_1}^t \log s\,ds \geq \\
&\qquad  \geq \left[ \varphi(y_n^{\ep_m})(x,t)-K\log t\right](t-t_1) 
+ K\left[t\log t-t\right]-K\left[t_1\log t_1-t_1\right]\quad \text{if $H= 0$}\,.
\end{aligned}
$$

Letting $\ep_m\to 0$, by  \eqref{funt} we obtain that  for a.e.~$t\in (t_1,T)$ and a.e.~$x\in\R$ 
\begin{equation*}
\|u_0\|_{\mathcal{M}(\R)}\geq \Phi_n(x,t_1,t)\ge 
\begin{cases}
\;\frac{ H\varphi(u_n)(x,t)+K }{H\,t^H}\;\frac{t^{H+1}-t_1^{H+1}}{H+1}-\frac KH (t-t_1) &\text{if }H\ne 0\,,\smallskip\\
\;\left[\varphi(u_n)(x,t)-K\right]\,(t-t_1)  + K\,t_1\log\frac{t}{t_1} &\text{if }H=0\,.
\end{cases}
\end{equation*}
Letting $t_1\to 0^+$ we find  in both cases that
\begin{equation}\label{f-1}
\varphi(u_n)(x,t) \le \frac{(H+1)  \|u_0\|_{\mathcal{M}(\R)}}{t} +K\quad \text{for a.e.~$t\in (t_1,T)$ and a.e.~$x\in\R$}\,
\end{equation}
(recall that we have assumed $H>-1$ if $\varphi$ is bounded; otherwise, if $\varphi$ is unbounded, there holds $H\ge 0$ since $\varphi'>0$ and $H\varphi+K>0$ in $[0,\infty)$ by $(\tilde H_2)$).
If $\lim_{u\to\infty}\varphi(u)=:\gamma<\infty$,  $K<\gamma$ and $H>-1$,   the sequence $\{u_n(\cdot,t) \}$ lies in a bounded subset of $L^{\infty}(\R)$ (thus, by \eqref{conv.un.tt} $u_s(\cdot,t)=0$, and $u_r(\cdot,t)\in L^{\infty}(\R)$) for  $a.e.$ $t\in(0,T)$ such that
$$
\frac{(H+1)\,\|u_0\|_{\mathcal{M}(\R)}}{t} +K < \gamma \quad \Leftrightarrow \quad
t>\frac{(H+1)\,\|u_0\|_{\mathcal{M}(\R)}}{\gamma-K}\,.
$$
This proves claim $(ii)$-$(a)$.

If $\gamma=\infty$,  there holds $H\ge 0$ since $H\varphi+K>0$ in $[0,\infty)$ (see $(\tilde H_2)$). Then by \eqref{f-1} the sequence $\{u_n(\cdot,t) \}$ lies in a bounded subset of $L^{\infty}(\R)$ for  $a.e.$ $t\in(0,T)$, hence by \eqref{conv.un.tt} as $n\to\infty$ we obtain that $t_0=0$. Hence claim $(ii)$-$(b)$ follows. This completes the proof. 
\hfill$\square$

\begin{remark}\label{phi_k bis} 
As we claimed in Remark \ref{phi_k}, in Theorem \ref{th.reg}-$(ii)$ we may relax hypothesis $(H_2)$ to $(H_{2,k})$, with
$k>0$. To prove this, for every $u_0\in\mathcal{M}^+(\Omega)$ let $\{u_{0n}\}$ be any sequence as in \eqref{stima.u0n}-\eqref{conv.qo.u0n}, and let $u_n$ be the entropy solution of problem $(P_n)$. Set $v_{0n}:= G_k(u_{0n})$, where $G_k(u):=(u-k)^+$ for every $u\geq 0$, and let $v_n$ be the entropy solution of problem
$$\left\{\begin{array}{ll}
v_{nt}+[\varphi_k(v_n)]_x=0&\mbox{in}\ \,S\\
v_n=v_{0n}&\mbox{in}\ \,\R\times\{0\}\,
\end{array}\right.$$
($\varphi_k(u)=\varphi(u+k)-\varphi(k)$).
A standard calculation shows that $G_k(u_n)$ is an entropy subsolution of the above problem, whence 
\begin{equation}\label{est.Gkun}
\mbox{$G_k(u_n)\leq v_n$ $a.e.$ in $S$.}
\end{equation}  
Following the proof of Theorem \ref{th.exi1}, the sequence $\{v_n\}$ converges to an entropy solution $v$ of problem $(P)$ with initial datum $v_0=u_{0s}+G_k(u_{0r})$. Moreover, by assumption  $(H_{2,k})$, $\varphi_k$ satisfies $(H_2)$ and we may apply Theorem \ref{th.reg}-$(ii)$ to $v$. Therefore the conclusion follows from \eqref{est.Gkun}.\end{remark}

%\smallskip

\noindent {\em Proof of Proposition \ref{wai2}.} By the proof  of Proposition \ref{dsu},
 inequality \eqref{Psi} is satisfied  for a.e.~$t\in (0,T)$ and all $x_0\in \supp u_s(\cdot,t)$. 
We fix such $t$.
 Let $x_1\in \supp u_s(\cdot,t)$ and set 
$\mathcal I_1:=(x_1-\ep,x_1+\ep)$ with $\ep>0$. By \eqref{Psi}, 
$$
\int_{\mathcal I_1}u_r(x,t)\,dx\!\ge\!\! \int_{x_1}^{x_1+\ep}\! \!  \Psi^{-1}\left(\! \Psi(\infty)\! - \! \frac 1t (x\! -\! x_1)\right)\,dx
\! =\! \! \int_{0}^{\ep}\! \!  \Psi^{-1}\left(\Psi(\infty)\! -\!  \frac yt\right)dy=:B_\ep\,.
$$
If $\supp u_s(\cdot,t)\not\subset \mathcal I_1$, let $x_2\in \supp u_s(\cdot,t)\setminus \mathcal I_1$ and set 
$\mathcal I_2:=(x_2-\ep,x_2+\ep)$. Since $(x_1,x_1+\ep)\cap(x_2,x_2+\ep)=\emptyset$ 
we have that 
$$
\int_{\mathcal I_1\cup\mathcal I_2}u_r(x,t)\,dx \ge 
\int_{x_1}^{x_1+\ep} u_r(x,t)\,dx+\int_{x_2}^{x_2+\ep}u_r(x,t)\,dx \ge 2B_\ep\,.
$$
We continue this construction recursively as long as $\supp u_s(\cdot,t)\not\subset \mathcal I_1\cup\cdots\cup\mathcal I_{n-1}$, with $\mathcal I_{n-1}:=(x_{n-1}-\ep,x_{n-1}+\ep)$:
there exists $x_n\in \supp u_s(\cdot,t)\setminus \{\mathcal I_1\cup\cdots\cup\mathcal I_{n-1}\}$ such that, setting
$\mathcal I_n:=(x_n-\ep,x_n+\ep)$, 
$$
nB_\ep\le \int_{\mathcal I_1\cup\cdots\cup\mathcal I_n }u_r(x,t)\,dx\le \|u_0\|_{\mathcal M(\R)}\,.
$$
Hence this construction stops at some $n=n_\ep$, and $n_\ep B_\ep\le \|u_0\|_{\mathcal M(\R)}$. Therefore,
$$
\supp u_s(\cdot,t)\subset \mathcal I_1\cup\cdots\cup\mathcal I_{n_\ep},
\qquad| \supp _s(\cdot,t)|\le |\mathcal I_1\cup\cdots\cup\mathcal I_{n_\ep}|\le 2n_\ep \ep\le \frac {2\ep}{B_\ep}\|u_0\|_{\mathcal M(\R)}.
$$
Since $B_\ep/\ep\to \infty$ as $\ep\to 0$, the claim follows. 
\hfill$\square$

%%%%%%%%%%%%%%%%%%%%%%%%%%%%%%%%%%%%%%%%%%%%%%%%%%%%%%%%%%%%%%%%%%%%%%%%%%%%%%%%%%%%%%%%%%%%%%%%%%%%%%%%%%%%%%%%%%

\section{Uniqueness: Proofs}\label{sec.uni}
\setcounter{equation}{0}

Again, without loss of generality we may assume that $C_\varphi=0$ in the following proofs (see Remark \ref{trasl}). 
\smallskip

\noindent {\em Proof of Proposition \ref{preco}.} 
$(i)$ The first step of  the proof consists in showing that
\begin{equation}\label{cur}
{\rm ess}\lim_{t\to 0^+}\|u_r(\cdot,t) - u_{0r}\|_{L^1(\R)}=0\,.
\end{equation}
Let $\{u_n^{\ep}\}$ be the sequence of solutions to problems $(Q_n^{\ep})$ considered in Section \ref{sec.ap}, and let $\{x_l\}$ $(l=1,\dots,N)$ be as in \eqref{assu0}.
We set $I_l:=(x_l,x_{l+1})$, $Q_l:=I_l\times(0,\tau)$ $(l=1,\dots,N-1)$, $I_-:=(-\infty,x_1)$, $I_+:=(x_N,\infty)$, 
and $Q_{\pm}:=I_{\pm}\times(0,\tau)$. 

Let $1\le l\le N-1$  and $\rho\in C^2_c(I_l)$, $\rho\geq 0$. 
Let $h_0>0$ be such that $x+h\in I_l$ if $x\in\,$supp$\,\rho$ and $|h|<h_0$. 
Let $\delta>0$. Setting $v_n^\ep(x,t):=u_n^\ep(x+h,t)$ and $z:=(u_n^\ep-v_n^\ep)(\rho+\delta)$, we apply 
the $L^1$-contraction property to the parabolic equation 
$$
z_t\!+\!\left[\left(R\!+\!\frac{2\ep \rho'}{\rho+\delta}\right)z\right]_x\!\!-\ep z_{xx}
=\left(\frac {R\rho'}{\rho+\delta}\!+\!\frac{\ep\rho''}{\rho+\delta}\right)z
=(\varphi_\ep(u_n^\ep)-\varphi_\ep(v_n^\ep))\rho'+\ep[u_n^\ep-v_n^\ep]\rho'',
$$
where $R:=\frac {\varphi_\ep(u_n^\ep)-\varphi_\ep(v_n^\ep)}{u_n^\ep-v_n^\ep}$ if $u_n^\ep\ne v_n^\ep$, 
and $R:=\varphi_\ep'(u_n^\ep)$ otherwise. Hence
$$
\begin{aligned}
&\int_{I_l} |z(x,\tau)|dx \leq \! \int_{I_l} |z(x,0)| dx \!
+\!\int_0^{\tau}\!\!\!\int_{I_l} |\varphi_{\ep}(u_n^{\ep}(x,t))\!-\!\varphi_{\ep}(u_n^{\ep}(x+h,t))||\rho'(x)|dxdt+\\
&\qquad \qquad \qquad +\,\ep\int_0^{\tau}\!\!\!\int_{I_l} |u_n^{\ep}(x,t))-u_n^{\ep}(x+h,t))|\,|\rho''(x)|\,dxdt
\quad \text{for $\tau\in (0,T)$}.
\end{aligned}
$$
First we let $\delta\to 0$ and then  $\ep=\ep_m\to 0$, where $\{\ep_m\}$ is as 
in Lemma \ref{prop.conv}. Hence  
\begin{eqnarray}\label{ec3}
&&\int_{I_l} |u_n(x,\tau)-u_n(x+h,\tau)|\,\rho(x)\,dx \leq  \int_{I_l} |u_{0n}(x)-u_{0n}(x+h)|\,\rho(x)\,dx + \\
&&\qquad +\int_0^{\tau}\!\!\!\int_{I_l} |\varphi(u_n(x,t))-\varphi(u_n(x+h,t))|\,|\rho'(x)|\,dxdt\quad\text{for a.e.~$\tau\in (0,T)$,}\nonumber
\end{eqnarray}
where $u_n$ is the entropy solution of problem $(P_n)$ $(n\in\N)$. 
Since, by \eqref{conv.un.tt},  
$u_n(\cdot,t)\stackrel{*}\rightharpoonup u(\cdot,t)$ in $\mathcal{M}(I_l)$  for a.e.~$t\in(0,T)$ and, by \eqref{absco22} and \eqref{assu0}, 
$u_s(\cdot,t)\lefthalfcup I_l\leq u_{0s}\lefthalfcup I_l=0$, 
the lower semicontinuity of the total variation implies that  for a.e.~$\tau\in (0,T)$ 
\begin{equation*}
\int_{I_l} |u_r(x,\tau)-u_r(x+h,\tau)|\,\rho(x)\,dx \leq \liminf_{n\to\infty}\int_{I_l} |u_n(x,\tau)-u_n(x+h,\tau)|\,\rho(x)\,dx\,.
\end{equation*}
By \eqref{conv.qo.u0n}, $\int_{I_l} |u_{0n}(x)-u_{0n}(x+h)|\,\rho(x)\,dx \to \int_{I_l} |u_{0r}(x)-u_{0r}(x+h)|\,\rho(x)\,dx$. In addition $\varphi(u_{n_j})\to \varphi(u_r)$ in $L^1(Q_l)$ for a subsequence $\{u_{n_j}\}$ of $\{u_n\}$ 
 (see Remark \ref{come}). Letting $n=n_j\to\infty$ in \eqref{ec3}, we obtain that  for a.e.~$\tau\in (0,T)$ 
\begin{eqnarray}\label{ec4}
&&\int_{I_l} |u_r(x,\tau)-u_r(x+h,\tau)|\,\rho(x)\,dx \leq \int_{I_l} |u_{0r}(x)-u_{0r}(x+h)|\,\rho(x)\,dx + \\
&&\qquad \qquad \quad+\int_0^{\tau}\!\!\!\int_{I_l}|\varphi(u_r(x,t))-\varphi(u_r(x+h,t))|\,|\rho'(x)|\,dxdt\,. \nonumber
\end{eqnarray}

Let $\{\tau_n\}\subset(0,T)$ be any sequence such that $\tau_n\to 0^+$ and \eqref{ec4} is satisfied with $\tau=\tau_n$. Since $u_{0r}\in L^1(\R)$ and $\varphi(u_r)\in L^1(S)$, it follows from \eqref{ec4} and the Fr\' echet-Kolmogorov Theorem that the sequence $\{u_r(\cdot,\tau_n)\,\rho\}$ is relatively compact in $L^1(\R)$. Then, by \eqref{cini}  and  a standard argument, 
\begin{equation}\label{ec5}
u_r(\cdot,\tau_n)\,\rho  \to u_{0r}\,\rho \quad\mbox{in}\ \,L^1(\R)\,.
\end{equation}

It follows from \eqref{ewf} and \eqref{cini tris} that for each $n\in\N$
\begin{equation}\label{fdpn}
\int_{I_l} [u_r(x,\tau_n)- u_{0r}(x)]\,\rho(x)\,dx=\int_0^{\tau_n}\!\!\! \int_{I_l}  \varphi(u_r)(x,t)\rho'(x)\,dxdt \,.
\end{equation} 
For sufficiently small $\delta>0$, the characteristic function $ \chi_{(x_l,x_l+\delta)\cup(x_{l+1}-\delta,x_{l+1})}$
can be approximated by functions $\rho_k\in C^2_c(I_l)$, $\rho_k\geq 0$ such that $\int_{I_l}|\rho_k'(x) |dx\le 4$
for all $k\in \N$. Setting $\rho=\rho_k$ in \eqref{fdpn} and letting $k\to\infty$, we find that
\begin{eqnarray}\label{est.resi}
&&\int_{x_l}^{x_l+\delta} u_r(x,\tau_n)\,dx + \int_{x_{l+1}-\delta}^{x_{l+1}} u_r(x,\tau_n)\,dx
\leq \\&&\qquad 
\le  \int_{x_l}^{x_l+\delta}u_{0r}(x)\,dx + \int_{x_{l+1}-\delta}^{x_{l+1}} u_{0r}(x)\,dx +4\|\varphi\|_{L^{\infty}(0,\infty)} \tau_n  \,.\nonumber
\end{eqnarray}
Since $u_{0r}\in L^1(\R)$, for every $\sigma>0$ there exists $\delta>0$ such that
\begin{equation}\label{uff}
\int_{x_l}^{x_l+\delta}u_{0r}(x)\,dx + \int_{x_{l+1}-\delta}^{x_{l+1}} u_{0r}(x)\,dx\leq \sigma\,.
\end{equation}
If $\rho\in C_c(I_l)$ is such that $0\le\rho\le 1$ in $I_l$, $\rho=1$ in $[x_l+\delta,x_{l+1}-\delta]$, then 
\begin{eqnarray*}
&&|u_r(\cdot,\tau_n)- u_{0r}| =|u_r(\cdot,\tau_n)- u_{0r}|\rho +
|u_r(\cdot,\tau_n)- u_{0r}| (1-\rho)\chi_{(x_l,x_l+\delta)\cup(x_{l+1}-\delta,x_{l+1})} 
\end{eqnarray*}
in $I_l$. Hence,   
by \eqref{est.resi} and \eqref{uff},
\begin{eqnarray*}
&&\int_{I_l}|u_r(\cdot,\tau_n)- u_{0r}|dx \le  2 \left\{ \int_{x_l}^{x_l+\delta}u_{0r}dx + \int_{x_{l+1}-\delta}^{x_{l+1}} u_{0r}dx\right\} 
+4\|\varphi\|_{L^{\infty}(0,\infty)} \tau_n+ \\ 
&&\qquad  +\int_{I_l}|u_r(\cdot,\tau_n)- u_{0r}| \rho\, dx
\le \int_{I_l}|u_r(\cdot,\tau_n)- u_{0r}| \rho\, dx+4\|\varphi\|_{L^{\infty}(0,\infty)} \tau_n+2\sigma\,. \nonumber
\end{eqnarray*}
Letting $n\to\infty$ in the above inequality, by \eqref{ec5} we obtain that
$\limsup\limits_{n\to\infty}\int_{I_l}|u_r(\cdot,\tau_n)-u_{0r}|\,dx \leq 2\sigma$,
whence, by the arbitrariness of $\sigma$, 
$$
\lim_{n\to \infty}  \int_{I_l}|u_r(x,\tau_n)-u_{0r}(x)|\,dx=0 \qquad (l=1,\dots,N-1)\,.
$$
A similar argument shows that $\int_{I_{\pm}}|u_r(x,\tau_n)-u_{0r}(x)|\,dx\to 0$ as $n\to\infty$, thus \eqref{cur} follows.

To complete the proof of  \eqref{epc.tot} observe that by \eqref{absco22} there holds $u_s(\cdot,t)\leq u_{0s}$ in $\mathcal{M}(\R)$ (recall that $C_{\varphi}=0$ by assumption). Hence
$$ 
\left\langle u_{0s}-u_s(\cdot,t), \rho\right\rangle_{\R}
\geq \|u_s(\cdot,t)-u_{0s}\|_{\mathcal{M}(\R)}
$$ 
for all $\rho \in C_c(\R)$ such that $\rho(x)=1$ for every $x\in {\rm supp}\,u_{0s}$. From the previous inequality,  
\eqref{cini} and \eqref{cur} we get
\begin{eqnarray}\label{ec8}
&&{\rm ess}\lim_{t\to 0^+}  \|u_s(\cdot,t)-u_{0s}\|_{\mathcal{M}(\R)}\leq {\rm ess}\lim_{t\to 0^+} \left\langle u_{0s}-u_s(\cdot,t), \rho\right\rangle_{\R}= \\
&& = {\rm ess}\lim_{t\to 0^+} \left\{\left\langle u_{0}-u(\cdot,t), \rho\right\rangle_{\R}-\int_{\R}(u_r(x,t)-u_{0r})\,\rho(x)\,dx\right\}=0\,.\nonumber
\end{eqnarray}
Then \eqref{epc.tot} follows.

\smallskip

\noindent $(ii)$ Let $\zeta^{\pm}\in C^1_c(Q_{\pm})$, $\zeta^{\pm}\ge0$, and for every $1\le l\le N-1$ let $\zeta_l\in C^1_c(Q_l)$, $\zeta_l\geq 0$. Let $h_0>0$ be such that $(x+h,t)\in Q_l$ (respectively, $(x+h,t)\in Q_{\pm}$) if $(x,t)\in\,$supp$\,\zeta_l$ (respectively, if $(x,t)\in\,$supp$\,\zeta_{\pm}$) and $|h|<h_0$. 

Let $u$ be an entropy solution of problem $(P)$, thus $v(\cdot,t)=\mathcal{T}_{-h}(u(\cdot,t))$ is an entropy solution of problem $(P)$ with $u_0$ replaced by $v_0:=\mathcal{T}_{-h}(u_0)$ (see Remark \ref{trasl}). We shall prove that for all $l=1,\dots,N-1$ and $\zeta_l$ as above
\begin{equation}\label{contra5CONT}
\iint_{Q_l} \{|v_{r}-u_{r}|\,\zeta_{lt}+  \sgn(v_r-u_r)\left[\varphi(v_{r})-\varphi(u_{r}) \right]\, \zeta_{lx}\}\,dxdt \geq 0 \,, 
\end{equation}
and for all $\zeta^{\pm}$ as above 
\begin{equation}\label{contra6CONT}
\iint_{Q_\pm}  \{|v_{r}-u_{r}|\, \zeta_t^\pm + \sgn(v_r-u_r)\left[\varphi(v_{r})-\varphi(u_{r})\right]\, \zeta_x^\pm\}\, dxdt \geq 0\,.
\end{equation}

Relying on \eqref{contra5CONT}-\eqref{contra6CONT} we can conclude the proof by an argument similar to that used in $(i)$. Let $\rho \in C^1_c(I_l)$, $0\leq \rho\leq 1$, be such that $x+h\in I_l$ if $x\in\supp \rho$ and $|h|<h_0$. By a proper choice of the function $\zeta_l$ in \eqref{contra5CONT}, for $a.e.$ $0<t_0<t_1\leq T$ we get
\begin{eqnarray*}
&&  \int_{I_l} |u_r(x,t_1)-v_r(x,t_1)|\,\rho(x)\,dx \leq \int_{I_l} |u_{r}(x,t_0)-v_{r}(x,t_0)|\,\rho(x)\,dx +\\
&&+ \int_{t_0}^{t_1} |\varphi(u_r)-\varphi(v_r)|\,|\rho'(x)|\,dxd\,t. 
\end{eqnarray*}
Let $t_0>0$ be fixed. Then for every $\tau\in (t_0,T]$ there exists a sequence $\tau_n\to \tau$ such that $\tau_n\in (t_0,T]$ and the above inequality holds true with $t_1=\tau_n$ for every $n$: 
\begin{eqnarray}\label{eq.cont>0.1}
 &&  \int_{I_l} |u_r(x+h,\tau_n)-u_r(x,\tau_n)|\,\rho(x)\,dx \leq \int_{I_l} |u_{r}(x+h,t_0)-u_{r}(x,t_0)|\,\rho(x)\,dx +\\
&&+ \|\rho'\|_{\infty} \int_{0}^{T} |\varphi(u_r(x+h,t))-\varphi(u_r(x,t))|\,dxdt\,.\nonumber
\end{eqnarray}
Since $\varphi(u_r)\in L^1(S)$ and $u_r(\cdot,t_0)\in L^1(\R)$, inequality \eqref{eq.cont>0.1} and the Fr\' echet-Kolmogorov Theorem imply that the sequence 
$\{u_r(\cdot,\tau_n)\,\rho\}$ is relatively compact in 
$L^1(\R)$, whence, by Proposition \ref{condini}
and  a standard argument, 
\begin{equation}\label{ec5CONT}
u_r(\cdot,\tau_n)\,\rho  \to u_{r}(\cdot,\tau)\,\rho \quad\mbox{in}\ \,L^1(\R)\,.
\end{equation}
Moreover, by arguing as in \eqref{est.resi} and \eqref{uff} with $u_{0r}$ replaced by $u_r(\cdot,\tau)$, for every $\sigma>0$ there exists $\delta>0$ such that  
\begin{equation}\label{est.resiCONT}
\int_{x_l}^{x_l+\delta} u_r(x,\tau_n)\,dx + \int_{x_{l+1}-\delta}^{x_{l+1}} u_r(x,\tau_n)\,dx
\leq  \sigma +4\|\varphi\|_{L^{\infty}(0,\infty)} |\tau_n-\tau|  \,.
\end{equation}
As in the proof of claim $(i)$, combining \eqref{ec5CONT} and \eqref{est.resiCONT} gives 
$$
\lim_{n\to \infty}  \int_{I_l}|u_r(x,\tau_n)-u_{r}(x,\tau)|\,dx=0 \qquad (l=1,\dots,N-1)
$$
(by a similar argument, 
$\int_{I_{\pm}}|u_r(x,\tau_n)-u_{r}(x,\tau)|\,dx\to 0$ as $n\to\infty$), whence 
\begin{equation*}
{\rm ess}\lim_{t\to \tau} \|u_r(\cdot,t)-u_r(\cdot,\tau)\|_{L^1(\R)}=0\,.
\end{equation*}
Since $C_{\varphi}=0$ it follows from \eqref{absco2} that $u_s(\cdot,t_2)\leq u_{s}(\cdot,t_1)$ in $\mathcal{M}(\R)$ if $t_2> t_1$, whence by arguing as in \eqref{ec8} we also obtain
$${\rm ess}\lim_{t\to \tau^+}  \|u_s(\cdot,t)-u_{s}(\cdot,\tau)\|_{\mathcal{M}(\R)}={\rm ess}\lim_{t\to \tau^-}  \|u_s(\cdot,t)-u_{s}(\cdot,\tau)\|_{\mathcal{M}(\R)}=0$$
and claim $(ii)$ follows.

Finally, it remains  to prove \eqref{contra5CONT} (the proof of \eqref{contra6CONT} is analogous).  Let $1\le l\le N-1$  and $\zeta_l\in C^1_c(Q_l)$, $\zeta_l\geq 0$, be fixed as above. 
Since $C_{\varphi}=0$, it follows from \eqref{absco22} and \eqref{assu0} that 
$u_{s}(\cdot,t)=v_{s}(\cdot,t)=0$ on supp$\,\zeta_l(\cdot,t)$ for a.e.~$t\in (0,T)$, and from
 \eqref{miseikru} that,
for $k\in [0,\infty)$, 
\begin{eqnarray}\label{contra2CONT}
&& 
\iint_{Q_l} \big\{|u_r-k|\,\zeta_{lt}+\sgn(u_r-k)\left [\varphi(u_r)-\varphi(k)\right ]\zeta_{lx}\big\}\,dxdt \geq 0
\,,  
\end{eqnarray}
\begin{eqnarray}\label{contra2vCONT}
&& 
\iint_{Q_l} \big\{|v_r-k|\,\zeta_{lt}+\sgn(v_r-k)\left [\varphi(v_r)-\varphi(k)\right ]\zeta_{lx}\big\}\,dxdt \geq 0
\,.  
\end{eqnarray}
We apply Kru\v zkov's method of doubling variables. Let $Z_l=Z_l(x,t,y,s)\in C^1_c(Q_l\times Q_l)$, $Z_l\ge 0$.
It follows from \eqref{contra2CONT}-\eqref{contra2vCONT} that
\begin{eqnarray*}
&& \iint_{Q_l} \big\{|u_{r}(x,t)-v_{r}(y,s)|\,Z_{lt}(x,t,y,s)+ \\
&&\qquad +\sgn(u_r(x,t)-v_r(y,s))\left [\varphi(u_r(x,t))-\varphi(v_r(y,s))\right ]Z_{lx}(x,t,y,s)\big\}\,dxdt \geq 0
\nonumber 
\end{eqnarray*}
 and 
\begin{eqnarray*}
&& \iint_{Q_l} \big\{|v_{r}(y,s)-u_{r}(x,t)|\,Z_{ls}(x,t,y,s)+ \\
&&\qquad + \sgn(v_r(y,s)-u_r(x,t))\left [\varphi(v_r(y,s))-\varphi(u_r(x,t))\right ]Z_{ly}(x,t,y,s)\big\}\,dyds 
\geq 0\,,
\nonumber  
\end{eqnarray*} 
whence 
$$
\begin{aligned}
&
\iint\!\!\!\!\!\!\iint_{Q_l\times Q_l} \big\{|u_{r}(x,t)-v_{r}(y,s)|\,\left(Z_{lt}+Z_{ls}\right)(x,t,y,s)+ \\
&+\!\sgn\!(\!u_r(x,t)\!-\!v_r(y,s)\!)[\varphi(u_{r}\!)(x,t)\!-\!\varphi(v_{r}\!)(y,s)\!](Z_{lx}\!+\!Z_{ly})(x,t,y,s)\big\}\,dxdtdyds\geq 0. 
\end{aligned}
$$
We choose 
$$
Z_l(x,t,y,s)=Z_l^{\epsilon}(x,t,y,s):=\zeta_l(x,t)\,\zeta_{\epsilon}(x-y,t-s) \qquad(\epsilon>0)\,,
$$
where 
$\zeta_{\epsilon}$ is a smooth approximation of the Dirac mass $\delta_{(0,0)}$,
\begin{equation*}
\zeta_{\epsilon}(x,y)=\frac{1}{\epsilon^2}\theta\left( \frac{x}{\epsilon}\right)\,\eta\left( \frac{t}{\epsilon}\right)\ge0\,, \quad \text{ with ${\rm supp}\,\theta\subseteq(-1,1)$, ${\rm supp}\,\eta\subseteq(-1,1)$.}
\end{equation*}
Then
$Z_{lt}+Z_{ls}=\zeta_{lt}\,\zeta_{\epsilon}$ and $Z_{lx}+Z_{ly}=\zeta_{lx}\,\zeta_{\epsilon}$, 
 whence, for sufficiently small $\epsilon$, 
$$
\begin{aligned}
&
\iint\!\!\!\!\!\!\iint_{Q_l\times Q_l} \big\{|u_{r}(x,t)-v_{r}(y,s)|\,\zeta_{lt}(x,t)+ \\
&+\!\sgn\!(\!u_r(x,t)\!-\!v_r(y,s)\!)[\varphi(u_{r}\!)(x,t)\!-\!\varphi(v_{r}\!)(y,s)\!]\zeta_{lx}(x,t)\!\big\}\zeta_{\epsilon}(x\!-\!y,t\!-\!s)\,dxdtdyds \geq 0 .
\end{aligned}
$$
Now \eqref{contra5CONT} follows by letting $\epsilon\to 0^+$: we claim that
\begin{eqnarray}\label{conv.lg1CONT}
&&\lim_{\epsilon\to 0^+}\iint\!\!\!\!\!\!\iint_{Q_l\times Q_l} |u_{r}(x,t)-v_{r}(y,s)|\,\zeta_{lt}(x,t)\,\zeta_{\epsilon_n}(x-y,t-s)\,dxdtdyds= \\
&&\qquad = \iint_{Q_l} |u_{r}(x,t)-v_{r}(x,t)|\,\zeta_{lt}(x,t)\,dxdt\,.\nonumber
\end{eqnarray} 
Analogously, it can be proven that
\begin{equation*}
\lim_{\epsilon\to 0^+} \iint\!\!\!\!\!\!\iint_{Q_l\times Q_l} \sgn(u_r(x,t)-v_r(y,s))[\varphi(u_{r})(x,t)-\varphi(v_{r})(y,s) ]\zeta_{lx}(x,t)\zeta_{\epsilon_n}(x-y,t-s)\,dxdtdyds=
\end{equation*}
\begin{equation*}
=\iint_{Q_l} \sgn(u_r(x,t)-v_r(y,s)) \left [\varphi(u_{r})(x,t)-\varphi(v_{r})(x,t)\right ]\zeta_{lx}(x,t)\,dxdt\,.
\end{equation*}
In order to prove \eqref{conv.lg1CONT}, for every sequence $\{\epsilon_n\}$, $\epsilon_n\to 0$, we set 
$$
F_{n}(x,t)\!:=\! \iint_{Q_l} \! |u_{r}(x,t)-v_{r}(y,s)| \zeta_{\epsilon_n}(x-y,t-s)\,dyds\quad \text{for }(x,t)\in K_l:={\rm supp}\,\zeta_l,
$$
and observe that $F_n\to |u_{r}-v_{r}|$ a.e.~in $(x,t)\in K_l$ and  
$$
\begin{aligned}
|F_n(x,t)|
&\leq |u_{r}(x,t)| + \iint_{Q_l} |v_{r}(y,s)|\zeta_{\epsilon_n}(x-y,t-s)\,dyds\\ 
&=|u_{r}(x,t)| + (\zeta_{\epsilon_n} * |v_{r}|)(x,t) \to |u_{r}(x,t)| +  |v_{r}(x,t)|\quad \mbox{in}\ \,L^1(K_l)\,.
\end{aligned}
$$
Thus, by a variant of the Dominated Convergence Theorem ($e.g.$, see \cite[Theorem 4, Section 1.3]{EG}), 
$F_n\to |u_{r}-v_{r}|$ in $L^1(K_l)$, and we obtain \eqref{conv.lg1CONT}. This completes the proof of \eqref{contra5CONT}, thus the result follows. 
\hfill $\square$ 

\medskip

\noindent {\em Proof of Theorem \ref{uni}.} 
Without loss of generality we may assume that $\varphi$ is nondecreasing (see Remark \ref{remmi}). 
By Theorem \ref{th.reg}-$(i)$,
$$
\tau:= \sup\, \{t\in[0,T) \,|\,  u_{is}(\cdot,t)(\{x_l\})>0 \;\; \forall l=1,\dots,N;\;i=1,2\}>0\,.
$$ 
Let us first prove that
\begin{equation}\label{contra}
\text{ $u_{1r}=u_{2r}$ \quad a.e.~in $\R\times(0,\tau)$\,.}
\end{equation}
To this aim, let $x_1,\dots,x_N$ be the points in \eqref{assu0}. Set $I_l:=(x_l,x_{l+1})$, $Q_l:=I_l\times (0,\tau)$ $(l=1,\dots,N-1)$, and $I_-:=(-\infty,x_1)$, $I_+:=(x_N,\infty)$, $Q_{\pm}:=I_{\pm}\times (0,\tau)$.  
By arguing as in the last part of the proof of Proposition \ref{preco}-$(ii)$ (in particular, see the proof of \eqref{contra5CONT}-\eqref{contra6CONT}), it follows that for all
$l=1,\dots,N-1$ and $\zeta_l\in C^1_c(Q_l)$, $\zeta_l\geq 0$, 
\begin{equation}\label{contra5BIS}
\iint_{Q_l} \{|u_{1r}-u_{2r}|\,\zeta_{lt}+  |\varphi(u_{1r})-\varphi(u_{2r})(x,t) |\, \zeta_{lx}\}\,dxdt \geq 0  
\end{equation}
and, for all $\zeta^{\pm}\in C^1_c(Q_{\pm})$, $\zeta^{\pm}\ge0$,
\begin{equation}\label{contra6BIS}
\iint_{Q_\pm}  \{|u_{1r}-u_{2r}|\, \zeta_t^\pm + |\varphi(u_{1r})-\varphi(u_{2r})|\, \zeta_x^\pm(x,t)\}\, dxdt \geq 0
\end{equation}
(recall that $\varphi$ by assumption is increasing).
We must show  that \eqref{contra5BIS} and \eqref{contra6BIS} imply \eqref{contra}.
Let $h\in C^1_c(0, \tau_1)$, $h\ge0$, and
\begin{eqnarray*}%\label{ft.rholp}
\rho_{l,p}(x)=\!\!\!\!\!\!\!\!&&  p\left(x-x_l -\tfrac{1}{p}\right) \chi_{\left[x_l +1/p, x_l +2/p\right)}
+\chi_{\left[x_l +2/p, x_{l+1} -2/p\right)}(x)-\\
&&- \,p\left(x-x_{l+1}+\tfrac{1}{p}\right) \chi_{\left[x_{l+1} -2/p,x_{l+1} -1/p\right)}(x)
\qquad (l=1,\dots,N-1), \nonumber
\end{eqnarray*}
with $p\in\N$ sufficiently large. 
By standard approximation arguments we may choose  $\zeta_l=\zeta_{l,p}:=\rho_{l,p}(x)h(t)$ in \eqref{contra5BIS}: 
\begin{eqnarray}\label{contra8BIS}
&& 
0\le\iint_{Q_l} \big\{|u_{1r}-u_{2r}|\,\rho_{l,p}(x)h'(t)+ 
\left |\varphi(u_{1r})-\varphi(u_{2r})\right |\rho_{l,p}'(x)h(t)\big\}\,dxdt \,.
 \end{eqnarray}
By the Dominated Convergence Theorem, as $p\to \infty$,  
\begin{equation*}
\iint_{Q_l} \big\{|u_{1r}-u_{2r}|\,\rho_{l,p}(x)h'(t)\,dxdt \to  \int_0^{\tau} dt\, h'(t)\int_{x_l }^{x_{l+1} }|u_{1r}-u_{2r}|\,dx\,.
\end{equation*}
Since $\rho_{l,p}'(x)= p\chi_{\left(x_l +\frac 1p, x_l +2/p\right)}(x) - p\chi_{\left(x_{l+1} -2/p,x_{l+1} -1/p\right)}(x)$
and $\varphi$ is bounded, it follows from \eqref{crus} and the Dominated Convergence Theorem that
\begin{eqnarray*}
&&\limsup_{p\to \infty} \iint_{Q_l} |\varphi(u_{1r})-\varphi(u_{2r}) |\rho_{l,p}'(x)h(t)\,dxdt \le \\ 
&&\qquad\qquad \le \int_0^{\tau} h(t)\,\left(\lim_{p\to\infty} p\int_{x_l+\frac{1}{p}}^{x_l+\frac{2}{p}}|\varphi(u_{1r})-\varphi(u_{2r}) |\,dx \right)\,dt=0 \,.
\end{eqnarray*}
Hence, by \eqref{contra8BIS}, 
\begin{equation*}
 \int_0^{\tau} dt\, h'(t)\int_{x_l}^{x_{l+1} }|u_{1r}(x,t)-u_{2r}(x,t)|\,dx\ge0
\end{equation*} 
and, by a proper choice of $h$, 
\begin{equation}\label{eq.chiave.cont}
\|u_{1r}(\cdot,t) - u_{2r}(\cdot,t)\|_{L^1(I_l)} \le \|u_{1r}(\cdot,t_1) - u_{2r}(\cdot,t_1)\|_{L^1(I_l)} 
\quad\text{for every $0<t_1\le t\le\tau$} 
\end{equation}
(recall that $u_{ir} \in C((0,T];\mathcal{M}(\R))$, $i=1,2$, by Proposition \ref{preco}-$(ii)$).
Letting $t_1\to 0^+$ it follows from \eqref{epc} that $\|u_{1r}(\cdot,t) - u_{2r}(\cdot,t)\|_{L^1(I_l)} = 0$
for a.e. $t\in(0,\tau)$ and all $l=1,\dots,N-1$. The proof that
$\|u_{1r}(\cdot,t) - u_{2r}(\cdot,t)\|_{L^1(I_-\cup I_+)}=0$ 
for a.e.~$t\in(0,\tau)$ is similar, so we have proven
\eqref{contra}.

Next, let us prove that 
\begin{equation}\label{tesipro}
\text{ $u_1=u_2 \quad$ in $\mathcal{M}(\R\times (0,\tau))$\,.}
\end{equation}
By \eqref{ewf} and \eqref{contra}, for every $\zeta\in C^1([0,\tau];C^1_c(\R))$, $\zeta(\cdot,\tau)=0$ in $\R$
$$
\int_0^{T}\!\! \langle u_{1s}(\cdot,t) \! - \! u_{2s}(\cdot,t), \zeta_t(\cdot,t)\rangle_{\R}dt \!
= \! \! \iint_{S} \! \{(u_{1r} \! - \! u_{2r})\,\zeta_t+ [\varphi(u_{1r})\! - \! \varphi(u_{1r})] \zeta_x\}dxdt= 0\,.
$$
We argue as in the proof of Lemma \ref{lefdp}: there exists a null set $F_0\subset (0,\tau)$ such that 
$\langle u_{1s}(\cdot,t)-u_{2s}(\cdot,t),\rho\rangle_{\R}=0$ for all $t\in (0,\tau)\setminus F_0$ and
$\rho \in C^1_c(\R)$.
Hence $u_1=u_2$ in $L^{\infty}(0,\tau;\mathcal{M}(\R))$ and, by \eqref{contra}, equality \eqref{tesipro} follows.
\smallskip

If $\tau=T$ the proof is complete. Otherwise, there exist $N_1<N$ different points $x_{l_k} \in \{x_1,\dots,x_N\}$ such that $u_{is}(\cdot,\tau)(\{x_{l_k}\})>0$ for each $k=1,\dots,N_1$ and $i=1,2$; moreover, for every point $x_l \in \{x_1,\dots,x_N\}$, $x_l\neq x_{l_k}$ 
it follows from \eqref{dist1'}, with $x_0=x_1=x_l$, that $u_{1s}(\cdot,\tau)(\{x_l\})=u_{2s}(\cdot,\tau)(\{x_l\})=0$, since $\varphi(u_{1r})=\varphi(u_{2r})$ in $\R\times (0,\tau)$ by \eqref{contra}. 
Then we set
$$
\tau_1:= \sup\, \{t\in[\tau,T) \,|\,  u_{is}(\cdot,t)(\{x_{l_k}\})>0 
\;\; \forall  k=1,\dots,N_1;\;i=1,2\}\,.
$$ 
Arguing as in the proof of \eqref{contra} we obtain inequality \eqref{eq.chiave.cont} for every $\tau <t_1\leq t\leq \tau_1$. Since $u_{ir} \in C((0,T];\mathcal{M}(\R))$, $i=1,2$ (see Proposition \ref{preco}-$(ii)$), and $u_{1r}(\cdot,\tau)=u_{2r}(\cdot,\tau)$, letting $t_1\to \tau^+$ we get $u_{1r}=u_{2r}$ $a.e.$ in $\R\times (\tau,\tau_1)$ (whence, also $u_1=u_2 $ in $\mathcal{M}(\R\times (\tau,\tau_1)$) and the proof is completed 
in a finite number of steps.
\hfill$\square$

\smallskip

Let us finally prove Proposition \ref{wai1}.

\smallskip

\noindent {\em Proof of Proposition \ref{wai1}.}
A calculation proves that the solution defined by 
\eqref{sol1} if $p<0$, respectively by \eqref{sol2} 
if $0<p<1$
is an 
entropy solution of problem \eqref{ester}-\eqref{esterbis}. 
If $p<0$, the solution also satisfies \eqref{crus} for $0<t<1$ and \eqref{epc}, so claim $(i)$ 
follows from the uniqueness result in Theorem \ref{exiuni}. 
If $0<p<1$, uniqueness of entropy solutions such that $u_s(t)=0$ for $t>0$ and $u_r\in L^{\infty}(\R \times (\tau,T))$ for $\tau\in(0,T)$ can be used  
(the proof of this uniqueness result is very similar to that given in \cite{LP}, thus we omit the details; see also Remark \ref{uLP}). 
Hence claim $(ii)$ follows.
 \hfill$\square$

\begin{remark} It is instructive to describe the approximation procedure which gives the solutions mentioned in Proposition \ref{wai1}.  Consider the family of approximating problems
\begin{equation*}
\left\{\begin{array}{ll}
u_{nt}+ \left[\varphi(u_n)\right]_x=0 
& \quad\mbox{in $S$} 
\smallskip\\
u_n=\frac n2\,\chi_{\left(-\frac 1n,\frac 1n\right)}  &\quad\mbox{in}\  \R\times \{0\}\,,
\end{array}\right. \leqno{(R_n)}
\end{equation*}
with $\varphi$ given by \eqref{esterbis}. It is easily seen that the entropy solution of $(R_n)$ is 
\begin{equation*}%\label{p3}
u_n(x,t):=\left\{\begin{array}{ll}
0 & \quad\mbox{if \;$x\ge |p|t+\frac 1n$\,,} 
\smallskip\\
\left(\frac{n|p|t}{nx-1}\right)^\frac{1}{1-p}-1 &\quad\mbox{if \;$|p|t+\frac 1n> x \ge\left(\frac{2}{n+2}\right)^{1-p}|p|t+\frac 1n$\,,} \smallskip  \\
\frac n2  &\quad\mbox{if \;$\left(\frac{2}{n+2}\right)^{1-p}|p|t+\frac 1n> x\ge \frac{2 (\sgn p)}{n} \left[\left(\frac{n+2}{2}\right)^p-1\right] t - \frac 1n$\,,} \smallskip  \\
0  &\quad\mbox{if \;$\frac{2 (\sgn p)}{n} \left[\left(\frac{n+2}{2}\right)^p-1\right]t - \frac 1n> x$\,.}
\end{array}\right. 
\end{equation*}
for $0\le t\le t_n:=\frac{1}{\varphi\left(\frac n2\right)-\frac n2\varphi'\left(\frac n2\right)}\,.$ At $t=t_n$ a shock $x=\xi(t)$ stems from $x=x_n:=\frac 1n \, \frac{\varphi\left(\frac n2\right)+\frac n2\,\varphi'\left(\frac n2\right)}{\varphi\left(\frac n2\right)-\frac n2\,\varphi'\left(\frac n2\right)}$\,, which solves the problem
\begin{equation*}
\left\{\begin{array}{ll}
\xi_n'(t) =\dfrac{\varphi(u_n^{(1)}(\xi_n(t),t))}{u_n^{(1)}(\xi_n(t),t)}=\sgn p\;\dfrac{\left(\frac{n|p|t}{n\,\xi-1}\right)^{\frac{p}{1-p}}-1}{\left(\frac{n|p|t}{n\,\xi-1}\right)^{\frac{1}{1-p}}-1}
& \quad\mbox{if $t>t_n$\,,} 
\smallskip\\
\xi_n(t_n)=x_n\,.  &\
\end{array}\right. 
\end{equation*}
Hence for $t>t_n$ the entropy solution of $(R_n)$ is
\begin{equation*}%\label{p33}
u_n(x,t):=\left\{\begin{array}{ll}
0 & \quad\mbox{if \;$x\ge |p|t+\frac 1n$\,,} 
\smallskip\\
\left(\frac{n|p|t}{nx-1}\right)^\frac{1}{1-p}-1 &\quad\mbox{if \;$|p|t+\frac 1n> x\ge \xi_n(t)$\,,} \smallskip  \\
0  &\quad\mbox{if \;$\xi_n(t)\ge x$\,.}
\end{array}\right. 
\end{equation*}
Letting $n\to\infty$ 
we obtain the entropy solution defined in parts $(i)$ (if $p<0$) and $(ii)$ (if $0<p<1$) of Proposition \ref{wai1}. 
\end{remark}
%%%%%%%%%%%%%%%%%%%%%%%%%%%%%%%%%%%%%%%%%%%%%%%%%%%%%%%%%%%%%

%\section*{Ethical Statement}
%Conflict of Interest: The authors declare that they have no conflict of interest.

%%%%%%%%%%%%%%%%%%%%%%%%%%%%%%%%%%%%%%%%%%%%%%%%%%%%%%%%%%%%%%%%%%%%%%%%%%%%%%%%%%%%%%%%%%%%%%%%%%%%%%%%%%%%%%%%%%

%%%%%%%%%%%%%%%%%%

\end{document}